\documentclass[12pt,a4pape,reqno]{amsart}
\usepackage[top=25mm, bottom=25mm, left=30mm, right=30mm]{geometry}
\usepackage{mathptmx}
\usepackage{mathrsfs}
\usepackage{amsfonts,amsmath,amscd,amssymb,paralist,amsthm}
\usepackage{amscd}
\usepackage{verbatim}
\usepackage{url}
\usepackage[all]{xy}
\usepackage{color}
\usepackage[colorlinks=true,citecolor=blue, backref]{hyperref}
\usepackage{tikz}
\usepackage{tikz-cd}
\usetikzlibrary{arrows}
\usetikzlibrary{patterns}
\usetikzlibrary{lindenmayersystems,backgrounds}

\theoremstyle{plain}

\newtheorem{thm}{Theorem}[section]

\newtheorem{lem}[thm]{Lemma}
\newtheorem{prop}[thm]{Proposition}
\newtheorem{cor}[thm]{Corollary}
\newtheorem{example}[thm]{Example}

\theoremstyle{definition}
\newtheorem{defn}[thm]{Definition}
\newtheorem{rem}[thm]{Remark}

\newtheorem{conj}{Conjecture}

\numberwithin{equation}{section}

\newcommand{\vertiii}[1]{{\left\vert\kern-0.25ex\left\vert\kern-0.25ex\left\vert #1
    \right\vert\kern-0.25ex\right\vert\kern-0.25ex\right\vert}}

\newcommand{\Z}{{\mathbb{Z}}}
\newcommand{\N}{\mathbb{N}}

\newcommand{\ep}{\varepsilon}

\newcommand{\T}{\mathbb{T}}
\newcommand{\R}{\mathbb{R}}

\def \RP {{\bf RP}}

\def\R {\mathbb R}
\def\Z {\mathbb Z}
\def\Q {\mathbb Q}
\def\T {\mathbb T}
\def\O {\mathcal O}

\def\A {\mathcal A}

\def \a {\alpha }

\def \ep {\epsilon}
\def \d {\delta}
\def \D {\Delta}

\begin{document}

\title[Structure theorems of commuting transformations and minimal $\mathbb{R}$-flows]{Structure theorems of commuting transformations and minimal $\mathbb{R}$-flows}
\author[S.~Shao]{Song Shao}
\address[S. Shao]{School of Mathematical Sciences, University of Science and Technology of China, Hefei, Anhui 230026, China}
\email{songshao@ustc.edu.cn}

\author[H.~Xu]{Hui Xu}

\address[H. Xu]{Department of Mathematics, Shanghai Normal University, Shanghai, 200234, China}
\email{huixu@shnu.edu.cn}

\thanks{This research is supported by National Key R$\&$D Program of China (No. 2024YFA1013601, 2024YFA1013600), and National Natural Science Foundation of China (12426201, 12371196, 12201599). }

\maketitle

%%%%%%%%%%%%%%%%%%%%%%%%%%%%%%%%%%%%%%%%%%%%%%%%%%%%%%%%%%%%%%%%%%%%%%%%%%%%%%%%%%%%%%%%%%%%%%%%%%%%%%%%%%%%%%%%%%%%%%%%%%%%%%%%%%%%%%%%%%%%%%%%%%%%%%%%%%%%%%%%%%%%%
\begin{abstract}
In this paper, we develop several structure theorems concerning commuting transformations and minimal $\mathbb{R}$-flows.
Specifically, we show that if $(X,S)$,  $(X,T)$ are minimal systems with $S$ and $T$ being commutative, then they possess an identical higher-order regionally proximal relation. Consequently, both $(X, S)$ and $(X, T)$ share the same increasing sequence of pro-nilfactors. For minimal $\mathbb{R}$-flows, we introduce the concept of higher-order regionally proximal relations and nilfactors, and establish that nilfactors are characteristic factors for minimal $\mathbb{R}$-flows, up to almost one to one extensions.

\end{abstract}

\section{Introduction}
Inspired by Furstenberg's groundbreaking ergodic-theoretic proof of Szem\'{e}redi's theorem, numerous mathematicians have dedicated significant effort to investigating the limiting behavior of the averages
\begin{equation}\label{eq1}
\frac{1}{N}\sum_{n=0}^{N-1}\prod_{j=1}^{k}f_j(T^{jn}).
\end{equation}
The $L^2$ convergence of (\ref{eq1}) was finally established by Host-Kra in \cite{HK05} (later by Ziegler \cite{Zie05} independently). Both approaches rely on the use of characteristic factors to capture the limit of (\ref{eq1}). Specifically, for an ergodic system $(X,\mathcal{X},\mu, T)$, there is a sequence of nil-factors $X\rightarrow \cdots\rightarrow Z_{k-1}\rightarrow \cdots\rightarrow Z_{0}$ such that the $L^{2}$ limit of (\ref{eq1}) can be reduced to the limit in $Z_{k-1}$.

A counterpart of a characteristic factor in topological dynamical systems was first introduced in \cite{G94}. To get the corresponding factors in a minimal system, Host, Kra and Maass introduced the notion of regionally proximal relation of higher order in \cite{HKM}, denoted by ${\rm \bf RP}^{[d]}$. They further showed that ${\rm \bf RP}^{[d]}$ is a closed invariant equivalence relation for a minimal distal system. Later, Shao and Ye showed in \cite{SY12} that ${\rm \bf RP}^{[d]}$ is a closed invariant equivalence relation for any minimal system under an abelian group action. Recently, Glasner, Huang, Shao, Weiss and Ye showed in \cite{GHSWY} that $X_{d}=X/{\rm \bf RP}^{[d]}$ are characteristic factors for minimal systems, up to almost one to one extensions.

It is well established that the study of dynamical systems initially emerged from the consideration of continuous-time evolutions. Nowadays, most research in  dynamics restricts on discrete time due to the fact many results in continuous time can be derived from the discrete case. In \cite{BLM12}, Bergelson, Leibman and Moreira gave methods to derive continuous-time versions of various discrete-time ergodic theorems. In \cite{P11}, Potts provided the natural analogues of results for multiple ergodic averages along flows. Specifically, to establish the convergence of continuous-time ergodic averages of a product of functions evaluated at return times along polynomials, she demonstrated that the Host-Kra factors are sill characteristic for these averages. As a result, the problem can be reduced to the case where the system is an inverse limit of nilflows. For further results on the convergence of continuous-time polynomial multiple ergodic averages, we refer the reader to \cite{A12, BLM12, P11, Zie05, Zie1}.
In this paper, we are focusing on the continuous-time regionally proximal relation of higher order and its associated characteristic factors and their relations with discrete-time versions.
%In this paper, we focus on the natural analogues of some of these results for flows.

A continuous flow or $\mathbb{R}$-flow is a triple $(X, \mathbb{R}, \phi)$ where $X$ is a compact metric space and $\phi: \mathbb{R}\times X\rightarrow X$ is a continuous mapping such that $\phi(0,x)=x$ for all $x\in X$ and $\phi(s,\phi(t,x))=\phi(s+t,x)$ for all $s,t\in\mathbb{R}$ and $x\in X$. When there is no ambiguity, we denote the $\mathbb{R}$-flow by $(X,\mathbb{R})$ without explicitly mentioning the map $\phi$. For each $t\in\mathbb{R}$, the map $T^{t}:=\phi(t, \cdot)$ induces a discrete-time dynamical system.  We also use $(X, T^{t})$ to denote this discrete-time system (is also called a system) and use $(X,\{T^{t}\}_{t\in\R})$ to denote the original continuous flow. There are inheritances of many dynamical properties, such that as recurrence, transitivity, mixing etc, between $(X,\{T^t\}_{t\in \R})$ and $(X, T^{t})$, see \cite[Chapter II]{Vries93}.

A classical result due to Pugh and Shub asserts that the dynamical system $(X,\{T^t\}_{t\in \R})$ is ergodic with respect to some measure $\mu$ if and only if $(X,T^{t})$ is ergodic with respect to $\mu$ for all $t\in\mathbb{R}$, with the possible exception of a countable set of $t$-values \cite{PS71}. In the topological setting, an analogous statement holds: the system $(X,\{T^t\}_{t\in \R})$ is minimal if and only if $(X,T^{t})$ is minimal for all $t\in\mathbb{R}$ except possibly for a countable set (see, for example, \cite[4.24]{G03}). In this paper, we present a new proof of this fact through reducing to consider the minimality on the maximal equicontinuous factors.

\medskip

In \cite{FK05} Frantzikinakis and Kra studied the convergence of multiple ergodic averages for commuting transformations. They showed that if $T$ and $S$ are commuting measure-preserving transformations of a probability space $(X,{\mathcal X},\mu)$ and that both $T$ and $S$ are ergodic, then these transformations share identical Host-Kra seminorms and the same characteristic factors for multiple ergodic averages \cite[Proposition 3.1]{FK05}. In particular, for ergodic commuting
transformations the limiting behaviors of the corresponding linear multiple ergodic averages
are the same \cite[Theorem 1.2]{FK05}. We establish the corresponding topological results.
We show the following result (for the case of $d=1$, this result was given by Akin and Glasner \cite{AG98}).

\begin{thm}\label{thm-RPd}
Let $G, H$  be abelian groups. Let $(X,G)$ and $(X,H)$ be commutative actions and both transitive and $d\in \N$. Then $\RP^{[d]}(X,G)=\RP^{[d]}(X,H)$.
\end{thm}

Thus if $(X,S)$,  $(X,T)$ are minimal systems with $S$ and $T$ being commutative homeomorphisms, then for any $d\in\mathbb{N}$,  we have that ${\RP}^{[d]}(X, S)={\RP}^{[d]}(X, T)$. Thus both $(X, S)$ and $(X, T)$ have the same increasing sequence of pro-nilfactors:
$$\{pt\}= X_0  \longleftarrow  X_1 \longleftarrow \cdots  \longleftarrow X_n  \longleftarrow \cdots  \longleftarrow X,$$
where $X_d=X/\RP^{[d]}(X)$ is the the
maximal $d$-step pro-nilfactor, $d\in \N$.

Moreover, we demonstrate that two commutative minimal homeomorphisms possess the identical dynamical cube ${\bf Q}^{[d]}$ as defined by Host-Kra-Maass in \cite{HKM} and the structure $N_d$ as introduced by Glasner in \cite{G94} (for precise definitions, refer to Subsection 2.1).

\begin{thm}\label{cube and Nd}
 \begin{enumerate}
\item[(1)] Let $G$ and $H$ be abelian groups.  Let $(X,G)$ and $(X,H)$ be commutative actions and both minimal. Then  we have ${\bf Q}^{[d]}(X, G)={\bf Q}^{[d]}(X, H)$ for any $d\in\mathbb{N}$.
\item[(2)] Let $(X,S)$,  $(X,T)$ be  minimal systems. If $S$ and $T$ commutes, then $N_{d}(X,S)=N_{d}(X,T)$ for each $d\in\mathbb{N}$.
\end{enumerate}
 \end{thm}

A direct consequence of Theorem \ref{cube and Nd} is Theorem C in \cite{GHSWY}, which states that for a minimal system $(X,T)$ and $k\in\mathbb{N}$,  $(X, T^{k})$ is minimal if and only if $N_{d}(X, T)=N_{d}(X, T^{k})$ for each $d\in\mathbb{N}$.

\medskip

Let $(X,\{T^{t}\}_{t\in \R})$ be a minimal $\R$-flow. Then $(X,T^t)$ ia also minimal for all $t\in \R\setminus E$, where $E$ is a countable set. By Theorem \ref{thm-RPd}, ${\bf RP}^{[d]}(\R)={\bf RP}^{[d]}(T^{t})$, and the space $X/{\bf RP}^{[d]}(\mathbb{R})$ is the same to $X/{\bf RP}^{[d]}(T^{t})$ for all $t\in\mathbb{R}\setminus E$. The discrete system $(X/{\bf RP}^{[d]}(T^{t}), T^{t})$ is a pro-nilsystem for any $t\in\R\setminus E$. Moreover, it can be shown that an analogous result holds for $\mathbb{R}$-flows. Specifically, $(X/{\bf RP}^{[d]}(\mathbb{R}), \{T^t\}_{t\in \R})$ is isomorphic to a $d$-step pro-nilflow (see Subsection \ref{subsection-nilflow} for details).

In the paper, we will define characteristic factors for minimal $\mathbb{R}$-flows, and show that $X_{d}=X/{\rm \bf RP}^{[d]}(\R)$ are characteristic factors for minimal $\mathbb{R}$-flows, up to almost one to one extensions.
To be precise, we present the following result, with the relevant concepts detailed in Section \ref{section-TCF}.  An equivalent form is given in Theorem \ref{equiv-thm-Real}.

\begin{thm}\label{thm-Real}
Let $(X,\{T^{t}\}_{t\in\R})$ be a minimal $\R$-flow and $d\in \N$. Let $\pi:X\rightarrow X_d$ be the factor map from $X$ to its maximal $d$-step pro-nilfactor $X_d$. Then there are minimal $\R$-flows $X^*$ and $X_d^*$ which are almost one to one
extensions of $X$ and $X_d$ respectively, and a commuting diagram below such that $\pi^*: X^*\rightarrow X^*_d$ is open, and $X_d^*$ is a topological characteristic factor of $X^*$ with respect to $\{\a_1t, \a_2t,\ldots, \a_dt\}$, where $\a_1,\ldots, \a_d$ are distinct nonzero real numbers.
\[
\begin{CD}
X @<{\varsigma^*}<< X^*\\
@VV{\pi}V      @VV{\pi^*}V\\
X_d @<{\varsigma}<< X_d^*
\end{CD}
\]

If in addition $\pi$ is open, then $X^*=X$, $X_d^*=X_d$ and $\pi^*=\pi$.
\end{thm}

On the one hand, one can obtain the discrete-time system from the continuous-time one. On the other hand, there is a natural way to get a continuous-time system from the discrete-time one by doing suspensions. For the nilfactors of the suspension flows, we have the following result.

\begin{thm}\label{mod RP vs suspension}
Let $(X,T)$ be a minimal system and $(\tilde{X},\{T^t\}_{t\in \R})$ be its suspension flow. Then for any $d\in\mathbb{N}$, $\tilde{X}_{d}:=\tilde{X}/{\bf RP}^{[d]}(\tilde{X},\mathbb{R})=\widetilde{X_{d}}$, where $\widetilde{X_{d}}$ is the suspension of $X_{d}=X/{\bf RP}^{[d]}(X,T)$.
\end{thm}

%A {\em $d$-step nilflow} is a system $(N/\Gamma,\mathbb{R})$, where $N$ is a $d$-step nilpotent Lie group, $\Gamma$ is a uniform lattice of $N$, and the action of $\mathbb{R}$ on $N/\Gamma$ is determined by an one-parameter subgroup $\{a(t):t\in\mathbb{R}\}$. A {\em $d$-step pro-nilflow} is an inverse limit of $d$-step nilsystems.

%\begin{thm}\label{RPd}
%Let $G$  be an abelian group and $(X,G)$ be a minimal system. If $H$ is a subgroup of $G$ such that $(X,H)$ is transitive. Then for any $d\in\mathbb{N}$,  we have ${\rm\bf RP}^{[d]}(X, G)={\rm\bf RP}^{[d]}(X, H)$.
%\end{thm}

\subsection*{Organization of the paper} In section \ref{section-pre} we give some basic notions and results used in the sequel. In section \ref{section-minimal}, we study the regionally proximal relation for commutative minimal actions and consequently  we show the fact that $(X,\mathbb{R})$ is minimal if and only if $(X,T^{t})$ is minimal  for all $t\in\mathbb{R}$ except a countable set.  In section \ref{section-RPd}, we discuss the regionally proximal relation of higher order for commutative minimal actions. In section \ref{section-TCF},  we discuss the topological characteristic factors of continuous flows. In section \ref{section-suspension}, we discuss the suspension flows.  In the appendix, we show the equivalent characterizations of regionally proximal relation of higher order.

\subsection*{Acknowledgments} We would like to thank Professor Wen Huang, Professor Xiangdong Ye  and Jiahao Qiu for their very useful comments.

\section{Preliminaries}\label{section-pre}

In this section we recall some basic definitions and results used in the paper. In the article, the sets of real numbers, rational numbers, integers, and positive integers are denoted by $\R$, $\Q$, $\Z$ and $\N$ respectively.

\subsection{Topological dynamical system}

By a {\em topological dynamical system} or a {\em flow} we mean a triple $(X, G, \phi)$ where $X$ is compact Hausdorff space, $G$ is a topological group with the unit $e$ and $\phi: G\times X\rightarrow X$ is a continuous mapping such that $\phi(e, x)=x$ for all $x\in X$ and $\phi(g, \phi(h,x))=\phi(gh, x)$ for all $g,h\in G$ and $x\in X$. When there is no confusion, we use $(X,G)$ short for $(X, G, \phi)$ and say it is a {\it $G$-system}. In addition, we use $gx$ short for $\phi(g,x)$. In this paper,
we always assume that $X$ is a compact metric space with the metric $\rho(\cdot, \cdot)$, and we assume that all groups $G$ involved are abelian.

When $G=\mathbb{Z}$, a $G$-system is determined by a single homeomorphism $T:X\rightarrow X$. In this case, we usually denote $(X,\Z)$ by $(X,T)$, and also call it a {\em discrete system}. When $G=\mathbb{R}$, we usually denote $(X,\R)$ by $(X,\{T^t\}_{t\in \R})$, where $\phi: \R\times X\rightarrow X$, $T^t:=\phi(t,\cdot)$, and also call it a {\em continuous flow} or an $\mathbb{R}$-{\it flow}.

\medskip

Let $(X,G)$ be a $G$-system and $x\in X$. Let $\O(x,G)=\{gx: g\in G\}$ be the
{\em orbit} of $x$, which is also denoted by $G x$. We usually denote the closure of $\O(x,G)$ by $\overline{\O}(x,G)$, or $\overline{Gx}$. A subset
$A\subseteq X$ is called {\em invariant} if $g a\in A$ for all
$a\in A$ and $g\in G$. When $Y\subseteq X$ is a closed and invariant subset of the flow $(X, G)$ we say that the system $(Y, G)$ is a {\em subsystem} of $(X, G)$. If $(X, G)$ and $(Y, G)$ are two $G$-systems, their {\em product system} is the $G$-system $(X \times Y, G)$, where $g(x, y) = (gx, gy)$ for any $g\in G$ and $x\in X, y\in Y$. For $n \geq 2$ we write $(X^n,G)$ for the $n$-fold product flow $(X\times
\cdots \times X, G)$.

\medskip

A $G$-system $(X,G)$ is called {\em minimal} if $X$ contains no proper non-empty closed invariant subsets. A point $x\in X$ is called a {\em minimal point} or an {\em almost periodic point} if $(\overline{\O}(x,G), G)$ is a minimal system. A $G$-system $(X,G)$ is called {\em transitive} if every invariant non-empty open subset of $X$ is dense; and it is {\em point transitive} if there is a point with a dense orbit (such a point is called a {\em transitive point}). It is easy to verify that a system is minimal if and only if every orbit is dense.
The $G$-system $(X, G)$ is {\em weakly mixing} if the product system $(X \times X, G)$ is transitive.

A subset $A\subset G$ is {\it syndetic} if there is a compact set $K\subset G$ such that $G=KA$. It is well-known that for a $G$-system $(X,G)$, a point $x$ is a minimal point if and only if for any neighborhood $U$ of $x$, the set $\{g\in G: gx\in U\}$ is syndetic (e.g. see \cite[Chapter 1, Theorem 7]{Aus88}).

A {\it factor map} $\pi: X\rightarrow Y$ between the $G$-systems $(X, G)$ and $(Y, G)$ is a continuous surjective map which intertwines the actions; we say that $(Y, G)$ is a {\it factor} of $(X, G)$ and
that $(X, G)$ is an {\it extension} of $(Y, G)$. The systems are said to be {\it isomorphic} if $\pi$ is bijective.

\medskip

Let $(X,T)$ be a $G$-system. Fix $(x,y)\in X^2$. It is a {\em proximal} pair if $\inf_{t\in G} \rho(tx, ty)=0$; it is a {\em distal} pair if it is not proximal. Denote by $P(X,G)$ or $P(X)$ the set of proximal pairs of $(X,G)$. $P(X,G)$ is also called the {\em proximal relation} of $(X, G)$. A $G$-system $(X, G)$ is {\it distal} if $P(X, G)= \D_X$, where $\D_X=\{(x,x)\in X^2: x\in X\}$ is the diagonal of $X\times X$. A $G$-system $(X, G)$ is {\em equicontinuous} if for any $\ep>0$, there is a $\d>0$ such that whenever $x,y\in X$ with  $\rho(x,y)<\d$, then $\rho(tx,ty)<\ep$ for all $t\in G$.
Any equicontinuous system is distal.

Let $(X,G)$ be a $G$-system. There is a smallest invariant equivalence relation $S_{eq}(X)$ such that the quotient flow $(X/S_{eq}(X), G)$ is equicontinuous \cite[Theorem 1]{EG60}. The equivalence relation $S_{eq}(X)$ is called the {\em equicontinuous structure relation} and the factor $(X_{eq}=X/S_{eq}(X), G)$ is called the {\em maximal equicontinuous factor} of $(X, G)$.

\subsection{Nilmanifolds and nilflows}\label{subsection-nilsystem}

Let $N$ be a group. For $g, h\in N$ and $A,B \subseteq N$, we write $[g, h] =
ghg^{-1}h^{-1}$ for the commutator of $g$ and $h$ and
$[A,B]$ for the subgroup spanned by $\{[a, b] : a \in A, b\in B\}$.
The commutator subgroups $N_j$, $j\ge 1$, are defined inductively by
setting $N_1 = N$ and $N_{j+1} = [N_j ,N]$. Let $d \ge 1$ be an
integer. We say that $N$ is {\em $d$-step nilpotent} if $N_{d+1}$ is
the trivial subgroup.

\medskip

Let $d\in\N$, $N$ be a $d$-step nilpotent Lie group and $\Gamma$ be a discrete
cocompact subgroup of $N$. The compact manifold $X = N/\Gamma$ is
called a {\em $d$-step nilmanifold}. The group $N$ acts on $X$ by
left translations and we write this action as $(g, x)\mapsto gx$.
The Haar measure $\mu$ of $X$ is the unique Borel probability measure on
$X$ invariant under this action. Fix $a\in N$ and let $T$ be the
transformation $x\mapsto a x$ of $X$, i.e. $a(g\Gamma)=(ag)\Gamma$ for each $g\in N$. Then $(X, \mu, T)$ is
called a {\em $d$-step nilsystem}. In the topological setting we omit the measure
and just say that $(X,T)$ is a $d$-step nilsystem. For more details on nilsystems, refer to \cite[Chapter 11]{HK18}.

%Here are some basic properties of nilsystems.

%\begin{thm}\cite{Leibman051, P}\label{thm-ParryLeibman}
%Let $(X = G/\Gamma,\mu , T )$ be a $d$-step nilsystem with $T$ the translation by the element $t\in G$. Then:

%\begin{enumerate}
%\item $(X, T )$ is uniquely ergodic if and only if $(X,\mu , T )$ is ergodic if and only if $(X, T )$ is minimal if and only if $(X, T )$ is transitive.

%\item Let $Y$ be the closed orbit of some point $x\in X$. Then $Y$ can be given the structure of a nilmanifold, i.e. $Y = H/\Lambda$, where $H$ is a closed subgroup of $G$ containing $t$ and $\Lambda$ is a closed cocompact subgroup of $H$.
%\end{enumerate}
%\end{thm}

%One can generalize the above results to the action of several translations. For example, let $X= G/\Gamma$ be a nilmanifold with the Haar measure $\mu$ and let $t_1,\cdots , t_k$ be commuting elements of $G$. If the group spanned by the translations $t_1, \cdots , t_k$ acts ergodically on $(X,\mu)$, then $X$ is uniquely ergodic for this group. For more details, refer to \cite{Leibman051, Leibman052, HK18}.

If $\{a_t\}_{t\in \R}$ is a one-parameter subgroup of $N$ then $\{a_t\}_{t\in \R}$ induces
a flow $\{T^t\}_{t\in \R}$ on $X$ defined by $T^t(g\Gamma)=(a_t\cdot g)\Gamma$ for all $g\in N$ and for all $t\in \R$. A flow defined in this manner is called a {\em $d$-step nilflow}. Generally, for A topological group $G$, a {\em $d$-step $G$-nilflow} is defined by $T_g(h\Gamma)=(a_gh)\Gamma$, where $G\rightarrow N, g\mapsto a_g$ is a homomorphism.

\medskip

We will need to use inverse limits of nilsystems, so we recall the
definition of a sequential inverse limit of systems. If
$(X_i,T_i)_{i\in \N}$ are systems %with $diam(X_i)\le 1$
and $\pi_i: X_{i+1}\rightarrow X_i$ is a factor map for all $i\in\N$, the {\em inverse
limit} of these systems is defined to be the compact subset of
$\prod_{i\in \N}X_i$ given by $\{ (x_i)_{i\in \N }: \pi_i(x_{i+1}) =
x_i, i\in\N\}$, and we denote it by
$\lim\limits_{\longleftarrow}(X_i,T_i)_{i\in\N}$.
It is a compact metric space.
% endowed with the distance $\rho((x_{i})_{i\in\N}, (y_{i})_{i\in \N}) = \sum_{i\in \N} 1/2^i \rho_i(x_i, y_i )$, where $\rho_{i}$ is the metric in $X_{i}$.
We note that the maps $\{T_i\}_{i\in \N}$ induce naturally a transformation $T$ on the inverse
limit such that $T(x_1,x_2,\ldots)=(T_1(x_1),T_2(x_2),\ldots)$. Similarly, we can define the inverse limit of $G$-systems.

%\medskip

%The following structure theorem characterizes inverse limits of nilsystems using dynamical parallelepipeds.

%\begin{thm}[Host-Kra-Maass]\cite[Theorem 1.2]{HKM}\label{HKM}
%Assume that $(X, T)$ is a transitive topological dynamical system and let $d \ge 2$ be an integer. The following properties are equivalent:
%\begin{enumerate}
%  \item If ${\bf x}, {\bf y} \in \Q^{[d]}$ have $2^d-1$ coordinates in common, then ${\bf x} = {\bf y}$.
%  \item If $x, y \in X$ are such that $(x, y,\ldots , y) \in  \Q^{[d]}$, then $x = y$.
%  \item $X$ is an inverse limit of $(d-1)$-step minimal nilsystems.
%\end{enumerate}
%\end{thm}

\begin{defn}% [{\cite[Definition 1.2]{HKM}}] %\label{HKM}[Host-Kra-Maass]
For $d\in \N$, a minimal system is called  a {\em $d$-step pro-nilsystem} or {\it system of order $d$} if $X$ is an inverse limit of $d$-step minimal nilsystems. A {\em $d$-step pro-nilflow} is an inverse limit of $d$-step nilflows.
\end{defn}

\subsection{Regionally proximal relation of order $d$}

Let $(X,G)$ be a $G$-system. A pair of points $(x,y)$ in $X$ is {\em regionally proximal} if  there are nets $\{x_n\}, \{y_n\}$ in $X$ and a net $\{g_n\}$ in $G$ such that $x_n\rightarrow x, y_n\rightarrow y$ and $\rho(g_nx_n, g_ny_n)\rightarrow 0$. Denote by ${\bf RP}(X,G)$ or $\RP(X)$ the set of regionally proximal pairs. $\RP(X,G)$ is closed and invariant. Recall that $G$ is abelian. It is well known that if $(X,G)$ is minimal, then the regionally proximal relation $\RP(X,G)$ is an equivalence relation and $S_{eq}(X)=\RP(X,G)$ ( see, for example, \cite[Chapter 9, Theorem 8]{Aus88}). Note that $\RP(X,G)=\D_X$ if and only if $(X,G)$ is equicontinuous (for example, see \cite[Chapter 7, Proposition 2]{Aus88}).

\begin{defn}
Let $(X,G)$ be a $G$-system and $d\in \N$. A pair $(x,y)\in X\times X$ is said to be {\it regionally proximal of order $d$} if for any $\delta>0$, there are $x',y'\in X$ with $\rho(x,x')< \delta, \rho(y,y')<\delta$ and ${\bf g}=(g^{(1)}, \ldots,g^{(d)})\in G^{d}$ such that for any ${\bf \epsilon}=(\epsilon_1,\ldots,\epsilon_{d})\in \{0,1\}^{d}\setminus\{{\bf 0}\}$,
\[ \rho \left({\bf  g}^{({\bf \epsilon})} x', {\bf  g}^{({\bf \epsilon})}y'\right)<\delta,\]
where ${\bf  g}^{({\bf \epsilon})} =\prod_{\epsilon_{i}=1}g^{(i)}$ and ${\bf 0}=(0,0,\cdots,0)$. Denote by ${\rm\bf RP}^{[d]}(X,G)$ the set of regionally proximal pairs of oder $d$. When $G=\mathbb{Z}$, we also use ${\rm \bf RP}^{[d]}(X,T)$ to denote ${\rm \bf RP}^{[d]}(X,\Z)$.
\end{defn}

It is easy to verify that
\[ {\rm\bf RP}^{[d]}(X, G)=\bigcap_{\alpha\in \mathcal{U}}\overline{ \bigcup_{{\bf g}\in G^{d}} \bigcap_{\epsilon\in\{0,1\}^{d}\setminus\{{\bf 0}\}}({\bf g}^{(\epsilon)} )^{-1}\alpha},\]
where $\mathcal{U}$ is the uniform structure on $X$.

The following theorems were proved in \cite{HKM} (for minimal distal systems) and
in \cite{SY12} (for general minimal systems).
%tell us conditions under which $(x,y)$ belongs to $\RP^{[d]}$ and the relation between $\RP^{[d]}$ and $d$-step pro-nilsystems, which are defined in Theorem \ref{HKM}.

\begin{thm}\label{thm-1}
Let $(X, T)$ be a minimal system. and let $d\in \N$. Then
\begin{enumerate}
%\item $(x,y)\in \RP^{[d]}$ if and only if $(x,y,\ldots,y)\in \Q^{[d+1]}$if and only if $(x,y,\ldots,y) \in \overline{\F^{[d+1]}}(x^{[d+1]})$.

\item $\RP^{[d]}(X,T)$ is an equivalence relation;

\item $(X,T)$ is a $d$-step pro-nilsystem if and only if $\RP^{[d]}(X,T)=\Delta_X$.
\end{enumerate}
\end{thm}

The regionally proximal relation of order $d$ allows us to construct the maximal $d$-step
pro-nilfactor of a system.
%That is, any factor of order $d$ (inverse limit of $d$-step minimal nilsystems) factorize through this system.

\begin{thm}\label{thm0}%\cite{SY}
Let $\pi: (X,T)\rightarrow (Y,T)$ be a factor map between minimal systems
and let $d\in \N$. Then
\begin{enumerate}
  \item $\pi\times \pi (\RP^{[d]}(X,T))=\RP^{[d]}(Y,T)$;
  \item $(Y,T)$ is a $d$-step pro-nilsystem if and only if $\RP^{[d]}(X,T)\subseteq R_\pi=\{(x_1,x_2)\in X^2: \pi(x_1)=\pi(x_2)\}$.
\end{enumerate}
In particular, $X_d:=X/\RP^{[d]}(X,T)$, the quotient of $(X,T)$ under $\RP^{[d]}(X,T)$, is the
maximal $d$-step pro-nilfactor of $X$. %(i.e. the maximal factor of order $d$).
\end{thm}

By Theorem \ref{thm-1} for any minimal system $(X,T)$,
$$\RP^{[\infty]}(X,T)=\bigcap\limits_{d\ge 1} \RP^{[d]}(X,T)$$
is a closed invariant equivalence relation (we write $\RP^{[\infty]}(X,T)$ in case of ambiguity). Now we formulate the
definition of $\infty$-step pro-nilsystems. % or systems of order $\infty$.

\begin{defn}[{\cite[Definition 3.4]{D-Y}}]
A minimal system $(X, T)$ is an {\em $\infty$-step
pro-nilsystem} or {\em a system of order $\infty$}, if the equivalence
relation $\RP^{[\infty]}(X,T)$ is trivial, i.e., coincides with the
diagonal.
\end{defn}

Similar to Theorem \ref{thm0}, one can show that the quotient of a
minimal system $(X,T)$ under $\RP^{[\infty]}(X)$ is the maximal
$\infty$-step pro-nilfactor of $(X,T)$ (for example, see \cite[Remark 3.4]{D-Y}). We denote the maximal
$\infty$-step pro-nilfactor of $(X,T)$ by $X_\infty=X/{\rm\bf RP}^{[\infty]}(X,T)$.

\subsection{Topological characteristic factors}\label{subsection-chara}

Let $X, Y$ be sets, and let $\phi : X\rightarrow Y$ be a map. A subset $L$ of $X$ is called
{\em $\phi$-saturated} if $$\{x\in L: \phi^{-1}(\phi(x))\subseteq L\}=L,$$ i.e., $L=\phi^{-1}(\phi(L))$.

Given a factor map $\pi: (X,T)\rightarrow (Y,T)$ and $d\ge 2$,
the system $(Y,T)$ is said to be a {\em $d$-step topological
characteristic factor
(for $\tau_d=T\times T^2\times \cdots \times T^d$)
of  $(X,T)$}, if there exists a dense $G_\d$
subset $\Omega$ of $X$ such that for each $x\in \Omega$ the orbit
closure $L_x=\overline{\O}((x, \ldots,x), \tau_d)$ is $\pi^{(d)}=\pi\times \cdots \times
\pi$ ($d$-times) saturated. That is, $(x_1,x_2,\ldots, x_d)\in L_x$
if and only if $(x_1',x_2',\ldots, x_d')\in L_x$, where $\pi(x_i)=\pi(x_i')$ for all $i\in \{1,2,\ldots, d\}$.
The following theorem was proved in \cite{GHSWY}.

\begin{thm}[{\cite[Theorem A]{GHSWY}}]\label{thm-GHSWY}
Let $(X,T)$ be a minimal system, and $\pi:X\rightarrow X_\infty$ be the factor map from $X$ to its maximal $\infty$-step pro-nilfactor $X_\infty$. Then there are minimal system $X^*$ and $X_\infty^*$ which are almost one to one
extensions of $X$ and $X_\infty$ respectively, and a commuting diagram below such that $\pi^*: X^*\rightarrow X^*_\infty$ is open, and $X_\infty^*$ is a $d$-step topological characteristic factor of $X^*$ for all $d\ge 2$.
\[
\begin{CD}
X @<{\varsigma^*}<< X^*\\
@VV{\pi}V      @VV{\pi^*}V\\
X_\infty @<{\varsigma}<< X_\infty^*
\end{CD}
\]
\end{thm}

\subsection{Some facts about hyperspaces} \label{sub:ellis}

Let $X$ be a compact metric space. Let $2^X$ be the collection of non-empty closed subsets of $X$.
Let $\rho$ be the metric on $X$.
One may define a metric on $2^X$ as follows:
\begin{equation*}
 \rho_H(A,B) = \inf \{\ep>0: A\subseteq B_\ep(B), B\subseteq B_\ep(A)\},
\end{equation*}
where $B_\ep (A)=\{x\in X: \rho(x, A)<\ep\}$ and $\rho(x,A)=\inf_{y\in A} \rho(x,y)$ .
The metric $\rho_H$ is called the {\em Hausdorff metric} of $2^X$.

Let $X,Y$ be two compact metric spaces. Let $F: Y\rightarrow 2^X$ be a map and $y\in Y$.
We say that $F$ is {\em upper semi-continuous (u.s.c.)} at $y$ if for each $\ep>0$ there exists a neighbourhood $U$ of $y$ such that $F(U)\subseteq B_\ep(F(y))$; and $F: Y\rightarrow 2^X$ is {\em lower semi-continuous (l.s.c.)} at $y\in Y$ if and only if for each $\ep>0$ there exists a neighbourhood $U$ of $y$ such that $F(y)\subseteq B_\ep(F(y'))$ for all $y'\in U$.
For more informations on $2^{X}$ and the semi-continuity, refer to \cite{Kur1}.

The following result is well known (for example, see \cite[Section 4]{G94}).
\begin{lem}\label{lsc of orbit map}
Let $(X,G)$ be a $G$-system. Then the map $\Phi: X\rightarrow 2^{X}, x\mapsto \overline{Gx}$ is lower semi-continuous.
\end{lem}

\begin{lem}\cite[p. 394 and p. 70-71]{Kur1,Kur2}
Let $X, Y$ be compact metric spaces. If $F: Y\rightarrow 2^{X}$ is upper semi-continuous (or lower semi-continuous), then the points of continuity of $F$ form a dense $G_{\delta}$ set in $Y$.
\end{lem}

\subsection{Kronecker-Weyl Theorem}
We use $\mathbb{T}^{n}=\mathbb{R}^{n}/\mathbb{Z}^{n}=\{(e^{2\pi i x_1},\ldots, e^{2\pi i x_n}): x_1,\ldots,x_n\in\mathbb{R}\}$ to denote the $n$-torus. The following is the well known Kronecker-Weyl theorem.

\begin{thm}\cite[Theorem 2.2.5]{QQ13}\label{KW}
Let $x_1,\ldots, x_n\in \mathbb{R}$.
\begin{enumerate}
\item[(1)] $\{(e^{2\pi i x_1t},\ldots, e^{2\pi i x_nt}): t\in\mathbb{R}\}$ is dense in $\mathbb{T}^{n}$ if and only if $x_1,\ldots, x_n$ are rationally independent.
\item[(2)] $\{(e^{2\pi i x_1k},\ldots, e^{2\pi i x_nk}): k\in\mathbb{Z}\}$ is dense in $\mathbb{T}^{n}$ if and only if $1, x_1,\ldots, x_n$ are rationally independent.
\end{enumerate}
\end{thm}

For $x_1,x_2,\ldots,x_n\in \R$, it defines  a continuous flow  on the torus $\mathbb{T}^{n}$ by
\[ T^{t}: \mathbb{T}^{n}\rightarrow \mathbb{T}^{n}, (z_1,z_2,\ldots, z_n)\mapsto (z_1e^{2\pi i x_1t}, z_2e^{2\pi i x_2t},\ldots,z_ne^{2\pi i x_nt}), \forall t\in \R.\]
Then  Theorem \ref{KW} implies that  the flow $(\mathbb{T}^{n},\{T^{t}\}_{t\in \R})$ is minimal if and only if $x_1,\ldots, x_n$ are rationally independent and the discrete system $(\mathbb{T}^{n}, T)$ (here $T:=T^{1}$) is minimal if and only if $1, x_1,\ldots, x_n$ are rationally independent.

\section{Minimal elements of a $\R$-flow}\label{section-minimal}

In this section, we show that $(X,\{T^t\}_{t\in \R})$ is minimal if and only if $(X,T^{t})$ is minimal  for all $t\in\mathbb{R}$ except possibly for a countable set. First we show some properties of the regionally proximal relation for commutative actions.

\subsection{Some properties about $\RP$}

Let $G$ be a topological abelian group and $A\subset G$ be a subset. Recall that $A$ is syndetic if there is a compact subset $K\subset G$ such that $G=KA$.

\begin{lem}\label{RP}
Let $(X,G)$ be a $G$-system. Suppose that $H$ is a syndetic subgroup of $G$. Then ${\rm\bf RP}(X, G)={\rm\bf RP}(X, H)$.
\end{lem}

\begin{proof}
Clearly, by definition we have that ${\rm\bf RP}(X, H)\subset{\rm\bf RP}(X, G)$. Next we show the reverse inclusion.
Since $H$ is syndetic, there is a compact subset $K\subset G$ with $G=HK$. Let $(x,y)\in {\rm\bf RP}(X, G)$ and we are going to show that $(x,y)\in {\rm\bf RP}(X, H)$. By the definition of $\RP(X,G)$, there are nets $\{x_{\alpha}\}, \{y_{\alpha}\}$ in $X$  with $x_{\alpha}\rightarrow x, y_{\alpha}\rightarrow y$ and a net $\{g_{\alpha}\}$ in $G$ such that
\[ g_{\alpha}x_{\alpha}\rightarrow z,\ \ \text{ and }\ \ \ g_{\alpha}y_{\alpha}\rightarrow z \]
for some $z\in X$. For each $\alpha$, write $g_{\alpha}=h_{\alpha}k_{\alpha}$ with $h_{\alpha}\in H$ and $k_{\alpha}\in K$. By the compactness of $K$, we may assume that $k_{\alpha}\rightarrow k\in K$. Since the map $G\times X\rightarrow X, (g,x)\mapsto gx$ is continuous, we have that $k_{\alpha}x_{\alpha}\rightarrow kx$ and $k_{\alpha}y_{\alpha}\rightarrow ky$. Now we have that $k^{-1}k_{\alpha}x_{\alpha}\rightarrow x, k^{-1}k_{\alpha}y_{\alpha}\rightarrow y$ and
\begin{align*}
 h_{\alpha}k^{-1}k_{\alpha}x_{\alpha} =k^{-1}h_{\alpha}k_{\alpha}x_{\alpha}=k^{-1} g_{\alpha}x_{\alpha}\rightarrow k^{-1}z,\\
 h_{\alpha}k^{-1}k_{\alpha}y_{\alpha} =k^{-1}h_{\alpha}k_{\alpha}y_{\alpha}=k^{-1}g_{\alpha}y_{\alpha}\rightarrow k^{-1}z.
 \end{align*}
It follows that $(x, y)\in {\rm\bf RP}(X, H)$. The proof is complete.
\end{proof}

\begin{cor}\label{cor-1}
Let $(X,G)$ be a $G$-system. Suppose that $H$ is a syndetic subgroup of $G$. Then $(X,G)$ is equicontinuous if and only if $(X,H)$ is equicontinuous.
\end{cor}

\begin{proof}
It suffices to show that if $(X,H)$ is equicontinuous, then $(X, G)$ is equicontinuous. Since $(X,H)$ is equicontinuous, we have that $\RP(X,H)=\D_X$. By Lemma \ref{RP}, ${\rm\bf RP}(X, G)={\rm\bf RP}(X, H)=\D_X$. Thus $(X,G)$ is equicontinuous.
\end{proof}

Corollary \ref{cor-1} is a well-known result. It applies to all group actions, not just abelian ones, as shown in \cite[Page 46]{Aus88}. If $H$ is not a syndetic subgroup of $G$, then we have the following result.

\begin{lem}\label{equi}
Let $(X,G)$ be a minimal system. Suppose that $H$ is subgroup of $G$ such that $(X,H)$ is also minimal. Then $(X,G)$ is equicontinuous if and only if $(X,H)$ is equicontinuous.
 \end{lem}

\begin{proof}
Clearly, if $(X,G)$ is equicontinuous then so is $(X,H)$. Suppose that $(X,H)$ is equicontinuous and we show that $(X,G)$ is also equicontinuous. Let $E(X,H)$ be the enveloping semigroup for $(X,H)$ (see \cite[Chapter 3]{Aus88} for definition of enveloping semigroup). Since $H$ is abelian, by \cite[Chapter 3, Theorem 3 and Lemma 4]{Aus88} $E(X,H)$ is a compact abelian group and $E(X,H)\subset {\rm Aut}(X,G)$, where ${\rm Aut}(X,G)=\{\phi\in {\rm Homeo}(X): \phi(gx)=g\phi(x), \forall g\in G, x\in X\}$ is the automorphism group of $(X,G)$. Since $E(X,H)$ acts on $X$ transitively in set theoretical sense, it follows from Gottschalk's Theorem \footnote{Gottschalk's Theorem says that for a minimal system $(X,G)$, if ${\rm Aut}(X,G)$ is algebraically transitive, i.e. for all $x,y\in X$, there is a $\phi\in {\rm Aut}(X,G)$ with $\phi x=y$, then $(X,G)$ is equicontinuous (see \cite[Chapter 2, Theorem 13]{Aus88}).} that $(X,G)$ is equicontinuous.
\end{proof}

Let $G,H$ be topological groups. Let $(X,G)$ and $(X,H)$ be systems. We say that $(X,G)$ and $(X,H)$ are {\it commutative} actions, if $ghx=hgx$ for all $g\in G, h\in H$ and all $x\in X$.

The following result was first given by Akin and Glasner in \cite{AG98}. Here we give a different proof.
\begin{thm}
Let $(X,G)$ and $(X,H)$ be minimal systems. If $(X, G)$ and $(X, H)$ are commutative, then ${\rm\bf RP}(X, G)={\rm\bf RP}(X, H)$.
\end{thm}

\begin{proof}
Since $G$ and $H$ are commutative, it is clear that ${\rm\bf RP}(X, G)$ is $H$-invariant. Thus it naturally induces the action of $H$ on $Y=X/{\rm\bf RP}(X, G)$ and the actions of $G$ and $H$ on $Y$ are commutative. Let $L=\langle G,H\rangle$, and it is an abelian group. By Lemma \ref{equi}, $(Y, L)$ is equicontinuous and hence $(Y, H)$ is equicontinuous. Thus by Theorem \ref{thm0} ${\rm\bf RP}(X, G)\supset {\rm\bf RP}(X, H)$. By the symmetry, we obtain the other inclusion and hence ${\rm\bf RP}(X, G)={\rm\bf RP}(X, H)$.
\end{proof}

\subsection{Minimal elements of an $\R$-flow}

We start with the following easy example.
\begin{example}
Let $T^t: \T\rightarrow \T, x\mapsto x+t$. Then $(\T, \{T^{t}\}_{t\in \R})$ is a minimal $\mathbb{R}$-flow, and $(\T, T^t)$ is minimal if and only if $t\not\in \Q$.
\end{example}

Next we show the general case.

\begin{thm}\label{min by eq}
Let $(X, \{T^{t}\}_{t\in \R})$ be a minimal $\mathbb{R}$-flow and let $(X_{eq}, \{T^{t}\}_{t\in \R})$ be the maximal equicontinuous factor of $(X, \{T^{t}\}_{t\in \R})$. Then for each $t\in \mathbb{R}\setminus\{0\}$, $(X, T^{t})$ is minimal if and only if $(X_{eq}, T^{t})$ is minimal.
\end{thm}

\begin{proof}
Fix $t\in\mathbb{R}\setminus\{0\}$. Let $\pi: X\rightarrow X_{eq}$ be the factor map. It is clear that $(X_{eq}, T^{t})$ is minimal whenever $(X, T^{t})$ is minimal. Now assume that $(X_{eq}, T^{t})$ is minimal and we are going to show that $(X, T^{t})$ is minimal.

On the contrary, suppose that $(X, T^{t})$ is not minimal. Since $\{nt: n\in\mathbb{Z}\}$ is a syndetic subgroup of $\mathbb{R}$, each point in $X$ is $T^{t}$-minimal and hence the system $(X, T^{t})$ is semisimple\footnote{A system is called semisimple if each point is minimal. Suppose $(X, G)$ is a system and $S$ is a syndetic normal
subgroup of $G$. Then $(X,G)$ is semisimple if and only if
$(X,S)$ is semisimple (see, for example, \cite[Chapter 1, Theorem 13]{Aus88}). Since we assume all groups are abelian in the paper, all subgroups are normal.}. Decompose $X=\bigcup_{\alpha\in \Gamma}X_{\alpha}$ into the disjoint union of minimal sets with each $X_{\alpha}$ being $T^{t}$-minimal. Since $(X,T^t)$ is not minimal, $|\Gamma|\ge 2$. By the minimality of $(X_{\alpha}, T^{t})$ and $(X_{eq}, T^{t})$, we have
that or each $\alpha\in \Gamma$, $\pi|_{X_\a}: X_{\alpha}\rightarrow X_{eq}$ is surjective. Thus for any $\alpha\neq \beta\in \Gamma$,  $R_{\pi}\cap (X_{\alpha}\times X_{\beta})\neq\emptyset$, where $R_\pi=\{(x_1,x_2)\in X^2: \pi(x_1)=\pi(x_2)\}$.

%\medskip

Let $\Phi: X\rightarrow 2^{X}, x\mapsto \overline{\mathcal{O}_{T^{t}}}(x)$ be the orbit closure map with respect to $T^{t}$. By Lemma \ref{lsc of orbit map}, the set $X_c$ of continuous points of $\Phi$ is residual in $X$. Thus there are $\alpha\neq \beta\in \Gamma$ and $x\in X_{\alpha}, y\in X_{\beta}$ such that $\Phi$ is continuous at $x$ and $y$. (Otherwise, there is some $\a\in\Gamma$ such that $X_0\subseteq X_\a$, and it follows that $X=X_\a$ is minimal contradicting with our assumption.)
By Lemma \ref{RP}, we have $R_{\pi}={\rm\bf RP}(X, \mathbb{R})={\rm\bf RP}(X, T^{t})$. Thus, for $\alpha\neq \beta$, by $R_{\pi}\cap (X_{\alpha}\times X_{\beta})\neq\emptyset$, there is $x_{\alpha}\in X_{\alpha}$ and $x_{\beta}\in X_{\beta}$ with $(x_{\alpha}, x_{\beta})\in {\rm\bf RP}(X, T^{t})$.

Let $\varepsilon=\rho_{H}(X_{\alpha}, X_{\beta})>0$, where $\rho_{H}$ is the Hausdorff metric on $2^{X}$. By the continuity of $\Phi$ at $x$ and $y$, there is some $\delta>0$ such that
\[ \bigcup_{z\in B_{\delta}(x)} \Phi(z) \subset B_{\varepsilon/2}(\Phi(x))\ \ \text{ and }\ \ \bigcup_{z\in B_{\delta}(y)} \Phi(z) \subset B_{\varepsilon/2}(\Phi(y)).\]
By the minimality of $(X_{\alpha}, T^{t})$ and  $(X_{\beta}, T^{t})$, there are $m,n\in\mathbb{Z}$ and open neighborhoods $U, V$ of $x_{\alpha}, x_{\beta}$ respectively such that $T^{mt}U\subset B_{\delta}(x)$ and $T^{nt}V\subset B_{\delta}(y)$. Since $(x_{\alpha}, x_{\beta})\in {\rm\bf RP}(X, T^{t})$, there is $u\in U, v\in V$  and $k\in \mathbb{Z}$ such that $\rho(T^{kt}u, T^{kt}v)<\varepsilon/3$. However, this is impossible by noting that
\[ \Phi(u) =\Phi(T^{mt}u)\subset B_{\varepsilon/2}(\Phi(x)) , \ \  \Phi(v) =\Phi(T^{nt}v)\subset B_{\varepsilon/2}(\Phi(y)), \]
 and $\rho_{H}(\Phi(x), \Phi(y))=\varepsilon$ as $\Phi(x)=X_{\alpha}$ and $\Phi(y)=X_{\beta}$.
This contradiction shows that $(X, T^{t})$ is minimal. The proof is complete.
\end{proof}

It is well-known that a minimal $G$-system $(X,G)$ is weakly mixing if and
only if $X_{eq}$ is trivial (for example, see \cite[Chapter 9, Theorem 13]{Aus88}). Thus we have the following corollary.

\begin{cor}\label{cor-2}
If $(X, \{T^{t}\}_{t\in \R})$ is a minimal $\mathbb{R}$-flow and it is weakly mixing, then for each $t\in \mathbb{R}\setminus\{0\}$, $(X, T^{t})$ is minimal and weakly mixing.
\end{cor}

\begin{lem}\label{min time for flow}
Let $(X,\{T^{t}\}_{t\in \R})$ be a minimal equicontinuous $\mathbb{R}$-flow. Then there is a countable set $E\subset \mathbb{R}$ such that for each $t\in \mathbb{R}\setminus E$, $(X,T^{t})$ is also minimal.
\end{lem}

\begin{proof}
Since $(X,\{T^{t}\}_{t\in \R})$ be a minimal equicontinuous, it is isomorphic to an $\R$-action on a  compact metric abelian group (see, for example, \cite[Chapter 3, Theorem 6]{Aus88}).
We may assume that $X$ is a compact connected metric abelian group. Then $(X,\mathbb{R})$ is conjugate to an inverse system $\underset{\longleftarrow}{\lim}(H_i, \mathbb{R})_{j\in \mathbb{N}}$, where each $H_j$ is a compact abelian Lie group. Since $(H_j, \mathbb{R})$ is minimal  for each $j\in\mathbb{N}$, $H_j$ is connected and hence is a torus. Further, we may assume that $H_j=\mathbb{T}^{n_j}$ and the flow is define by $\mathbb{R}\times H_j\rightarrow H_j, (t, z)\mapsto z(e^{2\pi i \alpha_{1}^{(j)}t},\ldots, e^{2\pi i \alpha_{n_{j}}^{(j)}t})$ with $(\alpha_1^{(j)},\ldots,\alpha_{n_j}^{(j)} )\in \mathbb{R}^{n_i}$. By Theorem \ref{KW}, one has that $\alpha_1^{(j)},\ldots,\alpha_{n_j}^{(j)} $ are rationally independent.

Now for $t\in\mathbb{R}$, by Theorem \ref{KW} $(H_j, T^{t})$ is minimal if and only if $1, \alpha_1^{(j)}t,\ldots,\alpha_{n_j}^{(j)}t $ are rationally independent. Let
\[ E_{j}=\left\{ \frac{r}{\alpha_1^{(i)}s_1+\cdots+\alpha_{n_j}^{(j)}s_{n_j} }: r\in\mathbb{Q},(s_1,\ldots,s_{n_j})\in\mathbb{Q}^{n_j}\setminus\{{\bf 0}\}\right\}.\]
Then $E_{j}$ is countable and for each $t\notin E_j$, $(H_j, T^{t})$ is minimal.

Set $E=\bigcup_{j=1}^{\infty} E_j$, which is countable. Since $(X,T^{t})$ is minimal if and only if $(H_j, T^{t})$ is minimal for each $j\in\mathbb{N}$, we have that for each $t\notin E$, $(X,T^{t})$ is minimal. The proof is complete.
\end{proof}

Now combining Theorem \ref{min by eq}, Lemma \ref{min time for flow} and Corollary \ref{cor-2}, we conclude that

\begin{thm}\label{minimality}
Let $(X, \{T^{t}\}_{t\in \R})$ be a minimal $\mathbb{R}$-flow.
\begin{enumerate}
  \item If $(X, \{T^{t}\}_{t\in \R})$ is weakly mixing, then $(X, T^{t})$ is minimal for each $t\in \R\setminus\{0\}$.
  \item If $(X, \{T^{t}\}_{t\in \R})$ is not weakly mixing, then there is a countable set $E\subset \mathbb{R}$ such that $(X, T^{t})$ is minimal for each $t\notin E$.
\end{enumerate}

\end{thm}

\subsection{A question}

In \cite{Gutschera},  the author investigates the conditions under which a connected Lie group, acting ergodically on a space with a finite invariant measure, contains a single element that acts ergodically.  Roughly speaking, the answer in \cite{Gutschera} is that such an element exists for most types of groups but not for compact non-Abelian groups or for certain types of solvable groups. We want to ask the same question in topological setting, replacing ergodicity by minimality.

\section{Regionally proximal relation of order $d$ for commutative actions}\label{section-RPd}

In this section we study the regionally proximal relation of order $d$ for commutative actions.
Let $(X,G)$ and $(X,H)$ be minimal systems, and assume that they are commutative actions. Then for any $d\in\mathbb{N}$,  we show that ${\RP}^{[d]}(X, G)={\RP}^{[d]}(X, H)$ and ${\bf Q}^{[d]}(X, G)={\bf Q}^{[d]}(X, H)$.
For two commutative minimal system $(X,T)$ and $(X,S)$, they have the same increasing sequence of pro-nilfactors:
$$\{pt\}= X_0  \longleftarrow  X_1 \longleftarrow \cdots  \longleftarrow X_n  \longleftarrow X_{n+1}
 \longleftarrow \cdots \longleftarrow X,$$
where $X_d=X/\RP^{[d]}(X,T)$ is the the
maximal $d$-step pro-nilfactor.
Further, we will show that $N_d(X,T)=N_d(X,S)$ for $d\in \N$.

\subsection{Regionally proximal relation of order $d$ for commutative actions}
The following lemma is a generation of Lemma \ref{RP}.
\begin{lem}\label{RPd for synd}
Let $(X,G)$ be a $G$-system. Suppose that $H$ is a  syndetic subgroup of $G$. Then ${\rm\bf RP}^{[d]}(X, G)={\rm\bf RP}^{[d]}(X, H)$ for any $d\in\mathbb{N}$.
\end{lem}

\begin{proof}
Fix $d\in\mathbb{N}$. Clearly, by the definition  ${\rm\bf RP}^{[d]}(X, H)\subset{\rm\bf RP}^{[d]}(X, G)$. Next we show the reverse inclusion. That is, we show that if $(x,y)\in {\rm\bf RP}^{[d]}(X, G)$, then $(x,y)\in {\rm\bf RP}^{[d]}(X, H)$.

Since $H$ is syndetic, there is a compact subset $K\subset G$ with $G=HK$. Fix a $\delta>0$. Since $K$ is compact, the map $K^{-1}\times X\rightarrow X, (k^{-1}, x)\mapsto k^{-1}x$ is uniformly continuous. Thus there is some $\eta_1\in (0,\d/2)$ such that whenever $p,q\in X$ with $\rho(p,q)<\eta_1$, we have that $\rho(k^{-1}p,k^{-1}q)<\d/3$ for all $k\in K$.
Since $(x,y)\in {\rm\bf RP}^{[d]}(X, G)$, there are $x',y'\in X$ with $\rho(x,x')<\eta_1< \delta/2, \rho(y,y')<\eta_1< \delta/2$ and ${\bf g}=(g^{(1)}, \ldots,g^{(d)})\in G^{d}$ such that for any ${\bf \epsilon}=(\epsilon_1,\ldots,\epsilon_{d})\in \{0,1\}^{d}\setminus\{{\bf 0}\}$,
$\rho \left({\bf  g}^{({\bf \epsilon})} x', {\bf  g}^{({\bf \epsilon})}y'\right)<\eta_1.$
Recall that here ${\bf g}^{\bf \epsilon}=\prod_{i=1}^{d}g_{i}^{\epsilon_i}$. Thus $$\rho \left(k^{-1}{\bf  g}^{({\bf \epsilon})} x', k^{-1}{\bf  g}^{({\bf \epsilon})}y'\right)<\d/3, \ \forall k\in K.$$

Write ${\bf g}={\bf h} {\bf k}$ with ${\bf h} \in H^{d}$ and $ {\bf k}\in K^{d}$. Again by the uniform continuity of $K^{-1}\times X\rightarrow X, (k^{-1}, x)\mapsto k^{-1}x$, there is an $\eta_2\in (0,\delta/2)$ such that whenever $p,q\in X$ with $\rho(p,q)<\eta_2$, we have
\[ \max_{{\bf \epsilon}\in \{0,1\}^{d}\setminus\{{\bf 0}\}}\max_{k\in K}\rho\left( {\bf g}^{({\bf \epsilon})}k^{-1} p, {\bf g}^{({\bf \epsilon})} k^{-1}q\right)<\delta/3.\]
Now choose $x'',y''\in X$ such that $\rho(x',x'')<\eta_2$ and $\rho(y',y'')<\eta_2$. Then $\rho(x,x'')<\delta$ and $\rho(y,y'')<\delta$. And
\begin{align*}
 & \rho \left({\bf  h}^{({\bf \epsilon})} x'', {\bf  h}^{({\bf \epsilon})}y''\right) =\rho \left({\bf  g}^{({\bf \epsilon})}( {\bf  k}^{({\bf \epsilon})})^{-1} x'', {\bf  g}^{({\bf \epsilon})}( {\bf  k}^{({\bf \epsilon})})^{-1}y''\right) \\
 \leq &\rho \left({\bf  g}^{({\bf \epsilon})}( {\bf  k}^{({\bf \epsilon})})^{-1} x'', {\bf  g}^{({\bf \epsilon})}( {\bf  k}^{({\bf \epsilon})})^{-1}x'\right)+\rho \left({\bf  g}^{({\bf \epsilon})} ( {\bf  k}^{({\bf \epsilon})})^{-1}x', {\bf  g}^{({\bf \epsilon})}( {\bf  k}^{({\bf \epsilon})})^{-1}y'\right)\\ \quad \quad \quad & +\rho \left({\bf  g}^{({\bf \epsilon})}( {\bf  k}^{({\bf \epsilon})})^{-1}  y', {\bf  g}^{({\bf \epsilon})}( {\bf  k}^{({\bf \epsilon})})^{-1}y''\right)\\
  <&\delta/3+\delta/3+\delta/3=\delta.
 \end{align*}
Thus we have $(x,y)\in {\rm\bf RP}^{[d]}(X, H)$. The proof is complete.
\end{proof}

Without the assumption on the syndeticity of $H$, we show Lemma \ref{RPd for synd} also holds provided that both $(X,G)$ and $(X,H)$ are transitive. First, we give the following lemma. 

\begin{lem}\label{RPd1}
Let $(X,G)$ and $(X,H)$ be commutative actions and both transitive. Then for any nonempty open set $W$ in $X$ and any ${\bf g}\in G^{d}$, there exists ${\bf h}\in H^{d}$ such that
\[ W\cap \bigcap_{\epsilon\in\{0,1\}^{d}\setminus\{{\bf 0}\}} ({\bf h}^{(\epsilon)})^{-1} {\bf g}^{(\epsilon)} W\neq\emptyset.\]
\end{lem}

\begin{proof}
Let ${\bf g}=(g_1,\ldots, g_d)\in G^d$. Now we construct ${\bf h}=(h_1,\ldots, h_d)\in H^d$ inductively.
First, by the transitivity of $(X, H)$ that there is some $h_1\in H$ with $h_1W\cap g_1W\neq\emptyset$. Set $W_1=W\cap (h_1^{-1}g_1W)$. Again by the transitivity of $(X, H)$, there is some $h_2\in H$ with $h_2W_1\cap g_2W_1\neq\emptyset$. Set $W_2=W_1\cap (h_2^{-1}g_2W_1)$. Inductively, for $i<d $, assume that we have found $h_1, \ldots, h_i\in H$ and open sets $W_1,\ldots, W_i$ in $X$ such that for each $j\in\{1,\ldots,i\}$,
\begin{itemize}
\item $W_{j}=W_{j-1}\cap( h_{j}^{-1}g_{j} W_{j-1})\neq\emptyset$, where we set $W_0=W$.
\end{itemize}

Then it follows from the  transitivity of $(X,H)$ that there is some $h_{i+1}\in H$ such that $h_{i+1}W_{i}\cap g_{i+1}W_{i}\neq\emptyset$ and we set $W_{i+1}=W_{i}\cap( h_{i+1}^{-1}g_{i+1} W_{i})$. In such inductive way, we construct $h_1,\ldots, h_d\in H$ and open sets $W_1,\ldots, W_d$ such that
\begin{itemize}
\item $W=W_0\supset W_1\supset\cdots\supset W_d$;
\item $W_{j}=W_{j-1}\cap( h_{j}^{-1}g_{j} W_{j-1})\neq\emptyset$, for all $j\in \{1,\ldots,d\}$.
\end{itemize}
Let ${\bf h}=(h_1,\ldots, h_d)$. Next we verify that for each $\epsilon\in\{0,1\}^{d}\setminus\{{\bf 0}\}$, $ W_{d}\subset ({\bf h}^{(\epsilon)})^{-1}{\bf g}^{(\epsilon)} W\neq\emptyset$. For this, we fix an $\epsilon\in\{0,1\}^{d}\setminus\{{\bf 0}\}$ and let $i_1<i_2<\cdots<i_r$ be the coordinates in $\epsilon$ with $\epsilon_{i_1}=\cdots=\epsilon_{i_r}=1$.  Then we have
\begin{align*}
{\bf h}^{(\epsilon)} W_{d}&= h_{i_1}\cdots h_{i_{r}} W_{d}\subset h_{i_1}\cdots h_{i_{r}}  W_{i_{r}}
= h_{i_1}\cdots h_{i_{r-1}} (h_{i_{r}}W_{i_{r}-1}\cap g_{i_{r}}W_{i_{r}-1})\\
&\subset h_{i_1}\cdots h_{i_{r-1}}   g_{i_{r}}W_{i_{r}-1}
\subset h_{i_1}\cdots h_{i_{r-1}}   g_{i_{r}}W_{i_{r-1}}\\
&= h_{i_1}\cdots h_{i_{r-1}}   g_{i_{r}}(W_{i_{r-1}-1}\cap h_{i_{r-1}}^{-1}g_{i_{r-1}} W_{i_{r-1}-1})\\
&\subset h_{i_1}\cdots h_{i_{r-1}}   g_{i_{r}} h_{i_{r-1}}^{-1}g_{i_{r-1}} W_{i_{r-1}-1}= h_{i_1}\cdots h_{i_{r-2}}   g_{i_{r}}g_{i_{r-1}} W_{i_{r-1}-1}\\
&\subset h_{i_1}\cdots h_{i_{r-2}}   g_{i_{r}}g_{i_{r-1}} W_{i_{r-2}} \subset \cdots \\
&\subset g_{i_r}\cdots g_{i_1} W= {\bf g}^{(\epsilon)} W.
\end{align*}
Thus $W_{d}\subset ({\bf h}^{(\epsilon)} )^{-1} {\bf g}^{(\epsilon)} W\neq\emptyset$ and the proof is complete.
\end{proof}

\begin{thm}\label{RPd}
Let $(X,G)$ be a transitive system. If $H$ is a subgroup of $G$ such that $(X,H)$ is transitive. Then for any $d\in\mathbb{N}$,  we have ${\rm\bf RP}^{[d]}(X, G)={\rm\bf RP}^{[d]}(X, H)$.
\end{thm}
\begin{proof}
Recall that the regionally proximal relation of oder $d$ is
\[ {\rm\bf RP}^{[d]}(X, G)=\bigcap_{\alpha\in \mathcal{U}}\overline{ \bigcup_{{\bf g}\in G^{d}} \bigcap_{\epsilon\in\{0,1\}^{d}\setminus\{{\bf 0}\}}({\bf g}^{(\epsilon)} )^{-1}\alpha},\]
where $\mathcal{U}$ is the uniform structure on $X$.
Next we show that for each $\a\in \mathcal{U}$, we have that
\begin{equation}\label{xs1}
  \overline{ \bigcup_{{\bf g}\in G^{d}} \bigcap_{\epsilon\in\{0,1\}^{d}\setminus\{{\bf 0}\}}({\bf g}^{(\epsilon)} )^{-1}\alpha} ~\subset~ \overline{ \bigcup_{{\bf h}\in H^{d}} \bigcap_{\epsilon\in\{0,1\}^{d}\setminus\{{\bf 0}\}}({\bf h}^{(\epsilon)} )^{-1}\alpha}.
\end{equation}
Then it will follows that ${\bf RP}^{[d]}(X,G)\subset {\bf RP}^{[d]}(X,H)$ whence ${\bf RP}^{[d]}(X,G)= {\bf RP}^{[d]}(X,H)$ since it is clear that ${\bf RP}^{[d]}(X,G)\supset {\bf RP}^{[d]}(X,H)$.

\medskip

Now we fix $\alpha\in \mathcal{U}$. Let $\displaystyle (x,y)\in \overline{ \bigcup_{{\bf g}\in G^{d}} \bigcap_{\epsilon\in\{0,1\}^{d}\setminus\{{\bf 0}\}}({\bf g}^{(\epsilon)} )^{-1}\alpha}$. Then for any open neighbourhood $E$ of $(x,y)$, $\displaystyle E\cap \Big( \bigcup_{{\bf g}\in G^{d}} \bigcap_{\epsilon\in\{0,1\}^{d}\setminus\{{\bf 0}\}}({\bf g}^{(\epsilon)} )^{-1}\alpha\Big)\neq \emptyset$. Choose non-empty open sets $U, V$ of $X$ such that $U\times V\subset E$ and
$$\bigcup_{\epsilon\in\{0,1\}^{d}\setminus\{{\bf 0}\} }{\bf g}^{(\epsilon)}(U\times V)\subset  \alpha \text{ for some } {\bf g}=(g_1,\ldots, g_{d})\in G^{d}.$$
Since $G$ acts on $X$ transitively, there is some $g\in G$ with $gV\cap U\neq\emptyset$. We set $W=gV\cap U$.
By Lemma \ref{RPd1}, there is some ${\bf h}\in H^{d}$ such that for any  $\epsilon\in\{0,1\}^{d}\setminus\{{\bf 0}\}$, $$O_{\epsilon}:={\bf h}^{(\epsilon)} W\cap {\bf g}^{(\epsilon)}W\neq\emptyset.$$
Since $W=gV\cap U$, we have
\[ O_{\epsilon}\subset{\bf h}^{(\epsilon)}W\subset {\bf h}^{(\epsilon)}U\ \ \text{ and } \ \ O_{\epsilon}\subset {\bf g}^{(\epsilon)}W\subset {\bf g}^{(\epsilon)}U.\]
And
\[g^{-1}O_{\epsilon}\subset g^{-1}{\bf h}^{(\epsilon)}W={\bf h}^{(\epsilon)}g^{-1}W\subset  {\bf h}^{(\epsilon)}V\ \ \text{ and }\ \ g^{-1}O_{\epsilon}\subset g^{-1}{\bf g}^{(\epsilon)} W={\bf g}^{(\epsilon)} g^{-1}W\subset  {\bf g}^{(\epsilon)} V.\]
These imply that
\[ \emptyset\neq O_{\epsilon}\times g^{-1}O_{\epsilon}\subset {\bf h}^{(\epsilon)}(U\times V)\cap {\bf g}^{(\epsilon)}(U\times V)\subset {\bf h}^{(\epsilon)}(U\times V)\cap \alpha.\]
Thus $(U\times V)\cap ({\bf h}^{(\ep)})^{-1}\a\neq\emptyset$ for all $\epsilon\in\{0,1\}^{d}\setminus\{{\bf 0}\}$. In particular,
$$E\cap \Big( \bigcup_{{\bf h}\in H^{d}} \bigcap_{\epsilon\in\{0,1\}^{d}\setminus\{{\bf 0}\}}({\bf h}^{(\epsilon)} )^{-1}\alpha\Big)\neq \emptyset.$$
Since $E$ is an arbitrary neighbourhood of $(x,y)$, $\displaystyle (x,y)\in \overline{ \bigcup_{{\bf h}\in H^{d}} \bigcap_{\epsilon\in\{0,1\}^{d}\setminus\{{\bf 0}\}}({\bf h}^{\epsilon} )^{-1}\alpha}$. That is,
\[\overline{ \bigcup_{{\bf g}\in G^{d}} \bigcap_{\epsilon\in\{0,1\}^{d}\setminus\{{\bf 0}\}}({\bf g}^{\epsilon} )^{-1}\alpha} ~\subset~ \overline{ \bigcup_{{\bf h}\in H^{d}} \bigcap_{\epsilon\in\{0,1\}^{d}\setminus\{{\bf 0}\}}({\bf h}^{\epsilon} )^{-1}\alpha}.\]
The proof is complete.
\end{proof}

As a corollary, we have Theorem \ref{thm-RPd}:
Let $(X,G)$ and $(X,H)$ be commutative actions and both transitive and $d\in \N$. Then $\RP^{[d]}(X,G)=\RP^{[d]}(X,H)$.

\begin{proof}[Proof of Theorem \ref{thm-RPd}]
Let $L=\langle G,H\rangle$. Then $L$ is an abelian group and both $G,H$ are its subgroups. Since $(X,G)$ is transitive, $(X,L)$ is also transitive. By Theorem \ref{RPd},
$$\RP^{[d]}(X,G)=\RP^{[d]}(X,L)=\RP^{[d]}(X,H).$$
The proof is complete.
\end{proof}

\subsection{${\bf Q}^{[d]}(X,G)$ }

The cube structure plays great role in the Host and Kra's work. We recall the notions of dynamical cubes. Let $(X,G)$ be a $G$-system.

A set of the form $F=\{ \epsilon\in\{0,1\}^{d}: \epsilon_{i_1}=\alpha_1,\ldots,\epsilon_{i_{d-1}}=\alpha_{d-1}\}$ for some $1\leq i_1<i_2<\cdots<i_{d-1}\leq d$ and $\alpha_{i}\in\{0,1\}$ is called a {\it hyperface} of the discrete cube $\{0,1\}^{d}$. Such a hyperface is called an {\it upper hyperface} if all $\alpha_1=\cdots=\alpha_{d-1}=1$. For a hyperface $F$ and $g\in G$, define $g_{F}\in G^{\{0,1\}^{d}}$ by
\[g_{F}(\epsilon)=\begin{cases}
g & \text{ if } \epsilon\in F,\\
e & \text{ if } \epsilon\notin F.
\end{cases}
\]
Now the $d$-dimensional {\it cube group} of $G$, denoted by $\mathcal{G}^{[d]}(G)$, is the subgroup of $G^{\{0,1\}^{d}}$ generated by $\{g_{F}: g\in G, F \text{ is a hyperface of } \{0,1\}^{d}\}$.   The $d$-dimensional face cube group of  $G$, denoted by ${\bf \mathcal{F}}^{[d]}(G)$, is the subgroup of $G^{\{0,1\}^{d}}$ generated by $\{g_{F}: g\in G, F$ is an upper hyperface of $\{0,1\}^{d}\}$.  There are  canonically defined surjective group homomorphisms from $G^{d+1}$ to $\mathcal{G}^{[d]}(G)$ by
\[(g_0,g_1,\ldots, g_{d})\mapsto \left(g_0g_1^{\ep_1}g_2^{\ep_2}\cdots g_{d}^{\ep_d}: (\ep_1,\ldots,\ep_d)\in\{0,1\}^{d}\right)\]
and  $G^{d}$ to ${\bf \mathcal{F}}^{[d]}(G)$ by
\[(g_1,\ldots, g_{d})\mapsto \left(g_1^{\ep_1}g_2^{\ep_2}\cdots g_{d}^{\ep_d}: (\ep_1,\ldots,\ep_d)\in\{0,1\}^{d}\right),\]
respectively. Now the action of $G$ on $X$ naturally induces the action of $\mathcal{G}^{[d]}(G)$ on $X^{\{0,1\}^{d}}$. The {\it $d$-dimensional cube} of $(X,G)$ is defined to be the closure of the orbit of $\Delta_{2^{d}}(X)=\{(x,x,\ldots,x)\in X^{\{0,1\}}:x\in X\}$ under $\mathcal{G}^{[d]}(G)$.
When $(X,G)$ is minimal,
\[ {\bf Q}^{[d]}(X,G)=\overline{\mathcal{G}^{[d]}(G)(x,x,\ldots,x)}=\overline{{\bf \mathcal{F}}^{[d]}(G)\Delta_{2^{d}}(X)},\]
for any $x\in X$.

\begin{thm}\label{Qd}
Let $(X,G)$ and $(X,H)$ be commutative actions and both transitive. Then for any $d\in\mathbb{N}$,  we have ${\bf Q}^{[d]}(X, G)={\bf Q}^{[d]}(X, H)$.
\end{thm}

\begin{proof}
Fix ${\bf x}=(x_{\epsilon})_{\epsilon\in\{0,1\}^{d}}\in {\bf Q}^{[d]}(X, G)$ and fix neighborhood $V_{\epsilon}$ of $y_{\epsilon}$ for each $\epsilon\in\{0,1\}^{d}$. Then there is a neighborhood $W$ of $x_{{\bf 0}}$ and ${\bf g}\in G^{d}$ with ${\bf g}^{(\epsilon)}W\subset V_{\epsilon}$, for each $\epsilon\in\{0,1\}^{d}\setminus \{{\bf 0}\}$. By Lemma \ref{RPd1}, there is some ${\bf h}^{(\epsilon)}\in H$ with
\[ U:=W\cap \bigcap_{\epsilon\in\{0,1\}^{d}\setminus\{{\bf 0}\}} ({\bf h}^{(\epsilon)})^{-1}{\bf g}^{(\epsilon)} W\neq\emptyset.\]
Thus for any $x\in U$ and each $\epsilon\in\{0,1\}^{d}\setminus \{{\bf 0}\}$, we have ${\bf h}^{\epsilon}(x)\in {\bf g}^{\epsilon}W\subset V_{\epsilon}$. It follows that ${\bf x}\in {\bf Q}^{[d]}(X, H)$ and hence ${\bf Q}^{[d]}(X, G)\subset{\bf Q}^{[d]}(X, H)$.  By the symmetry, we have ${\bf Q}^{[d]}(X, G)={\bf Q}^{[d]}(X, H)$.
\end{proof}

Consequently, we have shown the first assertion of Theorem \ref{cube and Nd}.

%Theorem \ref{RPd} also follows from Theorem \ref{Qd} and the following characterization  of regionally proximal relation of higher order via cube group action.

%\begin{thm}[\cite{SY12}]
%Let $(X,G)$ be a minimal system. Then the following are equivalent:
%\begin{enumerate}
%\item $(x,y)\in {\bf RP}^{[d]}(X,G)$;
%\item $(x, y,y,\ldots,y)\in {\bf Q}^{[d]}(X,G)$.
%\end{enumerate}
%\end{thm}

\subsection{$N_{d}(X, S)=N_{d}(X,T)$ for commuting transformations $T$ and $S$}\label{Nd for commuting}
Let $(X,T)$ be a system and $d\in\mathbb{N}$.
Set $\sigma_{d}(T)=T\times\cdots\times T$ and $\tau_{d}(T)=T\times T^2\times\cdots\times T^{d}$. Recall that $N_{d}(X,T)$ is defined to be the closure of
\[\{(T^{n}x, T^{2n}x,\ldots, T^{dn}x): x\in X, n\in\mathbb{Z}\}\]
in $X^{d}$. When $(X,T)$ is transitive and $x\in X$ is a transitive point, $N_d(X,T)=\overline{\O}(x, G_d(T))$, where $G_d(T)=\langle\sigma_{d}(T), \tau_{d}(T)\rangle$. Galsner showed that the $\Z^2$-system $(N_{d}(X,T),G_d(T))$ is minimal whenever $(X,T)$ is a minimal system \cite{G94}.

In this subsection, we show that if $S$ and $T$ are commutative minimal homeomorphisms on $X$, then $N_{d}(X, S)=N_{d}(X,T)$.

\medskip

First we need to study nilsystems.

\begin{thm}\cite[Theorem 2.2]{Zie05}
Let $(G/\Gamma, m_{1}, T)$ be a $k$-step ergodic nilsystem and let  $f_1,\ldots, f_{k}\in L^{\infty}(G/\Gamma)$. If $G$ is spanned by $\overset{\circ}{G}$ and $a$, where  $\overset{\circ}{G}$ is connected component of the identity and $a$ is the translating element determining $T$, then for almost all $x\in G$,
\begin{align*}
&\lim_{N\rightarrow\infty}\frac{1}{N}\sum_{n=1}^{N}f_1(ax\Gamma)f_{2}(a^2x\Gamma)\cdots f_{k}(a^{k}x\Gamma)\\
=&\int_{G/\Gamma}\left(\cdots\left(\int_{G_{k}/\Gamma_{k}}\prod_{j=1}^{k}f_{j}\left(x\prod_{i=1}^{j}y_{i}^{{j\choose i}}\Gamma\right) dm_{k}(y_{k}\Gamma_{k})\right)\cdots\right)dm_1(y_1\Gamma_1).
\end{align*}
\end{thm}

More precisely, Ziegler gave the structure $A_{x}=\overline{\mathcal{O}}(x^{(d)}, a\times a^2\times\cdots\times a^{k})$ as follows. Define  $\tilde{G}=G_1\times \cdots\times G_{k}$, which is endowed with the multiplication by
\[ (x_1,\ldots,x_k)\star (y_1,\ldots, y_k)=(z_1,\ldots,z_k),\]
where
\[ \prod_{j=1}^{i}z_{j}^{{i\choose j}}=\prod_{j=1}^{i}x_{j}^{{i\choose j}}\prod_{j=1}^{i}y_{j}^{{i\choose j}}, 1\leq i\leq k.\]
Let $\tilde{\Gamma}=\Gamma_1\times \Gamma_2\times\cdots\times \Gamma_{k}$. Then $(\tilde{\Gamma},\star)$ is a cocompact lattice of $\tilde{G}$. Now for each $x\in G$, define a transformation $S_{x}: \tilde{G}/\Gamma\rightarrow \tilde{G}/\Gamma$ by
\[S_{x}((y_1,\ldots,y_{k})\star\tilde{\Gamma})=((a[a,x],e,\ldots,e)\star(y_1,\ldots,y_{k}))\star\tilde{\Gamma}, \]
 and a mapping $I_{x}: \tilde{G}\rightarrow (G/\Gamma)^{k}$ by
\[ I_{x}(y_1,\ldots,y_{k})=\left(xy_1\Gamma, xy_{1}^{2}y_2\Gamma, \ldots, x\prod_{j=1}^{k}y_{j}^{{k\choose j}}\Gamma\right).\]
Then it is shown in \cite{Zie05} that $I_{x}$ is an isomorphism of the systems $(\tilde{G}/\tilde{\Gamma}, S_{x})$ and $(A_{x}, T\times\cdots\times T^{k})$ and  $S_{x}$ is minimal for almost every $x$.

Note that $I_{x}(\tilde{G})$ is independent of $a$ and hence $A_{x}$ is independent of the rotation $a$.  As a corollary, we have

\begin{cor}\label{Nd for nil}
If both $(G/\Gamma, S)$ and $(G/\Gamma, T)$ are   minimal nilsystems. Then $N_{d}(G/\Gamma, S)=N_{d}(G/\Gamma,T)$ for each $d\in\mathbb{N}$.
\end{cor}
%By passing to the inverse limit, we have
%\begin{cor}
%If both $(X, S)$ and $(X, T)$ are $k$-step minimal pro-nilsystems. Then $N_{k}(X,S)=N_{k}(X,T)$.
%\end{cor}
%{\color{red}add a proof, since $X$ is an inverse limits}

%\begin{lem}
%Let $(X, T)=\underset{\longleftarrow}{\lim}(X_i, T)$ be a minimal system. Then for each $d$, $N_{d}(X,T)=\underset{\longleftarrow}{\lim}  N_{d}(X_i,T)$.
%\end{lem}
%\begin{cor}\label{N_d for nilsystem}
%Let $(X,S)$ and  $(X,T)$ be minimal  pro-nilsystems. Then $N_{d}(X,S)=N_{d}(X,T)$ for each $d\in\mathbb{N}$.
%\end{cor}

\begin{lem}\label{N_d for proximal extension}
Let $\pi: X\rightarrow Y$ be a proximal extension between minimal systems both for $S$ and $T$. If $S$ and $T$ are commuting and $N_{d}(Y,S)=N_{d}(Y, T)$, then we have  $N_{d}(X,S)=N_{d}(X, T)$.
\end{lem}

\begin{proof}
It suffices to show that $N_d(X, S)$ is $\tau_d(T)=T\times T^2\times\cdots\times T^{d}$-invariant. In fact, if $N_d(X, S)$ is $\tau_d(T)$-invariant, then
$N_d(X,T)\subset N_d(X,S)$ since $(N_{d}(X,T),G_d(T))$ is minimal. Then $N_{d}(X,S)=N_{d}(X, T)$ is followed by the symmetry.

Let $(x_1,\ldots, x_d)\in N_d(X, S)$.
By the assumption, $(T\pi(x_1),\ldots, T^{d}\pi(x_{d}))\in N_d(Y, S)$. By the minimality of $ (N_d(Y, S),G_d(S))$, there is a sequence $(g_i)\in G_d(S)$ such that
$$ g_i(T\pi(x_1),\ldots, T^{d}\pi(x_d))\rightarrow (y,\ldots,y),$$ for some $y\in Y$.
Without loss of generality, we may assume that
$$g_i(Tx_1,\ldots, T^{d}x_d)\rightarrow (x_1',\ldots, x_{d}')\in (\pi^{-1}(y))^d$$
for some $x_1',\ldots, x_d'$. Since $\pi$ is proximal, $\overline{\mathcal{O}}_{S\times \cdots\times S}(x_1',\ldots, x_d')\subset \pi^{-1}(\Delta_{d}(Y))$. Thus the orbit closure of $(Tx_1,\ldots T^{d}x_d)$ under $\langle \sigma_{d}(S),\tau_{d}(S)\rangle$ meets the diagonal $\Delta_{d}(X)$. Then we have $(Tx_1,\ldots T^{d}x_d)\in N_{d}(X, S)$, since $(N_d(X, S), G_{d}(S))$ is minimal. This shows that $N_d(X, S)$ is $\tau_d(T)$-invariant. The proof is complete.
\end{proof}

\begin{lem}\label{N_d for char factor}\cite[Lemma 3.3]{GHSWY}
Let $\pi: (X, T) \rightarrow (Y, T)$ be an open extension between minimal systems. Suppose that $(Y, T)$ is a $d$-step characteristic factors of $(X,T)$. Then  $N_{d+1}(Y,T)$ is $\pi^{(d+1)}$-saturated, i.e., $(\pi^{(d+1)})^{-1}(N_{d+1}(Y,T))=N_{d+1}(X,T)$.
\end{lem}

The following construction of $O$-diagram is shown by Veech in \cite{Veech70}. We need a version worked for both $S$ and $T$.

\begin{lem}[$O$-diagram]\label{O-diag}
Let $\pi: X\rightarrow Y$ be an extension between minimal systems both for $S$ and $T$. Then the following diagram holds both for $S$ and $T$.
\begin{equation*}
\xymatrix
{
X \ar[d]_{\pi}  &  X^* \ar[l]_{\sigma}\ar[d]^{\pi*} \\
Y &  Y^*\ar[l]^{\tau}
}
\end{equation*}
%\[\begin{tikzcd}
%X\ar[d,"\pi"left] & X^{*}\ar[l,"\sigma^{*}" above]\ar[d,"\pi^{*}"]\\
%Y&  Y^{*}\ar[l, "\tau^{*}"above]
%\end{tikzcd}\]
\begin{enumerate}
\item $\sigma^{*}$ and $\tau^{*}$ are almost one to one extensions;
\item $\pi^{*}$ is an open extension.
\end{enumerate}
\end{lem}
\begin{proof}
We recall Veech's construction of $X^{*}$ and $Y^{*}$.
Let $\pi^{-1}: Y\rightarrow 2^X, y\mapsto \pi^{-1}(y)$. Then  $\pi^{-1}$ is a u.s.c. map, and the set $Y_c$ of
continuous points of $\pi^{-1}$ is a dense $G_\d$ subset of $Y$.
Let $$\widetilde{Y}=\overline{\{\pi^{-1}(y): y\in Y\}}\ \text{and}\ Y^*=\overline{\{\pi^{-1}(y): y\in Y_c\}},$$
where the closure is taken in $2^X$. It is obvious that $Y^*\subseteq \widetilde{Y}\subseteq 2^X$. Note that
for each $A\in \widetilde{Y}$, there is some $y\in Y$ such that $A\subseteq \pi^{-1}(y)$, and hence $A\mapsto y$
define a map $\tau: \widetilde{Y}\rightarrow Y$. It is easy to verify that $\tau: \widetilde{Y}\rightarrow Y$
is a factor map both for $S$ and $T$. One can show that if $(Y,T)$ (resp. $(Y,S)$) is minimal then $(Y^*,T)$ (resp. $(Y^{*},S)$) is a minimal system and it is the unique minimal subsystem in $(\widetilde{Y},T)$ (resp. $(\widetilde{Y}, S)$), and $\tau: Y^*\rightarrow Y$ is an almost one to one extension such that $\tau^{-1}(y)=\{\pi^{-1}(y)\}$ for all $y\in Y_c$. It can be proved that $X^*=\{(\tilde x,\tilde y)\in X \times Y^* : \tilde x \in \tilde y\}$ is a minimal subset of $X\times Y^*$. Let $\sigma: X^*\rightarrow X$ and $\pi^*: X^*\rightarrow Y^*$ be the projections. One can show that $\sigma$ is almost one to one and $\pi^*$ is open.
See \cite[Subsection 2.3]{V77} for details. Note that the constructions of $Y^{*}$ and $X^{*}$ are independent of $S$ and $T$. Thus the commuting diagram holds both for $S$ and $T$.
\end{proof}

Let $(X,S)$,  $(X,T)$ be  minimal systems with $S$ and $T$ being commutative. Then by Theorem\ref{thm-RPd}, for any $d\in\mathbb{N}$,  we show that ${\RP}^{[d]}(X, S)={\RP}^{[d]}(X, T)$ and hence ${\RP}^{[\infty]}(X, S)={\RP}^{[\infty]}(X, T)$.
Thus both $(X, S)$ and $(X, T)$ have the same increasing sequence of pro-nilfactors:
$$\{pt\}= X_0  \longleftarrow  X_1 \longleftarrow \cdots  \longleftarrow X_n  \longleftarrow X_{n+1}
 \longleftarrow \cdots \longleftarrow X_\infty \longleftarrow \cdots \longleftarrow X,$$
where $X_d=X/\RP^{[d]}(X)$ is the the
maximal $d$-step pro-nilfactor, $d\in \N\cup\{\infty\}$. Thus by Theorem \ref{thm-GHSWY} and Lemma \ref{O-diag}, we have the following result.

\begin{thm}
Let $(X,S)$,  $(X,T)$ be  minimal systems with $S$ and $T$ being commutative, and $\pi:X\rightarrow X_\infty$ be the factor map from $X$ to its maximal $\infty$-step pro-nilfactor $X_\infty$ both for  $S$ and $T$. Then there are minimal system $X^*$ and $X_\infty^*$ which are almost one to one
extensions of $X$ and $X_\infty$ respectively, and a commuting diagram below such that $\pi^*: X^*\rightarrow X^*_\infty$ is open, and $X_\infty^*$ is a $d$-step topological characteristic factor of $X^*$  both for $S$ and $T$, for all $d\geq 2$.
\[
\begin{CD}
X @<{\varsigma^*}<< X^*\\
@VV{\pi}V      @VV{\pi^*}V\\
X_\infty @<{\varsigma}<< X_\infty^*
\end{CD}
\]
\end{thm}

\begin{thm}\label{N_d for commutative}
Let $(X,S)$,  $(X,T)$ be  minimal systems. If $S$ and $T$ commutes, then $N_{d}(X, S)=N_{d}(X, T)$ for each $d\in\mathbb{N}$.
\end{thm}

\begin{proof}
%Clearly, it is true for $d=2$. For $k\in\mathbb{N}\cup \{\infty\}$, note that ${\bf RP}^{[k]}(X, S)={\bf RP}^{[k]}(X,T)={\bf RP}^{[k]}(X,\langle S,T\rangle)$. We denote this common set by ${\bf RP}^{[k]}$ and let $X_{k}=X/{\bf RP}^{[k]}$. Note that $(X_{\infty}, \langle S,T\rangle)=\underset{\longleftarrow}{\lim} (X_{k},  \langle S,T\rangle)_{k\in\N} $. We first claim that $N_{d}(X_{\infty}, S)=N_{d}(X_{\infty}, T)$. For this, it suffices to show that $N_{d}(X_{k}, S)=N_{d}(X_{k}, T)$ for each $k\in\N$. For this,  we fix $k\in\N$. Recall that it was shown in \cite{FK05} that $(X_{k},\mathcal{X}_{k}, \mu, \langle S,T\rangle)$ is measurably isomorphic to an inverse of $k$-step nilsystems $\underset{\longleftarrow}{\lim}  (Y_{i}, \mathcal{Y},\mu_i, \langle S,T\rangle)$. Further, it was shown in \cite{HKM} that this measurably isomorphic can be identified with a topological isomorphism (also see Section 5). Thus $(X_{k}, \langle S,T\rangle)=\underset{\longleftarrow}{\lim}  (Y_{i}, \langle S,T\rangle)$ in the topological sense. By Corollary \ref{Nd for nil}, one has $N_{d}(Y_{i}, S)=N_{d}(Y_i, T)$ for each $i\in\N$ and hence $N_{d}(X_k, S)=N_{d}(X_k, T)$, since $N_{d}(X_k, R)=\underset{\longleftarrow}{\lim}N_{d}(Y_i, R)_{i} $, for $R=S,T$.
By Corollary \ref{Nd for nil} we have $N_{d}(X_{\infty}, S)=N_{d}(X_{\infty}, T)$.
Then by Lemma \ref{N_d for proximal extension}, we have $N_{d}(X^{*}_{\infty}, S)=N_{d}(X^{*}_{\infty}, T)$ since $\varsigma$ is proximal. This together with Lemma \ref{N_d for char factor} shows that \[N_{d}(X^{*},S)=(\pi^{*(d)})^{-1}N_{d}(X^{*}_{\infty}, S)=(\pi^{*(d)})^{-1}N_{d}(X^{*}_{\infty}, T)=N_{d}(X^{*}, T).\]
Then, we have $N_{d}(X,S)=\varsigma^{*(d)}(N_{d}(X^{*},S))=\varsigma^{*(d)}(N_{d}(X^{*},T))=N_{d}(X,T)$.
The proof is complete.
\end{proof}

Thus we complete the proof of the second assertion of Theorem \ref{cube and Nd}. As an immediate corollary, we have
\begin{cor}[{\cite[Theorem C]{GHSWY}}]
Let $(X,T)$ be a minimal system and $k\geq 2$. Then $(X, T^{k})$ is minimal if and only if $N_{d}(X, T)=N_{d}(X, T^{k})$ for each $d\in\mathbb{N}$.
\end{cor}

\section{Pro-nilfactors of $\mathbb{R}$-flows and Topological characteristic factors}\label{section-TCF}

In this section, we discuss the topological characteristic factors of $\R$-flows. To begin, we offer an overview of the pro-nilfactors associated with $\R$-flows.

\subsection{Pro-nilfactors of $\mathbb{R}$-flows}\label{subsection-nilflow}

Let $(X,\{T^{t}\}_{t\in \R})$ be a minimal $\R$-flow. By Theorem \ref{minimality} there is some countable subset $E$ of $\R$ such that $(X,T^t)$ is minimal for all $t\in \R\setminus E$. By Theorem \ref{RPd}, $$\RP^{[d]}(X,\mathbb{R})=\RP^{[d]}(X,T^{t}), \forall t\in \R\setminus E.$$
Thus we can define
$$X_d=X/\RP^{[d]}(X,\R)=X/ \RP^{[d]}(X,T^t), \ t\in \R\setminus E.$$

In the discrete setting, by Theorem \ref{thm0} $(X_{d}, T^{t})$ is a pro-nilsystem for any $t\in\R\setminus E$. The analogous result for flows holds as well. In other words, Let $(X, \{T^t\}_{t\in \R})$ be a minimal $\R$-flow and $d\geq 1$. Then $(X_d=X/{\bf RP}^{[d]}(\mathbb{R}), \{T^t\}_{t\in \R})$ is isomorphic to a $d$-step pro-nilflow.
To see this,  one approach is to employ \cite[Proposition 3.8]{Zie1}, which endows $(X_d,\{T^t\}_{t\in \R})$ with the structure of a $d$-step pro-nilflow in the measurable sense. Subsequently, the techniques developed in \cite{HKM} can be utilized to establish the isomorphism between $(X_d,\{T^t\}_{t\in \R})$ and a $d$-step pro-nilflow. Alternatively, the structure theory of nilspaces, as presented in \cite[Theorem 7.5]{GGY18}, provides another route to this conclusion. It is worth noting that the results in \cite{GGY18} are applicable to some more general groups.

\subsection{Topological characteristic factors for $(X,\{T^{t}\}_{t\in \R})$}\

Similar to the discrete case in Subsection \ref{subsection-chara}, we define the topological characteristic factors of an $\mathbb{R}$-flow as follows.

\begin{defn}
Let $\a_1,\ldots,\a_d$ be distinct nonzero real numbers. Given a factor map $\pi: (X,\{T^{t}\}_{t\in\R})\rightarrow (Y,\{T^{t}\}_{t\in\R})$ and $d\ge 2$, the $\R$-flow $(Y, \{T^{t}\}_{t\in\R})$ is said to be a {\em  topological characteristic factor of  $(X,\{T^{t}\}_{t\in\R})$} with respect to the family $\{\a_1t, \a_2t,\ldots, \a_dt\}$, if there exists a dense $G_\d$ subset $\Omega$ of $X$ such that for each $x\in \Omega$ the orbit closure $$\overline{\O}((x, \ldots,x), \{T^{\a_1t}\times\cdots\times T^{\a_dt}\}_{t\in \R})=\overline{\{ (T^{\a_1t}x,\cdots, T^{\a_dt}x):t\in \R\}}$$
is $\pi^{(d)}=\pi\times \cdots \times \pi$ ($d$-times) saturated.
\end{defn}

Similar to Theorem \ref{thm-GHSWY}, we have Theorem \ref{thm-Real}, which is restated as follows:
\begin{thm}%\label{thm-Real}
Let $(X,\{T^{t}\}_{t\in\R})$ be a minimal $\R$-flow and $d\in \N$. Let $\pi:X\rightarrow X_d$ be the factor map from $X$ to its maximal $d$-step pro-nilfactor $X_d$. Then there are minimal $\R$-flows $X^*$ and $X_d^*$ which are almost one to one
extensions of $X$ and $X_d$ respectively, and a commuting diagram below such that $\pi^*: X^*\rightarrow X^*_d$ is open, and $X_d^*$ is a topological characteristic factor of $X^*$ with respect to $\{\a_1t, \a_2t,\ldots, \a_dt\}$, where $\a_1,\ldots, \a_d$ are distinct nonzero real numbers.
\[
\begin{CD}
X @<{\varsigma^*}<< X^*\\
@VV{\pi}V      @VV{\pi^*}V\\
X_d @<{\varsigma}<< X_d^*
\end{CD}
\]

If in addition $\pi$ is open, then $X^*=X$, $X_d^*=X_d$ and $\pi^*=\pi$.
\end{thm}

In fact, we will show the following result, which combing with O-diagram implies Theorem \ref{thm-Real}.

%Now we are going to show the main result in this section, which implies Theorem \ref{thm-Real}.
\begin{thm}\label{thm-real-1}
Let $\pi:(X, \{T^{t}\}_{t\in\R})\rightarrow (Y,  \{T^{t}\}_{t\in\R})$ be an extension of minimal $\R$-flows and $d\in\N$. If $\pi$ is open and $X_{d}$ is a factor of $Y$, then $Y$ is a characteristic factor of $X$ with respect to $\{\a_1t, \a_2t,\ldots, \a_dt\}$, where $\a_1,\ldots, \a_d$ are distinct nonzero real numbers.
\end{thm}

In the discrete case, proving Theorem \ref{thm-GHSWY} relies on several results from \cite{BHK05} and \cite{HSY16}. However, extending these results to the continuous case for flows does not appear to be a straightforward consequence of the discrete case. We will shortly present the analogous result for flows and provide its proof in the appendix. Subsequently, we will present the proof of Theorem \ref{thm-real-1}. Our approach to prove Theorem \ref{thm-real-1} differs from the proof of Theorem \ref{thm-GHSWY} as presented in \cite{GHSWY}.

\medskip

First result we need is about nilfunctions. In \cite{BHK05}, one of the main result is stated as follows. Let $(X,\mu, T)$ be an ergodic system, let $f\in L^\infty(\mu)$ and let $d>1$ be an integer. The sequence $\{\int_{X} f(x)\cdot f(T^{n}x)\cdot\ldots\cdot f(T^{dn}x)d\mu(x)\}_{n\in\mathbb{Z}}$ is the sum of a sequence tending to zero in uniform density and a $d$-step nilsequence. We need a similar result for continuous flows.

\begin{defn}
Let $d\geq 1$ be an integer and $X=G/\Gamma$ be a $d$-step nilmanifold. Let $\phi$ be a continuous function on $X$. Let $\{a^{t}\}_{t\in\mathbb{R}}$ be an one-parameter subgroup of $G$ and $x\in X$. The function $\phi(a^{t}\cdot x)$ is called a {\it  basic $d$-step nilfunction}.  A {\it $d$-step nilfunction} is a uniform limit of basic $d$-step nilfunctions.
\end{defn}

\begin{defn}
Let $\phi: \mathbb{R}\rightarrow\mathbb{R}$ be a bounded continuous function. We say $\{\phi(t)\}_{t\in\mathbb{R}}$ tends to zero in uniform density, and we write ${\rm UD\textendash Lim} \phi(t)=0$, if
\[ \lim_{\rho\rightarrow \infty}\sup_{\sigma\in \mathbb{R}}\frac{1}{\rho}\int_{\sigma}^{\sigma+\rho} |\phi(t)|dt=0.\]
\end{defn}

\begin{thm}\label{nilfunction}
 Let  $\{T^{t}\}_{t\in\mathbb{R}}$ be an ergodic flow on $(X,\mathcal{X},\mu)$. Let $f\in L^{\infty}(\mu)$ and $\a_1,\ldots,\a_k$ be distinct nonzero real numbers. Then the function
 \[I_{f}(k,t):=\int_{X} f(x)\cdot f(T^{\a_1 t}x)\cdot\ldots\cdot f(T^{\a_k t}x)d\mu(x)\]
 is the sum of a function tending to zero in uniform density and $k$-step nilfunction.
\end{thm}

Theorem \ref{nilfunction} will be used to establish the following Theorem \ref{char of RP^{d}}, which is is essential for proving Theorem \ref{thm-real-1}.

\begin{defn}
\begin{enumerate}
\item[(1)] We say that $S\subset \mathbb{R}$ is a set of {\em $d$-Poincar\'{e} recurrence} for $\mathbb{R}$  if every measure preserving $\R$-flow $(X,\mathcal{X},\mu, \{T^{t}\}_{t\in\R})$ and for every $A\in\mathcal{X}$ with $\mu(A)>0$ and distinct $\a_1,\ldots,\a_d\in \R$, there is some $t\in S$ such that
\[ \mu(A\cap T^{-\a_1 t}A\cap \cdots\cap T^{-\a_d t}A)>0.\]
\item[(2)] We say that $S\subset \mathbb{R}$ is a set of {\em $d$-Birkhoff recurrence} for $\mathbb{R}$  if every minimal $\R$-flow $(X, \{T^{t}\}_{t\in\R})$ and for every nonempty open set $U\subset X$ and distinct $\a_1,\ldots,\a_k\in \R$, there is some $t\in S$ such that
\[ U\cap T^{-\a_1 t}U\cap \cdots\cap T^{-\a_d t}U\neq\emptyset.\]
\end{enumerate}
 We use $\mathcal{F}_{Poi_{d}}(\R)$ and $\mathcal{F}_{Bir_{d}}(\R)$ to denote the families consisting of all sets of $d$-Poincar\'{e} recurrence for $\mathbb{R}$ and all sets of $d$-Birkhoff recurrence for $\mathbb{R}$, respectively.
\end{defn}

Let $(X,\{T^{t}\}_{t\in\mathbb{R}})$ be an $\mathbb{R}$-flow. For $x\in X$ and $A\subset X$, define
\[ R(x,A):=\{t\in\mathbb{R}: T^{t}x\in A\}.\]

The following equivalent characterizations of regionally proximal relation of higher order are given in \cite{HSY16} for discrete systems. %However, the continuous version seems not directly followed from the discrete case. We provide the proof in the appendix.

\begin{thm}\label{char of RP^{d}}
Let $(X, \{T^{t}\}_{t\in\mathbb{R}})$ be a minimal $\R$-flow and $x,y\in X$. The following assertions are equivalent for $d\in\N$:
\begin{enumerate}
\item[(1)] $(x,y)\in {\bf RP}^{[d]}(X)$.
\item[(2)] $R(x, U)\in \mathcal{F}_{Poi_{d}}(\R)$ for every neighborhood $U$ of $y$.
\item[(3)] $R(x, U)\in \mathcal{F}_{Bir_{d}}(\R)$ for every neighborhood $U$ of $y$.
\end{enumerate}
\end{thm}

\subsection{Proof of Theorem \ref{thm-real-1}}

Let $(X, \{T^{t}\}_{t\in\R})$ be an $\R$-flow and $d\in\N$. For $A\subset X$, set
\[ \Delta_{d}(A)=\{(x,\ldots,x)\in X^{d}: x\in A\}.\]
 Let $\a_1,\ldots,\a_d$ be distinct nonzero real numbers and set $\vec{\a}=(\a_1,\ldots,\a_d)$. Set
\[ \tau_{\vec{\a}}^{t}=T^{\a_1 t}\times T^{\a_2 t}\times\cdots\times T^{\a_d t}, t\in\R,\]
and for $\beta\in \R$, set
\[ T^{\beta t}_{\Delta_d}=T^{\beta t}\times T^{\beta t}\times\cdots\times T^{\beta t}, t\in \R.\]
Define
\[ N_{\vec{\a}}(X)=\overline{\O}(\Delta_{d}(X),\{\tau_{\vec{\a}}^{t}\}_{t\in\R})=
\overline{\{(T^{\a_1 t}x,T^{\a_2 t}x,\ldots, T^{\a_d t}x): t\in\R, x\in X\}}.\]

The following result is shown in \cite[Proposition 4.4]{GKR19} for discrete-time systems. But it is not difficult to show this also holds for continuous flows from the proof there.
\begin{lem}\label{cts pt for F}
Let $(X, \{T^{t}\}_{t\in\R})$ be a minimal $\R$-flow and $d\in\N$. For any distinct nonzero real numbers $\a_1,\ldots,\a_d$, there is residual set $X_0\subset X$ such that the map
\[ F: N_{\vec{\a}}(X)\rightarrow 2^{N_{\vec{\a}}(X)}, (x_1,\ldots, x_d)\mapsto \overline{\O}((x_1,\ldots,x_d), \{\tau_{\vec{\a}}^{t}\}_{t\in \R})\]
is continuous at $x^{(d)}=(x,x,\ldots, x)$ for any $x\in X_0$.
\end{lem}

The following result is obtained in the same way as in \cite[Proposition 1.55]{G03}.
\begin{thm}
Let $(X, \{T^{t}\}_{t\in\R})$ be a minimal $\R$-flow and $d\in\N$. For any $\beta\in\R\setminus\{0\}$ and any distinct nonzero real numbers $\a_1,\ldots,\a_d$, the $\R^{2}$-system $\left( N_{\vec{\a}}(X), \{\tau_{\vec{\a}}^{s}, T_{\Delta_d}^{\beta t}\}_{(s,t)\in\R^{2}}\right)$ is minimal and the set of $\{\tau_{\vec{\a}}^{s}\}_{s\in \R}$-minimal points is dense in $N_{\vec{\a}}(X)$.
\end{thm}

\medskip

Now we are ready to prove \ref{thm-real-1}.

\begin{proof}[Proof of Theorem \ref{thm-real-1}]
Let
\[ F: N_{\vec{\a}}(X)\rightarrow 2^{N_{\vec{\a}}(X)}, (x_1,\ldots, x_d)\mapsto \overline{\O}((x_1,\ldots,x_d), \{\tau_{\vec{\a}}^{t}\}_{t\in \R}).\]
Then there is a residual set $X_0$ of $X$ such that for each $x\in X_0$, $x^{(d)}:=(x,\ldots,x)$ is a continuous point of $F$. We will show that for each $x\in X_0$, the orbit closure $L_{x}=\overline{\O}(x^{(d)}, \{\tau_{\vec{\a}}^{t}\}_{t\in \R})$ is $\pi^{(d)}$-saturated. This implies that $Y$ is a characteristic factor of $X$ with respect to $\{\a_1t, \a_2t,\ldots, \a_dt\}$.

\noindent{\bf Claim}. Let $x\in X_0$ and $j\in\{1,\ldots,d\}$. Then for each $y\in X$ with $\pi(y)=\pi(x)$, one has
\[(x^{(j-1)},y, x^{(d-j)})=(x,\ldots,x,\underset{j^{th}}{y},x,\ldots,x)\in L_{x}.\]

\begin{proof}[Proof of Claim]
Let $\rho$ be the metric on $X$ and $\rho_{d}$ be the metric on $X^{d}$ given by
\[ \rho_{d}({\bf x}, {\bf y})=\max_{1\leq j\leq d}\rho(x_j, y_j),\]
where ${\bf x}=(x_1,\ldots,x_d)$ and ${\bf y}=(y_1,\ldots,y_d)$.

For a given $\epsilon>0$, it follows from the continuity of $F$ at $x^{(d)}$ that there is some $\delta\in(0, \epsilon/3)$ such that whenever $\rho_{d}({\bf x}, x^{(d)})<\delta$ with ${\bf x}=(x_1,\ldots,x_d)\in N_{\vec{\a}}(X)$, one has that
\begin{equation}\label{6eq1}
H_{\rho} \left(\overline{\O}( {\bf x}, \{\tau_{\vec{\a}}^{t}\}_{t\in \R}), \overline{\O}(x^{(d)}, \{\tau_{\vec{\a}}^{t}\}_{t\in \R}) \right) <\frac{\epsilon}{2},
\end{equation}
where $H_{\rho}$ is the Hausdorff metric on $2^{N_{\vec{\a}}(X)}$.

Let
\[ S_{j}^{t}= T^{(\a_1-\a_j)t}\times\cdots\times T^{(\a_{j-1}-\a_j)t}\times id\times T^{(\a_{j+1}-\a_j)t}\times\cdots\times T^{(\a_d-\a_j)t}, t\in \R. \]

Note that
\[ \langle \tau_{\vec{\a}}^{s}, T^{t}_{\Delta_{d}}: (s, t)\in \R^{2}\rangle= \langle S_{j}^{s}, T^{t}_{\Delta_{d}}: (s, t)\in \R^{2}\rangle.\]
The system $(N_{\vec{\a}}(X), \{S_{j}^{s}, T^{t}_{\Delta_d}\}_{(s,t)\in \R^{2}})$ is minimal and the set of $\{S_j^{s}\}_{s\in\R}$-minimal points is dense in $N_{\vec{\a}}(X)$. Take an  $\{S_j^{s}\}_{s\in\R}$-minimal point ${\bf z}=(z_1,\ldots, z_d)\in N_{\vec{\a}}(X)$ with
\[ \rho_{d}({\bf z}, x^{(d)})<\frac{\delta}{2}.\]
In particular, $\rho(z_j, x)\leq \rho_{d}({\bf z}, x^{(d)})<\frac{\delta}{2}$. Since $\pi$ is open, we may assume that there is some $y_j\in X$ such that
\[ \rho(y_j, y)<\frac{\epsilon}{6} \text{ and } \pi(y_j)=\pi(z_j).\]
Set
\[ Z=\overline{\O}({\bf z}, \{S_j^{t}\}_{t\in\R}).\]
Then $(Z,  \{S_j^{t}\}_{t\in\R})$ is a minimal $\R$-flow. Note that by the $j^{th}$ component of $S_j$ is the identity, one has $x_j=z_j$ for any ${\bf x}=(x_1,\ldots,x_d)\in Z$. Let
\[ U=Z\cap \left(B_{\delta/2}(z_1)\times B_{\delta/2}(z_2)\times\cdots\times B_{\delta/2}(z_d) \right)\]
 be an open subset of $Z$.

 \medskip
 Since $X_d$ is a factor of $Y$ and $\pi(z_j)=\pi(y_j)$, one has $(z_j,y_j)\in{\bf RP}^{[d]}(X, \{T^{t}\}_{t\in\R})$. Thus $R(z_j, V)\in \mathcal{F}_{Bir_d}(\R)$ for each neighborhood $V$ of $y_j$. Clearly, $ {\bf RP}^{[d]}(X, \{T^{t}\}_{t\in\R}) ={\bf RP}^{[d]}(X, \{T^{\beta t}\}_{t\in\R})$ for any $\beta\neq 0$. Considering $\R$-flows $(X, \{T^{\a_j t}\}_{t\in\R})$ and $(Z, \{S_{j}^{t}\}_{t\in\R})$,  it follows from the definition of $ \mathcal{F}_{Bir_d}(\R)$ that  for any distinct $\beta_1,\ldots, \beta_d\in\R$, there is some $t>0$ such that
 \[ T^{\a_j t}z_j\in B_{\epsilon/6}(y_j)\]
 and
 \[ U\cap S_j^{-\beta_1t}U\cap S_j^{-\beta_2t}U\cap \cdots\cap  S_j^{-\beta_dt}U \neq \emptyset. \]
Particularly, we take $\beta_{j}=0$ and
\[ \beta_{i}=\frac{\a_i}{\a_i-\a_j}, \text{ for } i=1,\ldots, j-1, j+1,\ldots, d.\]
Take ${\bf w}=(w_1,\ldots,w_d)\in U\cap S_j^{-\beta_1t}U\cap S_j^{-\beta_2t}U\cap \cdots\cap  S_j^{-\beta_dt}U$. Then one has $S_{j}^{\beta_i t}{\bf w}\in U$. That is
\[ T^{\a_i t}w_i=T^{(\a_i-\a_j)\beta_{i}t}w_i\in B_{\delta/2}(z_i), \ \ \forall i=1,\ldots,j-1,j+1,\ldots, d.\]
Recall that $z_j=w_j$ and  $T^{\a_j t}z_j\in B_{\epsilon/6}(y_j)$. Thus
\[\rho_{d}(\tau_{\vec{\a}}^{t}{\bf w}, (z_1,\ldots,z_{j-1}, y_j, z_{j+1},\ldots,z_d))<\max\{\delta/2, \epsilon/6\}\leq \epsilon/6. \]
Since $\rho_{d}({\bf z}, x^{(d)})<\delta/2<\epsilon/6$ and $\rho(y_j, y)<\epsilon/6$, we have that
\[\rho_{d}( (z_1,\ldots,z_{j-1}, y_j, z_{j+1},\ldots,z_d), (x^{(j-1)}, y, x^{(d-j)}))< \epsilon/6. \]
Thus
\begin{equation}\label{6eq2}
\rho_{d}\left((x^{(j-1)}, y, x^{(d-j)}), \overline{\O}({\bf w},\{\tau_{\vec{\a}}^{t}\}_{t\in\R}) \right)<\epsilon/3.
\end{equation}
Since ${\bf w}\in U$ and $\rho_{d}({\bf w}, {\bf z})<\delta/2$, we have
\[\rho_{d}({\bf w}, x^{(d)})<\rho_{d}({\bf w},{\bf z})+\rho_{d}({\bf z}, x^{(d)})<\frac{\delta}{2}+\frac{\delta}{2}=\delta.\]
By (\ref{6eq1}), we have
\[H_{\rho} \left(\overline{\O}( {\bf w}, \{\tau_{\vec{\a}}^{t}\}_{t\in \R}), \overline{\O}(x^{(d)}, \{\tau_{\vec{\a}}^{t}\}_{t\in \R}) \right) <\frac{\epsilon}{2}.\]
This together with (\ref{6eq2}) implies that
\[ \rho_{d}\left((x^{(j-1)}, y, x^{(d-j)}), \overline{\O}({\bf x},\{\tau_{\vec{\a}}^{t}\}_{t\in\R}) \right)<\frac{\epsilon}{3}+\frac{\epsilon}{2}<\epsilon.\]
Since $\epsilon$ is arbitrary, we have
\[(x^{(j-1)},y, x^{(d-j)})=(x,\ldots,x,\underset{j^{th}}{y},x,\ldots,x)\in \overline{\O}({\bf x},\{\tau_{\vec{\a}}^{t}\}_{t\in\R})\]
and the claim is followed.
\end{proof}

\medskip

Now we will use Claim to show that for each $x\in X_0$ the orbit closure $L_x=\overline{\O}(x^{(d)}, \{\tau_{\vec{\a}}^{t}\}_{t\in\R})$ is $\pi^{(d)}$-saturated, i.e. $Y$ is $d$-step topological characteristic factor of $X$ with respect to $\{\a_1, \a_2,\ldots, \a_d\}$.

For $j\in \{1,2,\ldots, d\}$, let
$${\bf z}=(z_1,z_2,\ldots, z_d)\in L_x=\overline{\O}(x^{(d)}, \{\tau_{\vec{\a}}^{t}\}_{t\in\R}).$$
We show that $${\bf z'}=(z_1,z_2,\ldots, z_{j-1}, z_j',z_{j+1},\ldots, z_d)\in L_x,$$
where $\pi(z_j')=\pi(z_j)$.

Since ${\bf z}=(z_1,z_2,\ldots, z_d)\in L_x,$
there is some sequence $\{t_i\}_{i\in \N}$ such that
$$\tau_{\vec{\a}}^{t_i}x^{(d)}\xrightarrow{i\to\infty} {\bf z}.$$
By Claim we have that
$$\{x^{(j-1)}\}\times \pi^{-1}(\pi(x))\times \{x^{(d-j)}\}\subset L_x.$$
Since $\pi$ is open, it follows that
\begin{equation*}
  \begin{split}
   &\quad  \tau_{\vec{\a}}^{t_i}\Big (\{x^{(j-1)}\}\times \pi^{-1}(\pi(x))\times \{x^{(d-j)}\}\Big ) \\
    & =\{T^{\a_1 t_i}x\}\times \{T^{\a_2 t_i}x\}\times \ldots \times \{T^{\a_{j-1}t_i}x\}\times \pi^{-1}(\pi(T^{\a_j t_i}x))\times \{T^{\a_{j+1}t_i}x\}\times \ldots \times \{T^{\a_d t_i}x\}\\
    & \xrightarrow{i\to\infty} \{z_1\}\times \{z_2\}\times \ldots \times \{z_{j-1}\}\times \pi^{-1}(\pi(z_j))\times \{z_{j+1}\}\times \ldots \times \{z_d\}\\
    & \subset \overline{\O}(L_x, \{\tau_{\vec{\a}}^{t}\}_{t\in\R})=L_x.
   \end{split}
\end{equation*}

To sum up, we have that if
${\bf z}=(z_1,z_2,\ldots, z_d)\in L_x=\overline{\O}(x^{(d)}, \{\tau_{\vec{\a}}^{t}\}_{t\in\R}),$ then for each $j\in \{1,2,\ldots, d\}$,
$$(z_1,z_2,\ldots, z_{j-1}, z_j',z_{j+1},\ldots, z_d)\in L_x,$$
where $\pi(z_j')=\pi(z_j)$.
Thus,  $(z_1,z_2,\ldots, z_d)\in L_x$ if and only if $(z_1',z_2',\ldots, z_d')\in L_x$ whenever for all $j\in \{1,2,\ldots, d\}$, $\pi(z_i)=\pi(z_i')$. That is,
$L_x$ is $\pi^{(d)}$-saturated.
\end{proof}

For further applications, we give some equivalent forms of Theorem \ref{thm-Real}.  The following concept of weak topological characteristic factor is introduced in \cite{QXYY}. 
\begin{defn}
Let $\a_1,\ldots,\a_d$ be distinct nonzero real numbers. Given a factor map $\pi: (X,\{T^{t}\}_{t\in\R})\rightarrow (Y,\{T^{t}\}_{t\in\R})$ and $d\ge 2$, the $\R$-flow $(Y, \{T^{t}\}_{t\in\R})$ is said to be a {\em  weak topological characteristic factor of  $(X,\{T^{t}\}_{t\in\R})$} with respect to the family $\{\a_1t, \a_2t,\ldots, \a_dt\}$, if there exists a dense $G_\d$ subset $\Omega$ of $X$ such that for each $x\in \Omega$, one has
$$(\pi^{-1}(\pi x))^{d}\subset\overline{\{ (T^{\a_1t}x,\cdots, T^{\a_dt}x):t\in \R\}}.$$
\end{defn}

The following equivalent characterizations is given in \cite[Proposition 2.16]{QXYY} for discrete discrete systems. But the proof for  $\R$-flows are the same as in  and we omit the proof. 
\begin{prop}
Suppose that there exists a commuting diagram of minimal $\mathbb{R}$-flows given by $(X,\{T^{t}\}_{t\in\R}),(X^{*},\{T^{t}\}_{t\in\R})$ and $(Y,\{T^{t}\}_{t\in\R}),(Y^{*},\{T^{t}\}_{t\in\R})$ below such that $\varsigma$ and $\varsigma^{*}$ are almost one to one, while $\pi^{*}$ is open. 
\[
\begin{CD}
X @<{\varsigma^*}<< X^*\\
@VV{\pi}V      @VV{\pi^*}V\\
Y @<{\varsigma}<< Y^*
\end{CD}
\]
Let $d\in\mathbb{N}$ and $\a_1,\a_2,\ldots,\a_d$ be distinct nonzero real numbers. Then the following statements are equivalent.
\begin{enumerate}
\item[(1)] $Y^{*}$ is a topological characteristic factor of $X^{*}$ with respect to $\{\a_1 t, \a_2 t,\ldots, \a_d t\}$.
\item[(2)] For any open subsets $V_0, V_1,\ldots, V_{d}$ of $X$ with $\bigcap_{i=0}^{d}\pi(V_i)$ having nonempty interior, there is some $t\in \mathbb{R}$ such that $V_0\cap T^{-\a_1 t}V_1\cap\cdots\cap T^{-a_d t}V_d\neq\emptyset$.
\item[(3)] $Y$ is a weak topological characteristic factor of $X$ with respect to $\{\a_1 t, \a_2 t,\ldots, \a_d t\}$.
\item[(4)] $Y^{*}$ is a weak topological characteristic factor of $X^{*}$ with respect to $\{\a_1 t, \a_2 t,\ldots, \a_d t\}$.
\item[(5)] For any open subsets $U_0, U_1,\ldots, U_{d}$ of $X^{*}$ with $\bigcap_{i=0}^{d}\pi^{*}(U_i)$ having nonempty interior, there is some $t\in \mathbb{R}$ such that $U_0\cap T^{-\a_1 t}U_1\cap\cdots\cap T^{-a_d t}U_d\neq\emptyset$.
\end{enumerate}

\end{prop}

In particular, one has the following form. In the case of discrete system, this is given in \cite{YY}.  
\begin{thm}\label{equiv-thm-Real}
Let $(X,\{T^{t}\}_{t\in\R})$ be a minimal $\R$-flow and $d\in \N$. Let $\pi:X\rightarrow X_d$ be the factor map from $X$ to its maximal $d$-step pro-nilfactor $X_d$. Then for any distinct nonzero real numbers $\a_1,\a_2,\ldots,\a_d$, there is a dense $G_{\delta}$ set $\Omega$ in $X$ such that for any $x\in \Omega$, one has
\[(\pi^{-1}(\pi x))^{d}\subset \overline{\{(T^{\a_1 t}x,T^{\a_2 t}x,\ldots, T^{\a_d t}x): t\in\mathbb{R}\}}.\]
\end{thm}

\subsection{Polynomial saturation theorems}\
\medskip

A family of polynomials $\{p_1(t),\ldots, p_d(t)\}\subset \R[t]$ is said to be {\em essentially distinct} if $p_j-p_i$ is non-constant for all $i\neq j\in \{1,2,\ldots,d\}$.

%Let $(X,T)$ be a system and $d\in \N$. Let $\A=\{p_1,\ldots, p_d\}$ be a set of integral polynomials and define $$\O_\A((x_1,\ldots,x_d))=\O_{\{p_1,\ldots,p_d\}}((x_1,\ldots,x_d))=\{(T^{p_1(n)}x_1, \ldots, T^{p_d(n)}x_d): n\in \Z \},$$ and $$\overline{\O}_\A((x_1,\ldots,x_d))=\overline{\{(T^{p_1(n)}x_1, \ldots, T^{p_d(n)}x_d): n\in \Z \}}.$$

%We have the following polynomial version of Theorem \ref{thm-GHSWY}.

%\begin{thm}\cite{HSY22-1}\label{thm-poly-sat1}
%Let $(X,T)$ be a minimal system, and $\pi:X\rightarrow X_\infty$ be the factor map from $X$ to its maximal $\infty$-step pro-nilfactor $X_\infty$. Then there are minimal system $X^*$ and $X_\infty^*$ which are almost one to one extensions of $X$ and $X_\infty$ respectively, an open factor map $\pi^*$ and a commuting diagram below
%\[ \begin{CD} X @<{\varsigma^*}<< X^*\\@VV{\pi}V      @VV{\pi^*}V\\X_\infty @<{\varsigma}<< X_\infty^*\end{CD}\]
%such that there is a $T$-invariant residual subset $X^*_0$ of $X^*$ having the following property:  for all $x\in X^*_0$, all sets of essentially distinct non-constant integral polynomials $\A=\{p_1,\ldots, p_d\}$ and $d\in \N$, $\overline{\O}_\A(x)=\overline{\{(T^{p_1(n)}x, \ldots, T^{p_d(n)}x): n\in \Z \}}$ is ${\pi^*}^{(d)}$-saturated, that is, $$\overline{\O}_\A(x^{\otimes d})=({\pi^*}^{(d)})^{-1}\Big({\pi^*}^{(d)}(\overline{\O}_\A(x^{\otimes d}))\Big).$$

%If in addition $\pi$ is open, then $X^*=X$, $X_\infty^*=X_\infty$ and $\pi^*=\pi$.
%\end{thm}

Based on Theorem \ref{thm-Real} and using PET induction, one can give the polynomial saturation theorem for flows. The proof of Theorem \ref{thm-Real-poly} closely mirrors the proof of \cite[Theorem 3.6]{HSY22-1}, and thus, we choose to omit it for brevity.

\begin{defn}
Let $\a_1,\ldots,\a_d$ be distinct nonzero real numbers. Given a factor map $\pi: (X,\{T^{t}\}_{t\in\R})\rightarrow (Y,\{T^{t}\}_{t\in\R})$ and $d\ge 2$, the $\R$-flow $(Y, \{T^{t}\}_{t\in\R})$ is said to be a {\em  topological characteristic factor of  $(X,\{T^{t}\}_{t\in\R})$} with respect to $\A=\{p_1(t), p_2(t),\ldots, p_d(t)\}\subset \R[t]$, if there exists a dense $G_\d$ subset $\Omega$ of $X$ such that for each $x\in \Omega$
$$\overline{\O}_\A((x_1,\ldots,x_d))=\overline{\{(T^{p_1(t)}x, \ldots, T^{p_d(t)}x): t\in \R \}}.$$ is $\pi^{(d)}$ saturated.
\end{defn}

\begin{thm}\label{thm-Real-poly}
Let $(X,\{T^{t}\}_{t\in\R})$ be a minimal $\R$-flow and $d\in \N$. Let $\pi:X\rightarrow X_\infty$ be the factor map from $X$ to its maximal $\infty$-step pro-nilfactor $X_\infty$. Then there are minimal $\R$-flows $X^*$ and $X_\infty^*$ which are almost one to one
extensions of $X$ and $X_\infty$ respectively, and a commuting diagram below such that $\pi^*: X^*\rightarrow X^*_\infty$ is open, and $X_\infty^*$ is a topological characteristic factor of $X^*$ with respect to $\{p_1(t), p_2(t),\ldots, p_d(t)\}$, where $p_1,\ldots, p_d$ are essentially distinct nonconstant real polynomials.
\[
\begin{CD}
X @<{\varsigma^*}<< X^*\\
@VV{\pi}V      @VV{\pi^*}V\\
X_\infty @<{\varsigma}<< X_\infty^*
\end{CD}
\]
\end{thm}

Using Theorem \ref{thm-Real-poly}, one can prove the following result, which is the continuous version of results in \cite{Qiu}. Given the similarity in methodology to that employed in \cite{Qiu}, we choose to omit the detailed proof here.

\begin{thm}\label{polynomial orbit in R}
Let $(X, \mathbb{R})$ be a minimal $\mathbb{R}$-flow. Suppose that the family of polynomials $p_1,\ldots, p_d\in \R[t]$ is $\R$-independent, then there is some $x\in X$ such that
$$\overline{\{(T^{p_1(t)}x, \ldots, T^{p_d(t)}x): t\in \R \}}=X^d.$$

In particular, for any non-constant real polynomial $p\in \mathbb{R}[t]$, there is some $x\in X$ such that $\{T^{p(t)}x: t\in\mathbb{R}\}$ is dense in $X$.
\end{thm}

Given that every flow possesses ergodic invariant measures, an alternative approach to obtain Theorem \ref{polynomial orbit in R} is by using the following result presented in \cite{P11}.

\begin{thm}\cite[Theorem 1.2]{P11}\label{thm-potts}
Suppose $\{T^t\}_{t\in \R}$ is an ergodic measure preserving flow on
a Lebesgue space $(X,{\mathcal X}, \mu)$, the family of polynomials $p_1,\ldots, p_d\in \R[t]$ is
$\R$-independent, and $f_1,\ldots,f_d\in L^\infty(\mu)$. Then
as $R\to \infty$,
$$\frac{1}{R}\int_{0}^Rf_1(T^{p_1(t)})f_2(T^{p_2(t)})\cdots f_d(T^{p_d(t)})dt$$
converges in $L^2(\mu)$ to $\displaystyle \int_Xf_1d\mu\cdot\int_Xf_2d\mu\cdot \cdots\cdot \int_Xf_dd\mu.$
\end{thm}

\subsection{Topological characteristic factors for $\Z^{d}$-systems}

For $Z^d$-systems, the situation is different from the $\Z$-systems, and currently, we are still uncertain about how to approach it.

\begin{defn}
Given a factor map $\pi: (X,\langle T_1,\ldots, T_d\rangle)\rightarrow (Y,\langle T_1,\ldots, T_d\rangle)$ and $d\ge 2$,
the system $(Y, \langle T_1,\ldots, T_d\rangle)$ is said to be a {\em topological
characteristic factor of  $(X,\langle T_1,\ldots, T_d\rangle)$}, if there exists a dense $G_\d$ subset $\Omega$ of $X$ such that for each $x\in \Omega$ the orbit
closure $L_x=\overline{\O}((x, \ldots,x), T_1\times T_2\times \cdots \times T_d)$ is $\pi^{(d)}=\pi\times \cdots \times
\pi$ ($d$-times) saturated.
\end{defn}

\begin{conj}\label{conj-Zm}
Let $(X, \langle T_1,\ldots, T_d\rangle)$ be a system where $T_iT_j=T_jT_i$ for all $i,j\in \{1,2,\ldots, d\}$. If $T_1,\ldots, T_d$ and $T_iT^{-1}_j, i\neq j$ are all minimal on $X$. Then there are minimal system $X^*$ and $X_\infty^*$ which are almost one to one
extensions of $X$ and $X_\infty$ respectively, and a commuting diagram below such that $\pi^*: X^*\rightarrow X^*_\infty$ is open, and $(X_\infty^*, \langle T_1,\ldots, T_d\rangle)$ is a topological characteristic factor of $(X^*, \langle T_1,\ldots, T_d\rangle)$.
\[
\begin{CD}
X @<{\varsigma^*}<< X^*\\
@VV{\pi}V      @VV{\pi^*}V\\
X_\infty @<{\varsigma}<< X_\infty^*
\end{CD}
\]

In particular, for a distal system $X_{d-1}$ is the maximal topological characteristic factor of $(X,\langle T_1,\ldots, T_d\rangle)$.
\end{conj}

In Conjecture \ref{conj-Zm} the condition about $T_iT^{-1}_j, i\neq j$ being minimal is necessary.

\begin{example}
Let $(Y,S)$ be a weakly mixing minimal system and let $Z=\T$ and $R_1, R_2: Z\rightarrow Z$ be two minimal irrational rotations. Let $X=Y\times Z$ and let $T_i=S\times R_i$, $i=1,2$. Then $(X,T_1)$ and $(X,T_2)$ are minimal system and $T_1T_2=T_2T_1$. It is clear that $(X,T_1T^{-1}_2)$ is not minimal.

It is clear that $X_1(T_i)=Z, i=1,2$. Let $\pi\times \pi: X\times X\rightarrow Z\times Z$, $\left((y_1,z_1),(y_2,z_2)\right)\mapsto (z_1,z_2)$.
It is easy to check that for each $x\in X$, $\overline{\O}\left((x,x), T_1\times T_2\right)$ is not $\pi\times\pi$-saturated, and $X_1$ is not characteristic factor of $(X,\langle T_1, T_2\rangle)$.

\end{example}

\section{Suspension}\label{section-suspension}

In this section we study the suspension flows.
%\subsection{Some basics of suspension flow}
Let $(X,T)$ be a system.   The suspension flow of $(X,T)$ is defined as follows. Let $\tilde{X}=X\times\mathbb{R}/\sim$, where the equivalent relation is defined by
\[ (x,s)\sim (y,t) \Longleftrightarrow s-t\in\mathbb{Z} \text{ and } T^{s-t} (x)=y.\]
In particular, $(x,1)\sim (Tx, 0)$. Let $q: X\times \mathbb{R}\rightarrow \tilde{X}$ denote the quotient map.  In this way, there is an $\mathbb{R}$ action on $\tilde{X}$ by
\[ T^{t}\cdot [x, s]=[x, s+t],\]
where $[x, s]=q((x,s))$.

%Let $\rho$ be the metric on $\tilde{X}$ defined as
%\[ \rho([x,s],[y,t])=\rho_{X}(T^{[s]}x, T^{[t]}y)+|\{s\}-\{t\}|.\]
For $(x,s)\in X\times\mathbb{R}$, $[x,s]$ has a local neighborhood basis in $\tilde{X}$ consisting of all sets of the form
\[ [U, s,\varepsilon]:=q(U\times (s-\varepsilon, s+\varepsilon))\]
with $U$ a neighborhood of $x$ in $X$ and $0<\varepsilon\leq 1/2$.

\begin{lem}\cite[Theorem 5.12]{Vries93}\label{min Z vs R}
Let $(X,T)$ be a system and $(\tilde{X},\mathbb{R})$ be the suspension flow. Then $(X,T)$ is minimal if and only if $(\tilde{X},\mathbb{R})$ is minimal.
\end{lem}

\begin{lem} \label{RP for suspension}
 Let $(X,T)$ be a minimal system and $(\tilde{X},\mathbb{R})$ be its suspension flow. Then for any $d\in\mathbb{N}$,
  $$([x_1,s_1], [x_2,s_2])\in {\bf RP}^{[d]}(\tilde{X},\mathbb{R})   \Longleftrightarrow s_1-s_2\in\mathbb{Z} \text{ and } (T^{s_1-s_2}x_1,x_2)\in{\bf RP}^{[d]}(X,T) .$$
 \end{lem}

\begin{proof}
Suppose that  $(x_1,x_2)\in{\bf RP}^{[d]}(X,T) $  and $s\in[0,1)$. We show $([x_1,s], [x_2,s])\in {\bf RP}^{[d]}(\tilde{X},\mathbb{R})$. By definition, it suffices to show there is $z\in X$ such that  for any $\delta>0$, and neighborhoods $U_{1}, U_{2}, U_{3}$ of $x_1,x_2,z$ respectively,  there are  $y_1,y_2$ in $X$ with  $y_i\in U_i, i=1,2$ and   $\vec{n}=(n_1,\ldots,n_d)\in\mathbb{Z}^{d}$ such that for any ${\bf \epsilon}\in \{0,1\}^{d}\setminus\{{\bf 0}\}$,
\[T^{{\bf \epsilon}\cdot \vec{n}}[ y_i,s]=[y_i,s+{\bf \epsilon}\cdot \vec{n}]=[T^{{\bf \epsilon}\cdot \vec{n}}y_i,s] \in [U_3, s, \delta], i=1,2.  \]
But this can be done since $(x_1,x_2)\in{\bf RP}^{[d]}(X,T) $. Indeed, there is some $z\in X, y_1\in U_1,y_2\in U_2$ and $\vec{n}=(n_1,\ldots,n_d)\in\mathbb{Z}^{d}$ such that for any ${\bf \epsilon}\in \{0,1\}^{d}\setminus\{{\bf 0}\}$, $T^{{\bf \epsilon}\cdot \vec{n}} y_i\in U_3, i=1,2$.

Now suppose that $([x_1,s_1], [x_2,s_2])\in {\bf RP}^{[d]}(\tilde{X},\mathbb{R})$.

{\bf Claim 1}. $s_1-s_2\in\mathbb{Z}$.

Let $\pi: \tilde{X}\rightarrow \mathbb{R}/\mathbb{Z}, [x,s]\mapsto s+\mathbb{Z}$. It is clear that $\pi$ commutes the $\mathbb{R}$ action on $\tilde{X}$ and natural action on $\mathbb{R}/\mathbb{Z}$. Thus, $(\mathbb{R}/\mathbb{Z}, \mathbb{R})$ is a factor of $(\tilde{X}, \mathbb{R})$. Thus $(s_1+\mathbb{Z},s_2+\mathbb{Z})\in{\bf RP}^{[d]}(\mathbb{R}/\mathbb{Z}, \mathbb{R})$. So we have $s_1-s_2\in\mathbb{Z}$.

By Claim 1, we may assume that $s_1=s_2=s\in[0,1)$.

Given $\eta>0$, we take $\delta\in (0,1/4)$ be such that $\max\limits_{-d-3\leq i\leq d+3} \rho_{X}(T^{i}x, T^{i}y)<\eta$ whenever $\rho_{X}(x,y)<2\delta$. By definition, there is $z\in X, s'\in[0,1)$,  $[y_1,s_1'], [y_2,s_2']$ in $\tilde{X}$ with $s_i'\in[0,1)$, $[y_i,s_i']\in[B(x_i,\delta), s,\delta], i=1,2$ and   $\vec{t}=(t_1,\ldots,t_d)\in\mathbb{R}^{d}$ such that for any ${\bf \epsilon}\in \{0,1\}^{d}\setminus\{{\bf 0}\}$,
\[[y_i, s_i'+{\bf \epsilon}\cdot \vec{t}]  \in [B(z', \delta), s',\delta], i=1,2.  \]
Thus there are $z_{i, \epsilon}\in B(z,\delta)$ and $\delta_{i,\epsilon}\in (-\delta,\delta)$ such that
\[ [y_i, s_i'+{\bf \epsilon}\cdot \vec{t}]=[z_{i,\epsilon}, s'+\delta_{i,\epsilon}], i=1,2. \]
This is equivalent to say
\[ s_i'+{\bf \epsilon}\cdot \vec{t}-(s'+\delta_{i,\epsilon})\in\mathbb{Z}\ \ \text{and}\ \ T^{s_i'+{\bf \epsilon}\cdot \vec{t}-(s'+\delta_{i,\epsilon})}y_i=z_{i,\epsilon},i=1,2. \]

Since
\[ s_1'+{\bf \epsilon}\cdot \vec{t}-(s'+\delta_{1,\epsilon})-\big( s_2'+{\bf \epsilon}\cdot \vec{t}-(s'+\delta_{2,\epsilon})\big)\leq |s_1'-s_2'|+|\delta_{1,\epsilon}-\delta_{2,\epsilon}|<4\delta<1,\]
 we have $ s_1'+{\bf \epsilon}\cdot \vec{t}-(s'+\delta_{1,\epsilon})=s_2'+{\bf \epsilon}\cdot \vec{t}-(s'+\delta_{2,\epsilon})$ and we denote this common value by $k_{\epsilon}\in\mathbb{Z}$. Further, note that $ T^{k_{\epsilon}}y_i=z_{i,\epsilon}\in B(z,\delta)$, whence we have
 \begin{equation}\label{eq 9.1}
 \rho_{X}(T^{k_{\epsilon}}y_1, T^{k_{\epsilon}}y_2)<2\delta.
 \end{equation}

 {\bf Claim 2}. $\left|\epsilon\cdot [\vec{t}]-k_{\epsilon}\right|\leq d+3$,  where $[\vec{t}]=([t_1],\ldots,[t_d])\in\mathbb{Z}^{d}$.

This is followed by the following estimation.
\begin{align*}
|\epsilon\cdot [\vec{t}]-k_{\epsilon}|=&|\epsilon\cdot [\vec{t}]-s_1'-{\bf \epsilon}\cdot \vec{t}+s'+\delta_{1,\epsilon}|\\
\leq &  |\epsilon\cdot [\vec{t}]-{\bf \epsilon}\cdot \vec{t}|+|s_1'|+|s'|+|\delta_{1,\epsilon}|\\
\leq & d+1+1+1=d+3.
\end{align*}

 Now by the relation between $\delta$ and $\eta$,  (\ref{eq 9.1}) implies that
 \[\rho_{X}(T^{\epsilon\cdot [\vec{t}]}y_1, T^{\epsilon\cdot [\vec{t}]}y_2)<\eta, \forall {\bf \epsilon}\in \{0,1\}^{d}\setminus\{{\bf 0}\} .\]
 This shows that $(x_1,x_2)\in{\bf RP}^{[d]}(X,T)$.
 \end{proof}

  Now we are ready to show Theorem \ref{mod RP vs suspension}.
 \begin{proof}[Proof of Theorem \ref{mod RP vs suspension}]
Let $\pi:X\rightarrow X_d$ and $\tilde{\pi}: \tilde{X}\rightarrow \tilde{X}_{d}$ be the factor maps. Define $\xi: \widetilde{X_d}\rightarrow \tilde{X}_{d}$ by $[\pi(x), s]\mapsto \tilde{\pi}([x, s])$.

 (1) $\xi$ is well defined.

 Suppose that $[\pi(x), s]=[\pi_{d}(y),t]\in \widetilde{X_{d}}$. Then $s-t\in\mathbb{Z}$ and $\pi(y)=T^{s-t}\pi(x)=\pi(T^{s-t}x)$. So we have $(T^{s-t}x, y)\in {\bf RP}^{[d]}(X,T)$. By Lemma \ref{RP for suspension}, we have $([x, s], [y,t])\in {\bf RP}^{[d]}(\tilde{X}, \mathbb{R})$. Thus $\tilde{\pi}([x,s])=\tilde{\pi}([y,t])$. This shows that $\xi$ is well defined.

 (2) $\xi$ is bijective.

 Suppose that $\tilde{\pi}([x,s])=\tilde{\pi}([y,t])$. Then we have $([x, s], [y,t])\in {\bf RP}^{[d]}(\tilde{X}, \mathbb{R})$. By Lemma \ref{RP for suspension},  $s-t\in\mathbb{Z}$ and $(T^{s-t}x, y)\in {\bf RP}^{[d]}(X,T)$. Thus $\pi(y)=\pi(T^{s-t}x)=T^{s-t}\pi(x)$.  It follows that $[x,s]=[y,t]$. This shows that $\xi $ is injective. Finally, it is clear that $\xi$ is surjective by the definition.

The continuity of $\xi$ is followed from the property of quotient topology. Next we show that $\xi$  is a conjugacy between $(\widetilde{X_d},\mathbb{R})$ and $(\tilde{X}_d,\mathbb{R})$. For this, we need to show that
\[T^{t}\tilde{\pi}([x,s])=\xi(T^{t}[\pi(x),s]).\]
But this is followed since
\[T^{t} \tilde{\pi}([x,s])=\tilde{\pi}(T^{t}[x,s])=\tilde{\pi}([x,s+t])=\xi([\pi(x),s+t])=\xi(T^{t}[\pi(x),s]).\]

 \end{proof}

 \begin{prop}\label{pass from suspension}
Let $(X,T)$ be a  system and $(\tilde{X},\mathbb{R})$ be its suspension flow.  If there is $ A\subset \mathbb{R}$ and $[x,s]\in \tilde{X}$  such that $\{T^{t}[x,s]: t\in A\}$ is dense in $\tilde{X}$, then $\{T^{[t]}x: t\in A\}$ is dense in $X$.
 \end{prop}

\begin{proof}
Fix $y\in X$ and a neighborhood $U$ of $y$. Since $\{T^{t}[x,s]: t\in A\}$ is dense in $\tilde{X}$, there is some $t\in A$ such that
$$T^{t}[x,s]=[x, s+t]\in [U, s, 1/2].$$
Thus there is some $z\in U$ and $\delta\in (-1/2,1/2)$ such that $q((x, s+t))=q((z, s+\delta))$. This implies that
\[ s+t-(s+\delta)=t-\delta \in\mathbb{Z} \ \ \text{ and }\ \ T^{t-\delta}x=z.\]
Note that $t-\delta=[t]+\{t\} -\delta$ and $\{t\}-\delta \in (-1/2,3/2)$. But it follows from $t-\delta\in\mathbb{Z}$ that $\{t\}-\delta=0$ or $1$. Thus we have either $T^{[t]}x=z$ or $T^{[t]+1}x=z$. Since $y$ and $U$ are arbitrary, we conclude that
$\{T^{[t]}x: n\in A\}\cup \{T^{[t]+1}x: n\in A\}$ is dense in $X$.

\medskip
Next we show that $\{T^{[t]}x: n\in A\}$ is dense in $X$. To the contrary, there is some  $y\in X$ and a neighborhood $U$ of $y$ such that $T^{[t]}x\notin U$ for any $t\in A$. WLOG, we assume $s\in[0,1)$. By the previous arguments, we have
\[ A_{U}:=\{t\in A: T^{t}[x,s]=[x, s+t]=[T^{[t]}x, s+\{t\}]\in [U, s,1/2]\}\neq\emptyset.\]
Since  $T^{[t]}x\notin U, \forall t\in A$, we have that $\{t\}\in(1/2,1), \forall t\in A_{U}$. On the other hand,
\[ B_{U}:=\{t\in A: T^{t}[x,s]\in [U, s+1/4,1/4]\}\neq\emptyset.\]
For $t\in B_{U}$, there is some $i\in\mathbb{Z}$  such that $T^{t}[x,s]=[T^{[t]+i}x, s+\{t\}-i]$ and $T^{[t]+i}x\in U, |\{t\}-i-1/4|<1/4$.  Since $\{t\}\in[0,1)$, we must have $i=0$. But then $T^{[t]}x\in U$, contradicts our assumption. This shows that $\{T^{[t]}x: n\in A\}$ is dense in $X$.
\end{proof}

 \begin{cor}
Let $(X,T)$ be a minimal system and $p(t)$ be a non-constant real polynomial. Then there is some $x\in X$ and uncountably many $\beta$'s in $\mathbb{R}$ such that $\{ T^{[p(\beta n)]}x: n\in\mathbb{Z}\}$ is dense in $X$.
\end{cor}
\begin{proof}
According to Theorem \ref{polynomial orbit in R} and its proof, there are uncountably many $\beta$'s in $\mathbb{R}$ and $[x,s]\in\tilde{X}$ such that  $\{T^{p(\beta n)}[x,s]: n\in\mathbb{Z}\}$ is dense in $\tilde{X}$. Then by Proposition \ref{pass from suspension}, $\{ T^{[p(\beta n)]}x: n\in\mathbb{Z}\}$ is dense in $X$.
\end{proof}

\appendix{}

\section{Proof of Theorem \ref{nilfunction}}

In the appendix we give the proofs of Theorem \ref{nilfunction}. We restate the result as follows.

%Let $G$ be a topological group. Let $G_1=G$ and $G_{i+1}=[G, G_i]$ for $i\in \N$, where $[G,G_i]$ is the closed subgroup generated by $\{[g,h]=g^{-1}h^{-1}gh: g\in G, h\in G_i\}$. If there is some $s\in\N$ such that $G_{s+1}=\{e\}$, then we say $G$ is a {\it $s$-step nilpotent group}.

\begin{thm}%\label{nilfunction}
 Let  $\{T^{t}\}_{t\in\mathbb{R}}$ be an ergodic flow on $(X,\mathcal{X},\mu)$. Let $f\in L^{\infty}(\mu)$ and $\a_1,\ldots,\a_k$ be distinct nonzero real numbers. Then the function
 \[I_{f}(k,t):=\int_{X} f(x)\cdot f(T^{\a_1 t}x)\cdot\ldots\cdot f(T^{\a_k t}x)d\mu(x)\]
 is the sum of a function tending to zero in uniform density and $k$-step nilfunction.
\end{thm}

For the notions of Host-Kra factors $\mathcal{Z}_{k}$ and Gowers' norms $\vertiii{\cdot}_{k}$, one can refer to \cite{HK05} or \cite{BHK05}.

\subsection{Reduction to nilflows}
To show Theorem \ref{nilfunction}, it suffices to consider the nilflows. Potts shows in \cite[Proposition]{P11}   that  for any $f_1,f_2,\ldots,f_k\in L^{\infty}(\mu)$, one has
\begin{equation}\label{eqA1}
\limsup_{R\rightarrow \infty}\left\|\frac{1}{R}\int_{0}^{R}T^{\a_1 t}f_1 \cdot\ldots\cdot T^{\a_kt}f_kdt\right\|_{L^{2}(\mu)}\leq \min_{1\leq i\leq k}\vertiii{f_i}_{k}.
\end{equation}
This also implies that
\begin{equation}\label{eqA1}
\lim_{\rho\rightarrow \infty}\sup_{\rho\in\R}\left\|\frac{1}{R}\int_{\sigma}^{\sigma+\rho}T^{\a_1 t}f_1 \cdot\ldots\cdot T^{\a_kt}f_kdt\right\|_{L^{2}(\mu)}\leq \min_{1\leq i\leq k}\vertiii{f_i}_{k}.
\end{equation}

Then it follows that
\begin{lem} \label{transfer to factor}
Let $f$ be a bounded function on $X$ and $g=\mathbb{E}(f\mid \mathcal{Z}_{k})$. Then $I_{f}(k,t)-I_{g}(k,t)$ converges to zero in uniform density.
\end{lem}
\begin{proof}
Let $\a_0=0$.  Then
\begin{align*}
&\lim_{\rho\rightarrow \infty}\sup_{\sigma\in \mathbb{R}}\frac{1}{\rho}\int_{\sigma}^{\sigma+\rho} |I_{f}(k,t)-I_{g}(k,t)|dt\\
=&\lim_{\rho\rightarrow \infty}\sup_{\sigma\in \mathbb{R}}\frac{1}{\rho}\int_{\sigma}^{\sigma+\rho}\left|\int_{X}\left(\prod_{i=0}^{k}f(T^{\a_i t}x)-\prod_{i=0}^{k}g(T^{\a_i t}x)\right)d\mu(x)\right|dt\\
\leq& \lim_{\rho\rightarrow \infty}\sup_{\sigma\in \mathbb{R}}\frac{1}{\rho}\int_{\sigma}^{\sigma+\rho}\int_{X} \sum_{i=1}\left|g(T^{\a_0 t}x)\cdots g(T^{\a_{i-1}t}x)f(T^{\a_i t}x)-g(T^{\a_i t}x)f(T^{\a_{i+1}t}x)\cdots f(T^{\a_k t}x) \right|\\
\leq&  \|f-g\|_{L^{2}(\mu)}\vertiii{f}_{k}+k\|g\|_{L^{2}(\mu)}\min\{\vertiii{f}_{k},\vertiii{f-g}_{k}\}=0.
\end{align*}

\end{proof}

\subsection{Multiple averages on nilmanifolds}

Let $G$ be a connected simply connected $k$-steped nilpotent Lie group and $\Gamma$ be a cocompact lattice of $G$. Let $\mathfrak{g}$ denote the Lie algebra of $G$ and $\exp: \mathfrak{g}\rightarrow G$ denote the exponential map. Then $\exp$ is a diffeomorphism. Let $\log$  denote the inverse of $\exp$. Now for any $g\in G$ and $t\in\mathbb{R}$, we define
\[ g^{t}:=\exp( t\log g),\]
which is the unique homomorphism $\mathbb{R}\rightarrow G$ such that $1\mapsto g$. It follows from the uniqueness that
\begin{equation}\label{eqA2}
hg^{t}h^{-1}=(hgh^{-1})^{t},\forall g,h\in G \text{ and } t\in \R.
\end{equation}
 This property will be used frequently.

\subsubsection{Constructions of $\widetilde{G}$ and $\widetilde{G}^{*}$}
For a real number $\alpha$ and $n\in\mathbb{N}$, define ${\alpha\choose n}=\frac{1}{n!}\alpha(\alpha-1)\cdots(\alpha-n+1)$. Let $\a_1,\ldots,\a_k$ be distinct nonzero real numbers. Define the map $j: G\times G_1\times G_2\times\cdots\times G_{k}\rightarrow G^{k+1}$ by
\[j(g,g_1,g_2,\ldots,g_k)=\left(g, gg_{1}^{{\a_1\choose 1}}g_2^{{\a_1\choose 2}}\cdots g_{k}^{{\a_1\choose k}},  \ldots, gg_{1}^{{\a_k\choose 1}}g_2^{{\a_k\choose 2}}\cdots g_{k}^{{\a_k\choose k}}\right), \]
and let $\widetilde{G}$ denote the range of the map $j$:
\[ \widetilde{G}=j(G\times G_1\times G_2\times\cdots\times G_{k}).\]
Note that when $\a_1=1,\a_2=2,\ldots,\a_k=k$, this coincides with the $\tilde{G}$ defined in subsection \ref{Nd for commuting}.

Similarly, define a map $j^{*}: G_1\times G_2\times\cdots\times G_{k}\rightarrow G^{k}$ by
\[j^{*}(g_1,g_2,\ldots,g_k)=\left(g_{1}^{{\a_1\choose 1}}g_2^{{\a_1\choose 2}}\cdots g_{k}^{{\a_1\choose k}},  \ldots, g_{1}^{{\a_k\choose 1}}g_2^{{\a_k\choose 2}}\cdots g_{k}^{{\a_k\choose k}}\right), \]
and let $\widetilde{G}^{*}$ denote the range of the map $j$:
\[ \widetilde{G}^{*}=j^{*}(G_1\times G_2\times\cdots\times G_{k}).\]

The following results are shown in \cite{Leib10}.
\begin{thm}\label{G tilde}
\begin{enumerate}
\item[(1)] $\widetilde{G}$ is a rational subgroup of $G^{k+1}$, that is $\widetilde{G}$ is a closed subgroup of $G^{k+1}$ and $\widetilde{G}\cap \Gamma^{k+1}$ is cocompact in $\widetilde{G}$.
\item[(2)] The commutator group $(\widetilde{G})_{2}$ of $\widetilde{G}$ is
\[(\widetilde{G})_{2}=\widetilde{G}\cap G_2^{k+1}=j(G_2\times G_2\times G_2\times G_3\times\cdots\times G_k).\]
\item[(3)] $\widetilde{G}^{*}$ is a subgroup of $G^{k}$.
\item[(4)] For $g\in G$ and $(g_1,g_2,\ldots,g_k)\in \widetilde{G}^{*}$, one has
\[(gg_1g^{-1},gg_2g^{-1},\ldots,gg_{k}g^{-1})\in \widetilde{G}^{*}.\]
\end{enumerate}
\end{thm}
The item (4) above is implied by (\ref{eqA2}).

Using the representation of Malcev's basis, one can show that $j$ and $j^{*}$ are injective. Let $\mathfrak{g}_{i}$ be the Lie algebra of $G_{i}$, $i=1,2,\ldots,k$. Then
\[ \mathfrak{g}_{1}\rhd\mathfrak{g}_2\rhd\cdots\rhd\mathfrak{g}_{k}.\]
There is  a Malcev's basis $\{X_1,X_2,\ldots, X_{n}\}$ for $G/\Gamma$ such that
\begin{enumerate}
\item[(1)] There are $1= \ell_1<\ell_2<\cdots<\ell_{k}< n$ such that $\mathfrak{g}_{i}={\rm Span}\{X_{\ell_i},X_{\ell_i+1}\ldots, X_n\}$ for every $i\in\{1,2,\ldots,k\}$;
\item[(2)] Every element $g\in G_{i}$ is uniquely represented as $\exp(t_{\ell_i}X_{\ell_{i}}+\cdots t_{n}X_n)$ with $t_{\ell_i},\ldots,t_n\in\R$;
\item[(3)] The elements in $\Gamma$ is consisting of $\{\exp(m_1X_1+\cdots+m_nX_n): m_1,\ldots,m_n\in\Z\}$.
\end{enumerate}
For each $i\in\{1,2,\ldots,k\}$, let ${\bf X}_{i}=(X_{\ell_i},\ldots, X_{\ell_{i+1}-1})$ and for ${\bf t}=(t_{\ell_i},\ldots,t_{\ell_{i+1}-1})\in\mathbb{R}^{\ell_{i+1}-\ell_i}$, write
\[{\bf t}{\bf X}_{i}=t{\ell_i} X_{\ell_i}+\cdots+t_{\ell_{i+1}-1}X_{\ell_{i+1}-1}.\]
Then each $g\in G$ is represented as $\exp({\bf t}_1{\bf X}_1+{\bf t}_2{\bf X}_2+\cdots+{\bf t}_k{\bf X}_k)$, and we write
\[ V_i(g)={\bf t}_{i}, \forall i\in\{1,2,\ldots,k\}.\]

We only show that $j^{*}$ is injective and then $j$ is clearly injective. For this, let $(g_1,g_2,\ldots,g_k)$, $(h_1,h_2,\ldots,h_k)\in G_1\times G_2\times\cdots\times G_k$ such that
\[j^{*}(g_1,g_2,\ldots,g_k)=j^{*}(h_1,h_2,\ldots,h_k).\]
 Further, we write $g_i=\exp({\bf t}_{i,i}{\bf X}_{i}+\cdots+{\bf t}_{i,k}{\bf X}_{k})$ and $h_i=\exp({\bf s}_{i,i}{\bf X}_{i}+\cdots+{\bf s}_{i,k}{\bf X}_{k})$, for each $i\in\{1,\ldots,k\}$. Then we have
 \begin{equation}
\prod_{i=1}^{k} \exp\left({\alpha_j\choose i}({\bf t}_{i,i}{\bf X}_{i}+\cdots+{\bf t}_{i,k}{\bf X}_{k})\right)=\prod_{i=1}^{k} \exp\left({\alpha_j\choose i}({\bf s}_{i,i}{\bf X}_{i}+\cdots+{\bf s}_{i,k}{\bf X}_{k})\right)
 \end{equation}
for any $j\in\{1,2,\ldots,k\}$. For each $j\in\{1,2,\ldots,k\}$, let
\[ a_{j}=g_{1}^{{\a_j\choose 1}}g_2^{{\a_j\choose 2}}\cdots g_{k}^{{\a_j\choose k}} \ \ \text{ and }\ \ \ b_{j}=h_{1}^{{\a_j\choose 1}}h_2^{{\a_j\choose 2}}\cdots h_{k}^{{\a_j\choose k}}.\]
By the  Baker-Campbell-Hausdorff formula, one has
\[ V_{i}(a_j)=V_i(g_{1}^{{\a_j\choose 1}}g_2^{{\a_j\choose 2}}\cdots g_{i}^{{\a_j\choose i}} ) \text{ and }
V_{i}(b_j)=V_i(h_{1}^{{\a_j\choose 1}}h_2^{{\a_j\choose 2}}\cdots h_{i}^{{\a_j\choose i}} ),\]
for each $i\in\{1,2,\ldots,k\}$.

We are going to show that ${\bf t}_{i,d}={\bf s}_{i,d}$ for any $i,d\in\{1,2,\ldots,k\}$ with $i\leq d$.

First, it is clear that $V_{1}(a_j)={\a_j\choose 1} {\bf t}_{1,1}$ and $V_{1}(b_j)={\a_j\choose 1} {\bf s}_{1,1}$. Thus ${\bf t}_{1,1}={\bf s}_{1,1}$. Suppose that we have shown that $V_{i}(a_j)=V_{i}(b_j)$, that is ${\bf t}_{c,d}={\bf s}_{c,d}$ for any $c\in\{1,2,\ldots, i\}$ and $d\in\{c,c+1,\ldots, i\}$ for some $i\in\{1,2,\ldots,k-1\}$. By the  Baker-Campbell-Hausdorff formula, we have
\[V_{i+1}(a_j)=\sum_{r=1}^{i+1}{\a_j\choose i+1} {\bf t}_{r,i+1}+ \phi_{i+1}( {\bf t}_{1,1}, \ldots,{\bf t}_{1,i}, {\bf t}_{2,2},\ldots,{\bf t}_{2,i},\ldots, {\bf t}_{i,i})\]
and
\[V_{i+1}(b_j)=\sum_{r=1}^{i+1}{\a_j\choose i+1} {\bf s}_{r,i+1}+ \phi_{i+1}( {\bf s}_{1,1}, \ldots,{\bf s}_{1,i}, {\bf s}_{2,2},\ldots,{\bf s}_{2,i},\ldots, {\bf s}_{i,i}),\]
where $\phi_{i+1}( {\bf t}_{1,1}, \ldots,{\bf t}_{1,i}, {\bf t}_{2,2},\ldots,{\bf t}_{2,i},\ldots, {\bf t}_{i,i})$ is a polynomial function. By the inductive hypothesis, we have
\[\phi_{i+1}( {\bf t}_{1,1}, \ldots,{\bf t}_{1,i}, {\bf t}_{2,2},\ldots,{\bf t}_{2,i},\ldots, {\bf t}_{i,i})=\phi_{i+1}( {\bf s}_{1,1}, \ldots,{\bf s}_{1,i}, {\bf s}_{2,2},\ldots,{\bf s}_{2,i},\ldots, {\bf s}_{i,i}).\]
Since $V_{i+1}(a_j)=V_{i+1}(b_j)$ for each $j\in\{1,2,\ldots,k\}$, one has
\[\sum_{r=1}^{i+1}{\a_j\choose i+1} {\bf t}_{r,i+1}=\sum_{r=1}^{i+1}{\a_j\choose i+1} {\bf s}_{r,i+1}.\]
Hence
\[\sum_{r=1}^{i+1}{\a_j\choose i+1} ({\bf t}_{r,i+1}- {\bf s}_{r,i+1})=0.\]
This implies that
\[ p(\a)=\sum_{r=1}^{i+1}{\a\choose i+1} ({\bf t}_{r,i+1}- {\bf s}_{r,i+1})\]
has $k$ distinct roots $\a_1,\a_2,\ldots,\a_k$. Since the degree of $p(\a)$ is no greater than $k$, one has $p(\a)\equiv 0$ and hence ${\bf t}_{r,i+1}= {\bf s}_{r,i+1}$ for any $r=1,2,\ldots,i+1$.

Therefore, it follows from the induction that $j^{*}$ is injective.

Next we show that $\widetilde{\Gamma}^{*}=\widetilde{G}^{*}\cap \Gamma^{k}$ is cocompact in $\widetilde{G}^{*}$. Then the cocompactness of other discrete subgroup such as $\widetilde{\Gamma}=\widetilde{G}\cap \Gamma^{k+1}$ in $\widetilde{G}$ are followed similarly. It suffices to show that there is a compact set $\Omega$ in $G_1\times G_2\times\cdots \times G_k$ such that for each ${\bf h}=(h_1,h_2,\ldots,h_k)\in \widetilde{G}^{*}$,
\[ {\bf h}j^{*}(\Omega) \cap \Gamma^{k}\neq \emptyset.\]

We will inductively define
\begin{align*}
&E_1:=[0, M_{1,1}]^{\ell_2-\ell_1}{\bf X}_1+\cdots +[0, M_{1, k}]^{n-\ell_{k}}{\bf X}_{k}\subset \mathfrak{g}_{1},\\
 &E_{2}:=[0, M_{2,2}]^{\ell_3-\ell_2}{\bf X}_2+\cdots +[0, M_{2, k}]^{n-\ell_{k}}{\bf X}_{k}\subset \mathfrak{g}_{2},\\ &\ldots,\\
 &E_{k}:=[0, M_{k,k}]^{n-\ell_k}{\bf X}_{k}\subset \mathfrak{g}_{k}
\end{align*}
  such that $\Omega=(\exp(E_1),\exp(E_2),\ldots,\exp(E_k))$ satisfies our requirements. For this, we may assume that ${\a_j\choose i}\neq 0$ for each $i,j\in\{1,2,\ldots,k\}$.

Firstly, define $M_{1,1}=1$. Then for any ${\bf h}\in \widetilde{G}^{*}$, $(g_2,\ldots,g_k)\in G_2\times\cdots\times G_k$  there is some ${\bf t}_{1}\in[0,M_{1,1}]^{\ell_2-\ell_1}$ such that for any ${\bf t_2}\in \R^{\ell_3-\ell_2},\ldots,{\bf t}_{k}\in\R^{n-\ell_k}$,
\[V_1\left({\bf h}j^{*}(\exp(\sum_{i=1}^{k}{\bf t}_i{\bf X}_{i}), g_2,\ldots,g_k)\right)\in \Z^{\ell_2-\ell_1}.\]

Suppose that we have defined $M_{1,1},\ldots,M_{1,i}, M_{2,2},\ldots,M_{2,i},\ldots,M_{i,i}$ such that for any ${\bf h}\in \widetilde{G}^{*}$, there is some $({\bf t}_{r,1}, {\bf t}_{r, 2},\ldots,{\bf t}_{r, i})\in [0, M_{r,1}]^{\ell_2-\ell_1}\times[0, M_{r,2}]^{\ell_3-\ell_2}\times\cdots\times[0, M_{r,i}]^{\ell_{i+1}-\ell_{i}} $, for each $r=1,2,\ldots,i$,  such that for any $(g_{i+1},\ldots,g_k)\in G_{i+1}\times\cdots\times G_k$  and any ${\bf t}_{r,i+1}\in \R^{\ell_{r,i+2}-\ell_{i+1}}$, $\ldots$, ${\bf t}_{r,k}\in\R^{n-\ell_k}$, one has
\[V_j\left({\bf h}j^{*}(\exp(\sum_{m=1}^{k}{\bf t}_{1,m}{\bf X}_{i}), \ldots,\exp(\sum_{m=1}^{k}{\bf t}_{i,m}{\bf X}_{i}),g_{i+1},\ldots,g_k)\right)\in \Z^{\ell_{j+1}-\ell_{j}},\]
for each $j=1,2,\ldots,i$.

Now let $M_{1,i+1}=\cdots=M_{i+1,i+1}=\max_{1\leq j\leq k}\frac{1}{{\a_j\choose i+1}}$. Then we claim that for any ${\bf h}\in \widetilde{G}^{*}$, there is some $({\bf t}_{r,1}, {\bf t}_{r, 2},\ldots,{\bf t}_{r, i})\in [0, M_{r,1}]^{\ell_2-\ell_1}\times[0, M_{r,2}]^{\ell_3-\ell_2}\times\cdots\times[0, M_{r,i+1}]^{\ell_{i+2}-\ell_{i+1}} $, for each $r=1,2,\ldots,i+1$,  such that for any $(g_{i+2},\ldots,g_k)\in G_{i+2}\times\cdots\times G_k$  and any ${\bf t}_{r,i+2}\in \R^{\ell_{r,i+3}-\ell_{i+1}},\ldots,{\bf t}_{r,k}\in\R^{n-\ell_k}$, one has
\[V_j\left({\bf h}j^{*}(\exp(\sum_{m=1}^{k}{\bf t}_{1,m}{\bf X}_{i}), \ldots,\exp(\sum_{m=1}^{k}{\bf t}_{i+1,m}{\bf X}_{i+1}),g_{2+1},\ldots,g_k)\right)\in \Z^{\ell_{j+1}-\ell_{j}},\]
for each $j=1,2,\ldots,i+1$. To this end, we fix ${\bf h}\in \widetilde{G}^{*}$. Then it follows from the inductive hypothesis that there is some $({\bf t}_{r,1}, {\bf t}_{r, 2},\ldots,{\bf t}_{r, i})\in [0, M_{r,1}]^{\ell_2-\ell_1}\times[0, M_{r,2}]^{\ell_3-\ell_2}\times[0, M_{r,i}]^{\ell_{i+1}-\ell_{i}} $, for each $r=1,2,\ldots,i$,  such that for any $(g_{i+1},\ldots,g_k)\in G_{i+1}\times\cdots\times G_k$  and any ${\bf t}_{r,i+1}\in \R^{\ell_{r,i+2}-\ell_{i+1}},\ldots,{\bf t}_{r,k}\in\R^{n-\ell_k}$, one has
\[V_j\left({\bf h}j^{*}(\exp(\sum_{m=1}^{k}{\bf t}_{1,m}{\bf X}_{i}), \ldots,\exp(\sum_{m=1}^{k}{\bf t}_{i,m}{\bf X}_{i}),g_{i+1},\ldots,g_k)\right)\in \Z^{\ell_{j+1}-\ell_{j}},\]
for each $j=1,2,\ldots,i$. Note that for each ${\bf t}_{r,i+1}\in [0, M_{r,i+1}]^{\ell_{i+2}-\ell_{i+1}}$, one has
\begin{align}
&V_{i+1}\left({\bf h}j^{*}(\exp(\sum_{m=1}^{k}{\bf t}_{1,m}{\bf X}_{i}), \ldots,\exp(\sum_{m=1}^{k}{\bf t}_{i,m}{\bf X}_{i}),g_{i+1},\ldots,g_k)\right)\label{eqAcoc} \\
=&\sum_{r=1}^{i+1}{\a_{j}\choose i+1}{\bf t}_{r,i+1}+\phi_{i+1}( {\bf t}_{1,1}, \ldots,{\bf t}_{1,i}, {\bf t}_{2,2},\ldots,{\bf t}_{2,i},\ldots, {\bf t}_{i,i}).\nonumber
\end{align}
Thus when ${\bf t}_{r,i+1}$ runs through $[0, M_{r,i+1}]^{\ell_{i+2}-\ell_{i+1}}$,  (\ref{eqAcoc}) must take some integral values since ${\a_j\choose i+1}[0, M_{r,i+1}]^{\ell_{i+2}-\ell_{i+1}}$ contains a cube with length greater than $1$.  In particular, such $$({\bf t}_{r,1}, {\bf t}_{r, 2},\ldots,{\bf t}_{r, i})\in [0, M_{r,1}]^{\ell_2-\ell_1}\times[0, M_{r,2}]^{\ell_3-\ell_2}\times[0, M_{r,i}]^{\ell_{i+1}-\ell_{i}} $$ and ${\bf t}_{r,i+1}\in [0, M_{r,i+1}]^{\ell_{i+2}-\ell_{i+1}}$ satisfy our requirements.

Then it follows from the induction that $\widetilde{G}^{*}\cap \Gamma^{k}$ is cocompact in $\widetilde{G}^{*}$.

\subsubsection{Nilmanifolds $\widetilde{X}$ and $\widetilde{X}_{x}$}

Let $\{T^{t}\}_{t\in\mathbb{R}}$ be a minimal  $\mathbb{R}$-flow on $X=G/\Gamma$ defined by
\[ T^{t}(g\Gamma)=b^{t}g\Gamma.\]

Now define two elements
\[ \tilde{b}=(1, b^{\a_1}, b^{\a_2},\ldots, b^{\a_{k}}) \text{ and } b^{\Delta}=(b,b,\ldots,b)\]
of $\widetilde{G}$ and the elmment
\[ \tilde{b}^{*}=(b^{\a_1}, b^{\a_2},\ldots, b^{\a_{k}})\in \widetilde{G}^{*}.\]
%Further, write $\tilde{T}^{t}$ and $T^{\Delta t}$ for the translations on $X^{k+1}$ by $(\tilde{b})^{t}$ and $(b^{\Delta})^{t}$, respectively, and $\tilde{T}^{*t}$ for the translation on $X^{k}$ by $(\tilde{b}^{*})^{t}$.

Define $\widetilde{\Gamma}=\widetilde{G}\cap \Gamma^{k+1}$. Then $\widetilde{\Gamma}$ is a discrete cocompact subgroup of $\widetilde{G}$.  Set
\[ \widetilde{X}=\widetilde{G}/\widetilde{\Gamma}\]
and let $\widetilde{\mu}$ be the Haar measure on this nilmanifold. Note that $\widetilde{X}$ is imbedded in $X^{k+1}$ in a natural way. Since $(\tilde{b})^{t}$ and $(b^{\Delta})^{t}$ are in $\widetilde{G}$, we can define translations  $\widetilde{T}^{t}$ and $T^{\Delta t}$  on the nilmanifold $\widetilde{X}$ by
\[ \widetilde{T}^{t}((h, h_1,\ldots,h_k)\widetilde{\Gamma})=(h, b^{\a_1 t}h_1,\ldots, b^{\a_k t}h_k)\widetilde{\Gamma}\]
and
\[ T^{\Delta t}((h, h_1,\ldots,h_k)\widetilde{\Gamma})=(b^{t}h, b^{ t}h_1,\ldots, b^{ t}h_k)\widetilde{\Gamma}.\]
\begin{lem}
The nilmanifold $\widetilde{X}$ is ergodic for the action spanned by  $\widetilde{T}^{s}$ and $T^{\Delta t}$.
\end{lem}

\begin{proof}
Since $G$ is connected, it suffices to show that the induced action on
 $\widetilde{G}/(\widetilde{G})_{2}\widetilde{\Gamma} $ by  $\widetilde{T}^{s}$ and $T^{\Delta t}$ is ergodic.

 By the part (2) in Theorem \ref{G tilde}, the map
 \[ (g_0,g_1,\ldots,g_k)\mapsto (g_0 \ {\rm mod}\ G_2, \ \ \ g_1 \ {\rm mod}\ G_2)\]
 induces an isomorphism form $\widetilde{G}/(\widetilde{G})_{2}$ onto  $G/G_2\times G/G_2$. Thus the compact abelian group $\widetilde{G}/(\widetilde{G})_{2}\widetilde{\Gamma}$ can be identified with $G/G_2\Gamma\times G/G_2\Gamma_1$, and the transformations induced by $\widetilde{T}^{s}$ and $T^{\Delta t}$ are $id\times T^{\a_1 s}$ and $T^{t}\times T^{t}$, respectively. Clearly, the action spanned by these transformations is ergodic.
\end{proof}

 For $x\in X$, define
\[ \widetilde{X}_{x}=\{(x_1,x_2,\ldots,x_k)\in X^{k}: (x,x_1,x_2,\ldots,x_k)\in\widetilde{X}\}.\]

Clearly, $\widetilde{X}_{x}$ is invariant under translations by elements of $\widetilde{G}^{*}$. We will give $\widetilde{X}_{x}$ the structure of a nilmanifold as a quotient of $\widetilde{G}^{*}$.

Fix $x\in X$ and let $a$ be a lift of $x$ in $G$. By the remark above,  we treat  $(x,x,\ldots,x)$ as an element in $\widetilde{X}_{x}$.

Now let $(x_1,x_2,\ldots,x_k)\in \widetilde{X}_{x}$. Since $(x,x_1,x_2,\ldots,x_k)\in \widetilde{X}$, we can lift it to an element of $\widetilde{G}$, saying $(g,gg_1,gg_2,\ldots,gg_{k})$ with $g\in G$ and $(g_1,g_2,\ldots,g_k)\in\widetilde{G}^{*}$. Set $\gamma=a^{-1}g$ and $h_i=a\gamma g_i\gamma^{-1}a^{-1}$ for $i=1,\ldots,k$. Then we have $\gamma \in \Gamma$ and $(1,h_1,h_2,\ldots,h_k)\in \widetilde{G}$ by Theorem \ref{G tilde} (4). Further we have
\[(g,gg_1,gg_2,\ldots,gg_{k})=(1,h_1,h_2,\ldots,h_k)\cdot (a,a,a,\ldots,a)\cdot(\gamma,\gamma,\gamma,\ldots,\gamma).\]
This shows that $(x_1,x_2,\ldots,x_k)$ is the image of $(x,x,\ldots,x)$ under the translation by $(h_1,h_2,\ldots,h_k)$, which belongs to $\widetilde{G}^{*}$. Therefore, the action of $\widetilde{G}^{*}$ on $\widetilde{X}_{x}$ is transitive. The stabilizer of $(x,x,\ldots,x)$ for this action is
\[\widetilde{\Gamma}_{x}:=\{ (a\gamma_1 a^{-1}, a\gamma_2 a^{-1},\ldots, a\gamma_k a^{-1}): (\gamma_1,\gamma_2,\ldots,\gamma_k)\in \widetilde{G}^{*}\cap \Gamma^{k}\},\]
Note that $\widetilde{G}^{*}\cap \Gamma^{k}$ is a cocompact lattice in $\widetilde{G}^{*}$.

Therefore, we can identify $\widetilde{X}_{x}$ with the nilmanifold $\widetilde{G}^{*}/\widetilde{\Gamma}_{x}$ via
\[(x_1,x_2,\ldots,x_k)\mapsto (h_1,h_2,\ldots,h_k)\widetilde{\Gamma}_{x}.\]

 Let $\widetilde{\mu}_{x}$ be the Haar measure on $\widetilde{X}_{x}$.

\begin{lem}
$\widetilde{\mu}=\int_{X}\delta_{x}\otimes \widetilde{\mu}_{x}d\mu(x)$.
\end{lem}
\begin{proof}
Let $\nu=\int_{X}\delta_{x}\otimes \widetilde{\mu}_{x}d\mu(x)$, which is concentrated on $\widetilde{X}$. Since $\widetilde{\mu}$ is the unique invariant measure on $\widetilde{X}$ under the action spanned by  $\tilde{T}^{t}$ and $T^{\Delta t}$, it suffices to show that $\nu$ is invariant under $\tilde{T}^{t}$ and $T^{\Delta t}$.

For any $(x_1,x_2,\ldots,x_k)\in\widetilde{X}_{x}$ and $t\in \mathbb{R}$, it is clear that $(x, b^{\a_1 t}x_1,b^{\a_2 t}x_2,\ldots, b^{\a_k t}x_k)\in \widetilde{X}$ since $(1,b^{\a_1 t},b^{\a_2 t},\ldots,b^{\a_k t})\in \widetilde{G}$. Thus one can define the translation $\widetilde{T}^{*t}$ on $\widetilde{X}_{x}$  by $(\tilde{b}^{*})^{t}=(b^{\a_1 t},b^{\a_2 t},\ldots,b^{\a_k t})$. As $ (\tilde{b}^{*})^{t}$ belongs to $\widetilde{G}^{*}$, $\widetilde{\mu}_{x}$ is invariant under $\widetilde{T}^{*t}$. Therefore, $\delta_{x}\otimes \widetilde{\mu}_{x}$ is invariant under $\widetilde{T}^{t}=id \times \widetilde{T}^{*t}$.

Let $x\in X$. Consider the image of $\widetilde{\mu}_{x}$ under $T^{t}\times T^{t}\times\cdots\times T^{t}$ ($k$ times). This measure is concentrated on $X_{T^{t}x}$ and is invariant under the action of $\widetilde{G}^{*}$ by noting that  if $(x_1,x_2,\ldots,x_k)\in \widetilde{X}_{x}$ is identified with $(h_1,h_2,\ldots,h_k)\widetilde{\Gamma}_{x}\in \widetilde{G}^{*}/\widetilde{\Gamma}_{x}$ then $(T^{t}x_1,T^{t}x_2,\ldots, T^{t}x_k)$ belongs to $\widetilde{X}_{T^{t}x}$ and is identified with $(b^{t}h_1b^{-t}, b^{t}h_2b^{-t},\ldots, b^{t}h_kb^{-t})\widetilde{\Gamma}_{x}\in \widetilde{G}^{*}/\widetilde{\Gamma}_{x}$. Thus it is equal to the Haar measure $\widetilde{\mu}_{T^{t}x}$. Therefore, the image of $\delta_{x}\otimes \widetilde{\mu}_{x}$ under $T^{\Delta t}$ is $\delta_{T^{t}x}\otimes \widetilde{\mu}_{T^{t}x}$. This shows that $\nu$ is invariant under $T^{\Delta t}$.
\end{proof}

The following result is shown in \cite{Leib10} and also can be proved in the way as in \cite[Section 5.4]{BHK05}.
\begin{thm}\cite[Theorem 5.4]{BHK05}
Let $f_1,f_2,\ldots, f_{k}$ be continuous functions on $X$ and let $\{\sigma_i\}$ and $\{\rho_i\}$ be two sequences of reals such that $\rho_i\rightarrow +\infty$. For $\mu$-almost every $x\in X$,
\[\frac{1}{\rho_i}\int_{\sigma_i}^{\sigma_i+\rho_i}f_1(T^{\a_1 t}x)f_2(T^{\a_2 t}x)\cdots f_{k}(T^{\a_k t}x)\rightarrow \int_{\widetilde{X}_{x}} f_1(x_1)f_2(x_2)\cdots f_{k}(x_k)d\widetilde{\mu}_{x}(x_1,x_2,\ldots, x_k)\]
as $i\rightarrow \infty$.
\end{thm}

\begin{cor}\cite[Corollary 5.5]{BHK05}
For $\mu$-almost every $x\in X$, the nilflow $(\widetilde{X}_{x}, \widetilde{\mu}_{x}, \{\widetilde{T}^{*t}\}_{t\in \R})$ is ergodic.
\end{cor}

\begin{cor}\cite[Corollary 5.6]{BHK05}\label{multi convergence}
Let $f_1,f_2,\ldots, f_{k}$ be continuous functions on $X$ and let $\{\sigma_i\}$ and $\{\rho_i\}$ be two sequences of reals such that $\rho_i\rightarrow +\infty$. For $\mu$-almost every $x\in X$ and for every $(g_1,g_2,\ldots,g_k)\in\widetilde{G}^{*}$,
\begin{align*}
\frac{1}{\rho_i}\int_{\sigma_i}^{\sigma_i+\rho_i} &f_1(T^{\a_1 t}g_1x)f_2(T^{\a_2 t}g_2x)\cdots f_{k}(T^{\a_k t}g_kx)\\
&\rightarrow \int_{\widetilde{X}_{x}} f_1(x_1)f_2(x_2)\cdots f_{k}(x_k)d\widetilde{\mu}_{x}(x_1,x_2,\ldots, x_k)
\end{align*}
as $i\rightarrow \infty$.
\end{cor}

\subsubsection{Construction of nilfunction $J_{f}(k,t)$}
Define
\[ H=\{(g,h)\in G\times G: hg^{-1}\in G_{2}\}.\]
Then $H$ is a connected closed subgroup of $G\times G$ and is a connected $k$-step nilpotent Lie group. The commutator subgroups $H_{j}$ of $H$ are given by
\[ H_{j}=\{(g,h)\in G_{j}\times G_{j}: hg^{-1}\in G_{j+1}\}.\]
We build the groups $\widetilde{H}$ and $\widetilde{H}^{*}$ in the same way that the groups  $\widetilde{G}$ and $\widetilde{G}^{*}$ are built from $G$.

The ergodic decomposition
\[ \mu\times \mu=\int_{Z} \mu_{s} d m(s)\]
of $\mu\times \mu$ under $T^{t}\times T^{t}$, where $m$ is the Haar measure on the Kronecker factor $Z$ of $X$. Actually, $Z=G/G_2\Gamma$ and the factor map $\pi: X\rightarrow Z$ is the natural projection $G/\Gamma\rightarrow G/G_2\Gamma$.

For each $s\in Z$, the measure $\mu_{s}$ is concentrated on the closed subset
\[ X_{s}=\{(x,y)\in X\times X: \pi(y)\pi(x)^{-1}=s\}\]
of $X\times X$.

$X_{s}$ can be given the structure of a nilmanifold, as the quotient of $H$. Then we can built the nilmanifold $\widetilde{X_{s}}$ in the way we define $\widetilde{X}$ and its Haar measure is denoted by $\widetilde{\mu_{s}}$. Further, for $(x,y)\in X_{s}$, the nilmanifold $\widetilde{X_{s}}_{(x,y)}$ is defined in the way we defined $\widetilde{X}_{x}$ and its Haar measure is denoted by $\widetilde{\mu_{s}}_{(x,y)}$. Then the following is an immediate corollary of Corollary \ref{multi convergence}.

\begin{prop}
Let $f$ be a continuous functions on $X$ and let $\{\sigma_i\}$ and $\{\rho_i\}$ be two sequences of reals such that $\rho_i\rightarrow +\infty$. For $m$-almost every $s\in Z$ and  for $\mu_{s}$-almost every $(x,y)\in X_{s}$ and for every $((g_1,g_2,\ldots,g_k),(h_1,h_2,\ldots,h_k))\in\widetilde{H}^{*}$,
\begin{align*}
\frac{1}{\rho_i}\int_{\sigma_i}^{\sigma_i+\rho_i} &\prod_{j=1}^{k} f(T^{\a_j t}g_j \cdot x)f(T^{\a_j t} h_j \cdot y)dt\\
&\rightarrow \int_{\widetilde{X_{s}}_{(x,y)}} \prod_{j=1}^{k} f(x_j)f(y_j)d\widetilde{\mu_{s}}_{(x,y)}\Big((x_1,x_2,\ldots, x_k),(y_1,y_2,\ldots, y_k)\Big)
\end{align*}
as $i\rightarrow \infty$.
\end{prop}

Define
\[ \vec{G}^{*}=\left\{{\bf g}=(g_1,g_2,\ldots,g_{k})\in\widetilde{G}^{*}: \big( (1,1,\ldots,1), (g_1,g_2,\ldots,g_{k})\big)\in\widetilde{H}^{*}\right\}. \]
Then $\vec{G}^{*}$ is a connected closed subgroup of $G^{k}$ and thus a connected Lie group. Since $j^{*}$ is injective and $H_{j}=\{(g,h)\in G_{j}\times G_{j}: hg^{-1}\in G_{j+1}\}$, one has
\[ \vec{G}^{*}=j^{*}(G_2\times G_3\times\cdots\times G_{k+1}).\]
Moreover, when ${\bf g}=(g_1,g_2,\ldots,g_{k})\in\widetilde{G}^{*}$ one has $({\bf g}, {\bf g})\in\widetilde{H}^{*}$. Thus
\[ \widetilde{H}^{*}=\left\{({\bf g}, {\bf h})\in\widetilde{G}^{*}\times\widetilde{G}^{*}: {\bf h}{\bf g}^{-1}\in\vec{G}^{*}\right\}\]
and $\vec{G}^{*}$ is a normal subgroup of $\widetilde{G}^{*}$.

Define
\begin{align*}
\vec{G}&=\left\{ (g,gh_1,gh_2,\ldots,gh_k): g\in G, (h_1,h_2,\ldots,h_k)\in\vec{G}^{*}\right\}\\
&=j(G\times G_2\times G_3\times\cdots\times G_{k+1}).
\end{align*}
Moreover, $\vec{G}$ is a normal subgroup of $\widetilde{G}$ and is a connected nilpotent Lie group.

Let $\vec{\Gamma}=\vec{G}\cap \Gamma^{k}$. Then $\vec{\Gamma}$ is a cocompact discrete subgroup of $\vec{G}$. Define
\[ \vec{X}=\vec{G}/\vec{\Gamma}\]
 and let $\vec{\mu}$ denote the Haar measure of this nilmanifold.  Clearly, $\vec{X}\subset \widetilde{X}$ as $\vec{G}\subset \widetilde{G}$.

 The following result is followed by the normality of $\vec{G}$ in $\widetilde{G}$. For detailed proof see \cite[Lemma 6.3]{BHK05}.
 \begin{lem}\label{inv of Gamma}
 This nilmanifold $\vec{X}$ and its Haar measure $\vec{\mu}$ are invariant under translation by any element of $\widetilde{\Gamma}$.
 \end{lem}

For $x\in X$, let
\[ \vec{X}_{x}=\{(x_1,x_2,\ldots,x_{k})\in X^{k}: (x,x_1,x_2,\ldots, x_{k})\in\vec{X}\}.\]
Then $\vec{X}_{x}$ can be identified with a nilmanifold and we let $\vec{\mu}_{x}$ denote its Haar measure.

\begin{lem}\cite[Lemma 6.4]{BHK05}
$\vec{\mu}=\int_{X} \delta_{x}\otimes \vec{\mu}_{x}d\mu(x)$.
\end{lem}

For a bounded function $f$ on $X$ and $t\in \R$, define

\[J_{f}(k, t)=\int_{\vec{X}} f(x_0)f(T^{\a_1 t}x_1)\cdots f(T^{\a_k t}x_k)d\vec{\mu}(x_0,x_1,\ldots, x_k).\]
 By Lemma , one has
\[J_{f}(k, t)=\int_{X}f(x)\int_{\vec{X}_{x}}f(T^{\a_1 t}x_1)\cdots f(T^{\a_k t}x_k)d\vec{\mu}_{x}(x_1,\ldots, x_k)d\mu(x).\]

 The following result is obtained in the same way as in \cite[Proposition 6.5]{BHK05} by replacing the sums with integrals.
\begin{prop}\label{I-J}
Let $f$ be a bounded function on $X$. Then the function $\{I_{f}(k, t)-J_{f}(k,t)\}$ converges to zero in uniform density.
\end{prop}

Next we are going to describe the nilflow $(Y, \nu, \{S^{t}\}_{t\in \R})$.

 Let $K$ denote the group $\widetilde{G}/\vec{G}$, let $p: \widetilde{G}\rightarrow K$ be the natural projection and let $\Lambda=p(\widetilde{\Gamma})$.

 Let $Y$ denote the nilmanifold $K/\Gamma$, $\nu$ be its Haar measure, $c=p(\tilde{b})\in K$ and $S^{t}$ be the translation by $c^{t}$ on $Y$.

\begin{lem}\cite[Lemma 7.1]{BHK05}
This nilflow $(Y, \nu, \{S^{t}\}_{t\in \R})$ is ergodic (and thus is uniquely ergodic and minimal).
\end{lem}

\begin{prop}\cite[Proposition 7.2]{BHK05}\label{nilfunction for f}
Let $f$ be a bounded function on $X$, then $J_{f}(k,t)$ is a nilfunction.
\end{prop}

\begin{proof}
Define the function $\psi$ on $\widetilde{G}$ by
\begin{equation}\label{eqAnilfunc}
\psi(g_0,g_1,\ldots, g_{k})=\int_{\vec{X}} \prod_{j=0}^{k}f(g_{j}\cdot x_{j}) d\vec{\mu}(x_0,x_1,\ldots,x_k).
\end{equation}
The function $\psi$ is continuous and satisfies
\[\psi(1, b^{\a_1 t},\ldots, b^{\a_k t})=J_{f}(k,t)\]
for every $t\in\mathbb{R}$.

The measure $\vec{\mu}$ is invariant under left translation by elements of $\vec{G}$ by definition, and by left translation by elements of $\widetilde{\Gamma}$ by Lemma \ref{inv of Gamma}. Thus the function $\psi$ is invariant under right translation by elements of $\vec{G}\widetilde{\Gamma}$.

Writing $q$ for the natural projection $\widetilde{G}\rightarrow Y=\widetilde{G}/\vec{G}\widetilde{\Gamma}$, we get that there is a continuous function $\phi$ on $Y$ satisfying $\psi=\phi\circ q$. Note that  for every $t\in \R$, $q(1, b^{\a_1 t},\ldots, b^{\a_k t})=S^{t}e_{Y}$ and hence
\[ J_{f}(k,t)=\phi(S^{t}e_{Y}).\]
This shows that $J_{f}(k,t)$ is a basic nilfunction.
\end{proof}

 \subsubsection{Proof of Theorem \ref{nilfunction}}

 We say a subset $A\subset \R$ is {\it discrete syndetic} if there are finitely many $t_1,t_2,\ldots,t_r\in \R$ such that
 \[ \R=(A-t_1)\cup(A-t_2)\cup\cdots\cup (A-t_r).\]
 The following result is clear from the definition.
 \begin{lem}\label{synd-->density}
 Let $A$ be a measurable subset of $\R$. If $A$ is discrete syndetic, then $A$ has positive lower Banach density, i.e.,
 \[ \lim_{\rho\rightarrow \infty}\inf_{\sigma\in\R} \frac{1}{\rho} \int_{\sigma}^{\sigma+\rho}{\bf 1}_{A}dt>0.\]
 \end{lem}

 \begin{defn}
 Let $\phi(t)$ be a real-valued function in $L^{\infty}(\R)$. The discrete syndetic supremum of $\phi$ is
 \[ {\rm Dsynd\textendash sup} \phi(t) :=\sup\{ c\in\R:\{ t\in \R: \phi(t)>c\} \text{ is discrete syndetic}\}.\]
 \end{defn}

 \begin{lem}\label{UD-lim vs Dsynd}
 Let $\phi(t)$ and $\psi(t)$ be real-valued functions on $L^{\infty}(\R)$. If ${\rm UD\textendash Lim}(\phi(t)-\psi(t))=0$, then ${\rm Dsynd\textendash sup} \phi(t)={\rm Dsynd \textendash sup} \psi(t)$.
 \end{lem}
\begin{proof}
 Suppose that ${\rm Dsynd\textendash sup} \phi(t)\neq {\rm Dsynd\textendash sup} \psi(t)$. We may assume that there are $c,d\in\R$ such that
 \[{\rm Dsynd\textendash sup} \phi(t)>c>d> {\rm Dsynd \textendash sup} \psi(t).\]
 Then $A:=\{ t\in \R: \phi(t)>c\}$ is measurable and discrete syndetic, and $B: =\{ t\in R: \psi(t)\leq d\}$ is thick, that is there is a sequences $(\sigma_i)$ and $(\rho_i)$ in $\R$ with $\rho_i\rightarrow \infty$ such that
 \[ B\supset \bigcup_{i=1}^{\infty}[\sigma_i, \sigma_i+\rho_i].\]
 Thus, by Lemma \ref{synd-->density}, we have
 \[ \limsup_{i\rightarrow\infty}\frac{1}{\rho_i}\int_{\sigma_i}^{\sigma_i+\rho_i} |\phi(t)-\psi(t)|dt\geq (c-d)\limsup_{i\rightarrow\infty}\frac{1}{\rho_i}\int_{\sigma_i}^{\sigma_i+\rho_i} {\bf 1}_{A}dt>0.\]
 This contradicts the condition.
 \end{proof}

 \begin{lem}\label{synd sup for nil}
 Let $(X, \{T^{t}\}_{t\in\R})$ be a minimal flow.  Let $f$ be a continuous function on $X$ and $x_0\in X$. Then
 \[{\rm Dsynd\textendash sup}f(T^tx_0) =\sup_{t\in \R} f(T^{t}x_0).\]
 \end{lem}
 \begin{proof}
Let $M=\sup_{t\in \R} f(T^{t}x_0)$ and $\epsilon>0$. Then
 \[ U=\{x\in X: f(x)>M-\epsilon\}\]
 is a nonempty open set in $X$.  Since $(X,\mu, \{T^{t}\}_{t\in\R})$ is minimal, $A:=\{t\in \R: T^{t}x_0\in U\}$ is discrete syndetic. Thus ${\rm Dsynd\textendash sup} f(T^{t}x_0)\geq M-\epsilon$. Since $\epsilon$ is arbitrary, we have  ${\rm Dsynd\textendash sup} f(T^{t}x_0)\geq M$. Clearly,  ${\rm Dsynd\textendash sup} f(T^{t}x_0)\leq M$. Thus we conclude that  \[{\rm Dsynd\textendash sup}f(T^tx_0) =\sup_{t\in \R} f(T^{t}x_0).\]
 This completes the proof.
  \end{proof}

 \begin{proof}[Proof of Theorem \ref{nilfunction}]
 Let $f\in L^{\infty}(\mu)$. Without loss of generality, we may assume that $\|f\|_{\infty}\leq 1$.

 Let $g=\mathbb{E}(f\mid \mathcal{Z}_{k})$. Then it follows from Lemma \ref{transfer to factor} that $\{I_{f}(k, t)-I_{g}(k,t)\}$ converges to zero in uniform density. Thus it suffices show the theorem for the function $g$ substituted for $f$.

$Z_{k}(X)$ is an inverse limit of a sequence of ergodic $k$-step nilsflows. Let $r$ be a positive integer. There is a factor $X'$ of $Z_{k}(X)$, which is a $k$-step nilflow, such that
\[\| g-\mathbb{E}(g\mid \mathcal{X}')\|_{L^{1}(\mu)}\leq \frac{1}{(k+1)r}.\]
Let $g'=\mathbb{E}(g\mid \mathcal{X}')$. Then for each $t\in \R$, $|I_{f}(k,t)-I_{g'}(k,t)|\leq \frac{1}{r}$.
Just as the explanations in \cite[Section 4.2]{P11}, writing $X'=G/\Gamma$ we may assume that $G$ is connected and simply connected.  Then by Proposition \ref{I-J} and Proposition \ref{nilfunction for f} , $\{I_{g'}(k,t)\}$ is decomposed as a sum of a $k$-step nilfunction and a function tending to zero in uniform density.  Thus we have
\[ I_{f}(k,t)=a_r(t)+b_r(t)+c_r(t),\]
where $|a_r(t)|\leq \frac{1}{r}$ for each $t\in \R$, ${\rm UD \textendash Lim}~ b_r(t)=0$ and $c_{r}(t)$ is a basic $k$-step nilfunction.

Now for $r\neq r'$, we have
\[ c_r(t)-c_{r'}(t)=\big( a_r(t)-a_{r'}(t)\big) +\big( b_r(t)-b_{r'}(t)\big).\]
We also have ${\rm UD\textendash Lim}\big( b_r(t)-b_{r'}(t)\big)=0 $ and $\sup_{t\in\R}| a_r(t)-a_{r'}(t)|\leq \frac{1}{r}+\frac{1}{r'}$. Thus by Lemma \ref{UD-lim vs Dsynd},
\[ {\rm Dsynd\textendash sup} |c_{r}(t)-c_{r'}(t)|\leq \frac{1}{r}+\frac{1}{r'}.\]
 Since $\{c_r(t)-c_{r'}(t)\}$ is a nilfunction, it follows from Lemma \ref{synd sup for nil} that $\sup_{t\in\R} |c_{r}(t)-c_{r'}(t)|\leq \frac{1}{r}+\frac{1}{r'}$.  Therefore, $\{c_{r}(t)\}$ is a Cauchy sequence in $L^{\infty}(\R)$ and it converges uniformly to some function $c(t)$. This function is a $k$-step nilfunction. It is easy to verify that $\{ I_{f}(k, t)-c(t)\}$ converges to zero in uniform density.
 \end{proof}

\section{Proof of Theorem \ref{char of RP^{d}}}
\begin{defn}
A subset $A\subset \R$ is a $Nil_{d} ~Bohr_{0}$-{\it set} if there exist a $d$-step nilflow $(X,\mu, \{T^{t}\}_{t\in\R})$, $x_0\in X$ and a neighborhood $U$ of $x_0$ such that
\[R(x_0, U):=\{ t\in \R: T^{t}x_0\in U\}\]
is contained in $A$.
\end{defn}
We use $\mathcal{F}_{d,0}(\R)$ to denote collection of all ${\rm Nil}_{d} ~{\rm Bohr}_{0}$-sets and  $\mathcal{F}_{d,0}^{*}(\R)$ to denote the dual of  $\mathcal{F}_{d,0}(\R)$, that is
\[  \mathcal{F}_{d,0}^{*}(\R)=\{ A\subset \R: \forall B\in \mathcal{F}_{d,0}(\R), A\cap B\neq \emptyset\}.\]

We follow the line given in \cite{HSY16} to show Theorem \ref{char of RP^{d}}.
\begin{thm} \label{equiv char of RPd}
Let $(X, \{T^{t}\}_{t\in\mathbb{R}})$ be a minimal $\R$-flow and $x,y\in X$. The following assertions are equivalent for $d\in\N$:
\begin{enumerate}
\item[(1)] $(x,y)\in {\bf RP}^{[d]}(X)$.
\item[(2)] $R(x, U)\in \mathcal{F}_{Poi_{d}}(\R)$ for every neighborhood $U$ of $y$.
\item[(3)] $R(x, U)\in \mathcal{F}_{Bir_{d}}(\R)$ for every neighborhood $U$ of $y$.
\item[(4)] $R(x, U)\in \mathcal{F}_{d,0}^{*}(\R)$ for every neighborhood $U$ of $y$.
\end{enumerate}
\end{thm}

In the same way as \cite[Corollary D]{HSY16}, one has
\[ \mathcal{F}_{Poi_{d}}(\R)\subset \mathcal{F}_{Bir_{d}}(\R)\subset \mathcal{F}_{d,0}^{*}(\R).\]
Then, it is clear that (2) $\Longrightarrow$ (3) $\Longrightarrow$ (4). Furthermore, (1) can be derived from (3) by the discrete version (\cite[Theorem 7.2.7]{HSY16}). 

In the same way as  \cite[Theorem 7.2.5]{HSY16}, it is easy to show that (4) $\Longrightarrow$ (1).Indeed, since the topology on $\R$ does not affect the regionally proximal relation of higher order, we can obtain the continuous version of  \cite[Theorem 7.2.5]{HSY16} by endowing $\R$ with discrete topology. Furthermore, it follows from the proof there that $R(x,U)\cap F$ is a discrete syndetic set for each $F\in \mathcal{F}_{d,0}(\R)$.

Next we show that (1) implies (2).

We need the following special case of continuous polynomial ergodic Szemer\'{e}di theorem shown in \cite{BLM12}.
\begin{lem}\cite[Theorem 8.15]{BLM12}\label{poly sez}
Let $(X, \mu, \{T^{t}\}_{t\in \R})$ be a measure preserving flow and let $p_i: \R\rightarrow \R, i=1,2,\ldots,k$ be polynomials with $p_i(0)=0$ for each $i$. If $A$ is a measurable subset of $X$ with $\mu(A)>0$, then for any F{\o}lner sequence $(\Phi_{n})$ in $\R$,
\[\liminf_{n\rightarrow \infty} \frac{1}{|\Phi_{n}|} \int_{\Phi_{n}} \mu\left(A\cap T^{-p_1(t)}A\cap\cdots\cap T^{-p_{k}(t)}A \right)dt>0,\]
where $|\Phi_{n}|$ is the Lebesgue measure of $\Phi_{n}$.
\end{lem}

Let $(X, \mu,\{T^{t}\}_{t\in\R})$ be an ergodic flow and  let $A$ be a measurable set in $X$ with $\mu(A)>0$. If we take $f={\bf 1}_{A}$ in Theorem \ref{nilfunction}, then we can write $I_{f}(k, t)=F_{k}(t)+N(t)$ where $F_{k}(t)$ is a $k$-step nilfunction and $N(t)$ tends to zero in uniform density. Suppose that we write $F_{k}(t)=\phi(T^{t}y_0)$, where $(Y, \nu,\{T^{t}\}_{t\in\R})$ is a $k$-step nilflow , $\phi$ is a continuous function on $Y$ and $y_0\in Y$. 

\noindent{\bf Claim}. We may assume that $\phi(y_0)>0$.

\begin{proof}[Proof of Claim]
Since $\{I_{f}(k,t)-\phi(T^{t}y_0)\}$ tends to zero in uniform density and by Lemma \ref{poly sez}, one has 
\[ \lim_{\rho\rightarrow \infty}\inf_{\sigma\in\R}\frac{1}{\rho}\int_{\sigma}^{\sigma+\rho} I_{f}(k,t)dt >0,\]
then we have
\[ \lim_{\rho\rightarrow \infty}\inf_{\sigma\in\R}\frac{1}{\rho}\int_{\sigma}^{\sigma+\rho} \phi(T^{t}y_0)dt >0.\]
In particular, there is some $y_1\in Y$ with $\phi(y_1)>0$. Recall that $(Y, \nu,\{T^{t}\}_{t\in\R})$ is uniquely ergodic, it follows from the Oxtoby Theorem that for any $y\in Y$,
\[ \lim_{\rho\rightarrow \infty}\inf_{\sigma\in\R}\frac{1}{\rho}\int_{\sigma}^{\sigma+\rho} \phi(T^{t}y_0)dt=\int_{Y}\phi(y)d\nu(y).\]
As a consequence, $\{\phi(T^{t}y_0)-\phi(T^{t}y_1)\}$ tends to zero in uniform density. By replacing $\phi(T^{t}y_0)$ with $\phi(T^{t}y_1)$, the assertion is followed.
\end{proof}

Let $0<\epsilon<\phi(y_0)$. Then we conclude that
\begin{lem}\label{almost nil Bohr}
Let $(X, \mu,\{T^{t}\}_{t\in\R})$ be an ergodic flow and  let $A$ be a measurable set in $X$ with $\mu(A)>0$. Then for any distinct $\a_1,\a_2,\ldots,\a_k$, there is some $\delta>0$ such that
\[I=\{ t\in \R: \mu(A\cap T^{-\a_1 t}A\cap T^{-\a_2 t}A\cap\cdots\cap T^{-a_k t}A)>\delta\}\]
is almost ${\rm Nil}_{k}~{\rm Bohr}_{0}$-set, i.e., there is some set $M\subset \R$ such that $I\Delta M$ is a ${\rm Nil}_{k}~{\rm Bohr}_{0}$-set and the characteristic function of $M$ tends to zero in uniform density.
\end{lem}
\begin{rem}\label{choice of delta}
We may choose $\delta= \frac{1}{3}\lim_{\rho\rightarrow \infty}\inf_{\sigma\in\R}\frac{1}{\rho}\int_{\sigma}^{\sigma+\rho} I_{f}(k,t)dt $, since we may assume that $\phi(y_0)\geq \lim_{\rho\rightarrow \infty}\inf_{\sigma\in\R}\frac{1}{\rho}\int_{\sigma}^{\sigma+\rho} I_{f}(k,t)dt $. Then
\[\{t\in\R: T^{t}y_0>2\delta\}=\{t\in \R: \phi(T^{t}y_0)>2\delta\}\]
is a ${\rm Nil}_{k}~{\rm Bohr}_{0}$-set and the upper Banach density of $\{t\in \R: |M(t)|\geq \delta\}$ is zero.
\end{rem}

 \begin{proof}[Proof of (1) $\Longrightarrow$ (2)]
Suppose that $(x,y)\in {\bf RP}^{[d]}$ and fix a neighborhood $U$ of $y$. We need to show $R(x,U)\in \mathcal{F}_{Poi_d}(\R)$. Let $\a_1,\a_2,\ldots,\a_d$ be distinct real numbers.

Let $(Y,\mathcal{Y}, \mu, \{T^{t}\}_{t\in \R})$ be a measure preserving flow and $A\in\mathcal{Y}$ with $\mu(A)>0$. Let $\mu=\int_{\Omega}\mu_{\omega}dm(\omega)$ be an ergodic decomposition of $\mu$. Then there is some $\Omega'\subset \Omega$ with $m(\Omega')>0$ such that for each $\omega\in\Omega'$, $\mu_{\omega}(A)>0$. For each $\omega\in \Omega'$, set $\delta_{\omega}>0$ as in Remark \ref{choice of delta}. Then $\delta_{\omega}$ is measurable with respect to $\omega$. Further, set
\[ F_{\omega}:=\{t\in\R: \mu_{\omega}(A\cap T^{-\a_1 t}A\cap T^{-\a_2 t}A\cap\cdots\cap T^{-\a_d t}A)>\delta_{\omega}\}.\]
By Lemma \ref{almost nil Bohr}, there is some $M_{\omega}\subset \R$ with zero upper Banach density such that $B_{\omega}=F_{\omega}\Delta M_{\omega}$ is a ${\rm Nil}_{d}$ ${\rm Bohr}_{0}$-set. Then by the remark below Proposition \ref{equiv char of RPd}, $R(x,U)\cap (F_{\omega}\Delta M_{\omega})$ is discrete syndetic. Then there is some $BD_{*} (R(x,U)\cap F_{\omega})>0$ since $BD^{*}(M_{\omega})=0$. Let $\Omega''\subset \Omega'$ with $m(\Omega'')>0$ be such that for each $\omega\in\Omega''$, $\delta_{\omega}\geq \delta>0$ for some $\delta>0$.
Now for each $\sigma\in\R$ and $\rho>0$, one has
\begin{align*}
&\liminf_{\rho\rightarrow\infty}\frac{1}{\rho}\int_{0}^{\rho} \int_{Y} {\bf 1}_{A\cap T^{-\a_1 t}A\cap \cdots \cap T^{-\a_d t}A}(z)\cdot {\bf 1}_{R(x,U)}(t)d\mu(z) dt \\
=& \liminf_{\rho\rightarrow\infty}\frac{1}{\rho}\int_{0}^{\rho} \int_{Y} \int_{\Omega}{\bf 1}_{A\cap T^{-\a_1 t}A\cap \cdots \cap T^{-\a_d t}A}(z)\cdot {\bf 1}_{R(x,U)}(t) d\mu_{\omega}(z)dm(\omega) dt \\
\geq&  \int_{\Omega}\liminf_{\rho\rightarrow\infty}\frac{1}{\rho}\int_{0}^{\rho} \int_{Y}{\bf 1}_{A\cap T^{-\a_1 t}A\cap \cdots \cap T^{-\a_d t}A}(z)\cdot {\bf 1}_{R(x,U)}(t) d\mu_{\omega}(z)dtdm(\omega)  \\
\geq&  \int_{\Omega'}\liminf_{\rho\rightarrow\infty}\frac{1}{\rho}\int_{0}^{\rho} \int_{Y}{\bf 1}_{A\cap T^{-\a_1 t}A\cap \cdots \cap T^{-\a_d t}A}(z)\cdot {\bf 1}_{R(x,U)}(t) d\mu_{\omega}(z)dtdm(\omega)   \\
=&  \int_{\Omega'}\liminf_{\rho\rightarrow\infty}\frac{1}{\rho}\int_{0}^{\rho} \mu_{\omega}(A\cap T^{-\a_1 t}A\cap \cdots \cap T^{-\a_d t}A)\cdot {\bf 1}_{R(x,U)}(t) dtdm(\omega)  \\
\geq& \delta  \int_{\Omega''}\liminf_{\rho\rightarrow\infty}\frac{1}{\rho}\int_{0}^{\rho} {\bf 1}_{R(x,U)\cap F_{\omega}}(t)dtdm(\omega)\\
\geq& \delta  \int_{\Omega''}BD_{*}(R(x,U)\cap F_{\omega})dm(\omega)>0.
\end{align*}
This implies that there is some $t>0$ such that
\[\mu(A\cap T^{-\a_1 t}A\cap T^{-\a_2 t}A\cap\cdots\cap T^{-\a_d t}A)>0.\]
Therefore, $R(x,U)\in\mathcal{F}_{Poi_{d}}(\R)$.
\end{proof}

\end{document}